\documentclass[11pt, a4paper]{article}
\usepackage[margin=.8in]{geometry}

\RequirePackage{amsthm,amsmath,amsfonts,amssymb,mathtools}
\usepackage[square,comma,sort&compress,numbers]{natbib}

\renewcommand{\citep}{\cite}

\renewcommand{\citet}{\cite}

\renewcommand{\citealp}{\cite}
\RequirePackage[colorlinks,citecolor=black,urlcolor=black,linkcolor=black]{hyperref}
\usepackage{enumitem}
\usepackage{bbm}
\usepackage{comment}
\usepackage{csquotes}
\usepackage{enumitem} 
\usepackage{stackengine}
\usepackage{booktabs}
\usepackage{makecell}
\usepackage{tablefootnote}
\usepackage{algorithm}
\usepackage{algorithmicx}
\usepackage[noend]{algpseudocode}
\algblockdefx{MRepeat}{EndRepeat}[1]{\algorithmicrepeat\ #1 \textbf{times}}{}
\algnotext{EndRepeat}
\usepackage{multirow}
\usepackage{xcolor}
\usepackage{todonotes}
\usepackage{caption}
\usepackage{siunitx}

\usepackage{setspace}
\allowdisplaybreaks

\usepackage[capitalize,noabbrev]{cleveref}

\theoremstyle{plain}
\newtheorem{theorem}{Theorem}[section]
\newtheorem{lemma}[theorem]{Lemma}
\newtheorem{corollary}[theorem]{Corollary}

\newtheorem{definition}[theorem]{Definition}

\newcommand\xrowht[2][0]{\addstackgap[.5\dimexpr#2\relax]{\vphantom{#1}}}

\newtheorem{remark}[theorem]{Remark}

\newtheorem{proposition}[theorem]{Proposition}

\newtheorem*{assumption}{Assumption}

\newtheorem*{AoK}{Assumptions on kernels}

\renewcommand{\P}{\mathbb{P}}
\newcommand{\N}{\mathcal{N}}
\newcommand{\E}{\mathbb{E}}
\newcommand{\R}{\mathbb{R}}
\newcommand{\V}{\mathbb{V}}
\newcommand{\1}{\mathbbm{1}}

\newcommand\norm[1]{\left\lVert#1\right\rVert}
\newcommand{\B}{\mathcal{B}}
\newcommand{\normm}[1]{\left\lVert#1\right\rVert_{2}}

\newcommand{\normf}[1]{\left\lVert#1\right\rVert_{F}}
\newcommand{\normop}[1]{\left\lVert#1\right\rVert_{\rm op}}

\DeclarePairedDelimiter\ceil{\lceil}{\rceil}
\DeclarePairedDelimiter\Ceil{\bigg\lceil}{\bigg\rceil}
\DeclarePairedDelimiter\floor{\lfloor}{\rfloor}

\DeclarePairedDelimiter\abs{\lvert}{\rvert}

\DeclareMathOperator{\sign}{sign}

\graphicspath{{./Figures/}}

\title{Adaptive monotonicity testing in sublinear time}

\author{Housen Li \and Zhi Liu \and Axel Munk}

\date{Institute for Mathematical Stochastics, University of G\"ottingen, Germany}

\begin{document}
\maketitle

\begin{abstract}
Modern large-scale data analysis increasingly faces the challenge of achieving computational efficiency as well as statistical accuracy, as classical statistically efficient methods often fall short in the first regard. In the context of testing monotonicity of a regression function, we propose FOMT (Fast and Optimal Monotonicity Test), a novel methodology  tailored to meet these dual demands. FOMT employs a sparse collection of \emph{local} tests, strategically generated at random, to detect violations of monotonicity scattered throughout the domain of the regression function. This sparsity enables significant computational efficiency, achieving sublinear runtime  in most cases, and quasilinear runtime (i.e.~linear up to a log factor) in the worst case. In contrast, existing statistically optimal tests typically require at least quadratic runtime. FOMT's statistical accuracy is achieved through the precise calibration of these local tests and their effective combination, ensuring both sensitivity to violations and control over false positives. More precisely, we show that FOMT separates the null and alternative hypotheses at minimax optimal rates over Hölder function classes of smoothness order in $(0,2]$. Further, when the smoothness is unknown, we introduce an adaptive version of FOMT, based on a modified Lepskii principle, which attains statistical optimality and meanwhile maintains the same computational complexity as if the intrinsic smoothness were known. Extensive simulations confirm the competitiveness and effectiveness of both FOMT and its adaptive variant.
\end{abstract}

\noindent%

\section{Introduction}
The monotonic relationship between a response variable and explanatory variables is often a reasonable assumption and exploited across various fields of science and technology. For instance, in medical studies, the white blood cell count and the DNA index are typically modeled as monotonic functions contingent upon the expression of leukemia antigens (cf.\ \citealp{Schell1997The}). In economics, monotonic regression is employed to analyze the relationship between customer satisfaction scores and the number of services utilized in customer support systems \citep{Graves1994A}. Besides important scientific insights, monotonicity can be exploited for a variety of statistical tasks, such as estimation or prediction. For instance, Obozinski, {\it et al.}~\citet{obozinski2008consistent} demonstrated that monotonic regression typically enhances prediction consistency without compromising precision. As an example, in protein function prediction, one seeks for probabilistic predictions for hierarchical protein function annotations, ensuring that these predictions align with the monotonic structure dictated by the gene ontology hierarchy. 

Popular techniques for estimation of such an monotone regression function include (nonparametric) maximum likelihood (for a comprehensive treatment, see \citealp{robertson1988orderb,groeneboom2014nonparametric}), monotone regression splines \citealp{ramsay1998estimating,meyer2008inference} and Bayesian estimation \citealp{lin2014bayesian,shively2009bayesian,okano2024locally}. However, an improper application of the monotonicity assumption may result in biased estimates and erroneous conclusions, as illustrated by \citet{swanson2015definition} in the context of instrumental variable studies. Therefore, it is highly advisable to conduct monotonicity checks or tests, prior to any data analysis based on a monotonicity assumption of the signal.

For the purpose of our theoretical analysis, in this paper, we focus on a nonparametric regression model with additive Gaussian noise, where the observations are given~by 
\begin{equation}
    \label{model}
    Y_i = f(x_i) + \varepsilon_i,\quad i \in [n] := \{1,\dots,n\}
\end{equation}
for equidistant sampling points $x_i \equiv i/n$. Here $f:[0,1]\to \mathbb{R}$ is an \emph{unknown} function and the random errors $\varepsilon_i$ are \textit{i.i.d.}~Gaussian distributed with mean zero and \emph{known} variance $\sigma^2$, for simplicity. We stress, however, that all our results can be transferred to unknown variance, see Remark~\ref{r: unknown variance}. We assume that the regression function (or signal) \(f\) is Hölder smooth of order~\(\beta\), i.e., 
\begin{enumerate}[wide]
    \item [{\rm (M1)}] \label{M1} The signal $f\in \Sigma(\beta,L)$ with \emph{known} $\beta \in (0,2]$ and \emph{known} $L>0$.
\end{enumerate}
The case of \emph{unknown} $\beta\in (0,2]$ is considered in Section~\ref{S: Adaptivity}. Here $\Sigma(\beta,L )$ denotes the  Hölder class on \([0,1]\) of (smoothness) order \(\beta>0\) and radius \(L >0\), which is defined as
\begin{equation*}    
    \Sigma(\beta,L) \coloneqq \left\{f:[0,1]\to \mathbb{R} \;\middle|\; \bigl|f^{(\ceil{\beta}-1)}(x)-f^{(\ceil{\beta}-1)}(y)\bigr|\leq L \abs{x-y}^{\beta+1-\ceil{\beta}} \text{ for all } x,y\right\}
\end{equation*}
where $\ceil{\beta} \coloneqq \min\{x\in \mathbb{Z}\,|\,x\geq\beta\}$. For example, the class $\Sigma(1,L)$ (i.e.\ $\beta= 1$) consists of all {Lipschitz} functions with the Lipschitz constant $L$.  
Our goal is to test whether the function $f$ is monotone increasing (i.e.\ isotonic) or not, using the observations $\{(x_i, Y_i)\}_{i = 1}^n$ from \eqref{model}. Formally, we consider the null hypothesis
\begin{equation}\label{eq:null}
    H: f \in H \coloneqq \Sigma(\beta,L)\cap \mathcal{M},
\end{equation}
where $\mathcal{M}$ denotes the set of monotone increasing functions $f:[0,1]\to \mathbb{R}$.  

\subsection{Related work}\label{ss:literature}
Early approaches to test the monotonicity of regression function include kernel-based methods \citep{schlee1982nonparametric,bowman1998testing}, local least-squares-based tests \citep{hall2000testing} and robust approaches based on sign-type statistics \citep{gijbels2000tests}. Somewhat similar in spirit is the test by \citet{ghosal2000testing}, which is based on a $U$-statistic that can be seen as locally estimating the degree of concordance as in Kendall's tau statistic. The separation rates of this test are established over Hölder class $\Sigma(\beta,L)$ when the smoothness order $\beta$ is \emph{known}, which aligns with the minimax lower bounds (see e.g.\ \citealp{dumbgen2001multiscale}). An alternative minimax optimal approach, via estimating the distance to the set of monotone functions, is discussed in \citet{JuNe02}. 

A next step beyond minimax has been taken by \citet{dumbgen2001multiscale}, who proposed testing procedures that aggregate test statistics across a collection of bandwidths and locations under a Gaussian white noise model, which can be shown to be \emph{adaptively} minimax optimal for Hölder smooth functions of order $\beta = 1$ or $2$.  {See \citet{SHMD13} for an extension to convolution models.} Baraud, {\it et al.}~\citet{baraud2005testing} introduced a test based on the differences of local averages of observations from consecutive intervals of various sizes, establishing adaptive minimax optimality over the Hölder class for $\beta \in (0,1]$ in terms of $L^{\infty}$-norm. Akakpo, {\it et al.}~\citet{akakpo2014testing} developed a multiscale test (as an extension of \citealp{durot2003kolmogorov}) by combing many local least concave majorants, and established a test that achieves adaptive minimax optimality over Hölder classes with $\beta \in(0,2]$. Additionally, Chetverikov~\citet{chetverikov2019testing} proposed an adaptive version of the Kendall's tau statistic by maximizing the statistic over various bandwidths (as an extension of \citealp{hall2000testing} and \citealp{ghosal2000testing}), and proved adaptive minimax optimality over Hölder classes with $\beta \in(1,2]$.

More recently, monotonicity testing was also investigated under the Bayesian framework. For example, Scott, {\it et al.}~\citet{scott2015nonparametric} introduced two testing procedures utilizing constrained smoothing splines with a hierarchical stochastic process prior, and regression splines with a prior over the regression coefficients. Salomond~\citet{salomond2018testing} put a posterior distribution on the largest absolute discrepancy between the parameter and the null model. The resulting test is shown to attain asymptotic frequentist optimality, being adaptively minimax optimal (up to a log-factor) for Hölder smooth alternatives of order $\beta\in (0,1]$, see also \citet{chakraborty2021convergence} for discrepancies between null and alternatives measured in the Hellinger distance.

\subsection{Our contribution}
Despite the statistical guarantees --- such as consistency and (adaptive) minimax optimality --- provided by many of the aforementioned testing procedures, there remains a lack of attention on computational efficiency. In the modern era of \emph{big data} (i.e.\ of enormous sample sizes), computational scalability has become  a primary criterion for statistical procedures (cf.\ \citealp{fan2014challenges}). From this perspective, the existing (adaptively) minimax optimal methods (see Section~\ref{ss:literature} for details) ask for at least a \emph{quadratic} runtime in terms of sample sizes, significantly limiting their applicability in large-scale data scenarios.  To overcome this computational burden, we propose in this paper a novel monotonicity testing procedure under the nonparametric regression model in \eqref{model}, which we call the \emph{Fast and Optimal Monotonicity Test} (FOMT; cf.\ Algorithm~\ref{alg:MC}). As its name suggests,  FOMT attains not only (adaptive) minimax optimality in separation of null and alternatives for $\beta \in (0,2]$ but also computational efficiency in the sense that it has a computational complexity \emph{sublinear} in sample sizes in most situations (specified by mild conditions).

FOMT is built upon a collection of \emph{local} tests of {local} hypotheses $H_{i,j}: f(x_i)-f(x_j)\leq 0$ for indices $i,j\in [n] \coloneqq\{1,\ldots, n\}$ and $i < j$. It rejects the null hypothesis that $f$ is monotone increasing if at least one  $H_{i,j}$ is rejected. Specifically, we reject $H_{i,j}$ if $\hat{f}_n(x_j)- \hat{f}_n(x_i)$ exceeds a properly specified critical value, where $\hat{f}_n$ is a {local polynomial estimator} of $f$. The success of FOMT in achieving both statistical and computational efficiencies relies on two key aspects. First,  $\hat{f}_n(x_i)$ and $\hat{f}_n(x_j)$ exhibit strong correlation when $x_i$ and $x_j$ are close (Theorem~\ref{Theorem: Bounds of variances} in the appendix), apart from the well-known minimax optimality of local polynomial estimators. Second,  a sparse and random selection of local tests ensures that, with high probability, every local property of $f$ is examined with minimal computational effort. This selection strategy is inspired by the \emph{spot-checkers} from the computer science literature, originally introduced by \citet{ergun2000spot}, as a fast method for verifying correctness of computer programs. 

\begin{figure}
\centering
\includegraphics[width=0.95\textwidth]{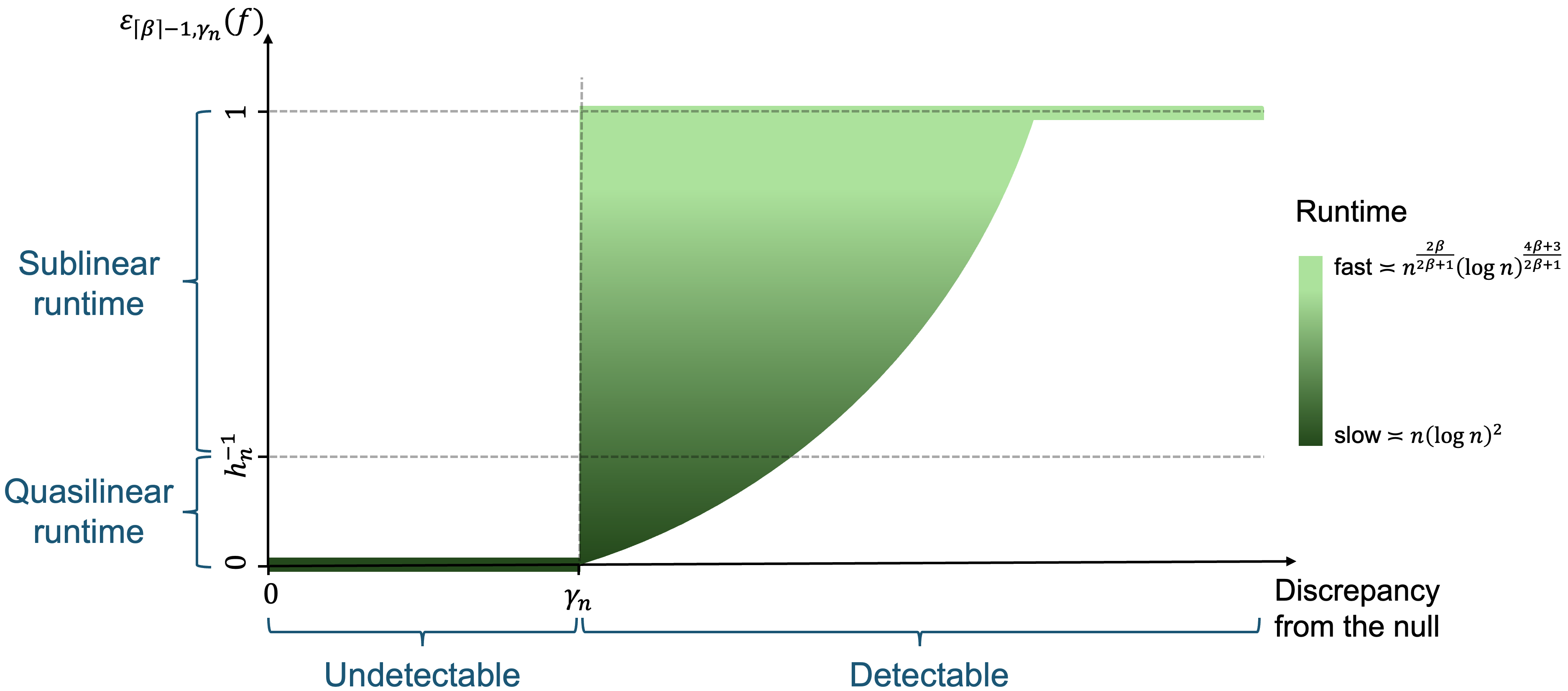}
\caption{Phase diagram of FOMT in computational and statistical efficiency for the Hölder class $\Sigma(\beta, L)$ with $\beta \in (0,2]$. Here 
$h_n \asymp (\log(n)/n)^{1/(2\beta+1)}$ is the optimal bandwidth of the underlying local polynomial estimators, and $\gamma_n \asymp (\log(n)/n)^{(\beta + 1-\ceil{\beta})/(2\beta+1)}$ the minimax estimation error in $L^{\infty}$-norm. The shaded (green) region shows all possible relations between the $\gamma$-exceedance fraction $\varepsilon_{\ceil{\beta}-1,\gamma_n}(f)$ and the discrepancy from the null hypothesis (i.e.\ monotone increasing functions). The discrepancy is measured by $\max \{f(a) - f(b)\;|\; 0 \le a \le b\le 1\}$ if $\beta \in (0,1]$, and by $\max(\{-f'(x)\;|\;0\le x \le 1\} \cup \{0\})$ if $\beta \in (1,2]$. This phase diagram not only delineates the phase transition in statistical detectability of alternatives, but also highlights computational complexity (indicated by varying degrees of darkness) interpolating between $O\bigl(n^{\frac{2\beta}{2\beta+1}} \left(\log n\right)^{\frac{4\beta+3}{2\beta+1}}\bigr)$ and $O\bigl(n(\log n)^2\bigr)$. See Section~\ref{SS: CA of Phi} for further details.} 
\label{F:Phase}
\end{figure}

Our analysis of FOMT suggests that the \emph{scale} of deviations from monotonicity affects computational complexity. To quantify this scale, we introduce the \emph{$\gamma$-exceedance fraction} of $f$, which measures the proportion of its domain on which the signal $f$ is ``$\gamma$-apart'' from any monotone functions on $[0,1]$.   
\begin{definition}[$\gamma$-exceedance fraction]
\label{e: exceedance 0}
Let $\mathcal{M}$ be the set of monotone increasing functions on \([0,1]\), $\lambda$ the Lebesgue measure and $\gamma>0$. 
        \begin{enumerate}[label=(\roman*), wide]
        \item Let $f\in  C\bigl([0,1]\bigr)$ the set of continuous functions on \([0,1]\). The \emph{$\gamma$-exceedance fraction of order zero} of $f$ to $\mathcal{M}$ is defined~as
        \begin{subequations}
        \begin{equation}
            \label{e:epsilon_f}
            \varepsilon_{0,\gamma}(f) = \underset{g\in \mathcal{M}}{\min}\lambda\left(\left\{x\in[0,1] \;| \;\abs{f(x)-g(x)} >\gamma \right\}\right).
        \end{equation}
        \item Let $f\in  C^1\bigl([0,1]\bigr)$ the set of continuously differentiable functions on \([0,1]\). The \emph{$\gamma$-exceedance fraction of order one} of $f$ to $\mathcal{M}$ is defined as
        \begin{equation}
            \label{e:epsilon_f_1}
            \varepsilon_{1,\gamma}(f) = \lambda\left(\{x\in [0,1]\;|\; f'(x)\leq -{\gamma}\}\right).
        \end{equation}
        \end{subequations}
    \end{enumerate}
\end{definition}
Basic properties of \(\gamma\)-exceedance fraction (incl.\ the existence of minimizers in \eqref{e:epsilon_f}) are provided in Section~\ref{s:gamma_exceedance} in the appendix. The \(\gamma\)-exceedance fraction allows to display computational and statistical efficiency in a phase diagram of FOMT in Figure~\ref{F:Phase}. More precisely, we show that FOMT has the following three favorable properties:
\begin{enumerate}[label=\roman*.]
    \item \emph{Minimax optimality.} Under the nonparametric regression model in \eqref{model}, FOMT achieves minimax optimal separation rates between monotone functions and alternatives given by Hölder classes of smoothness order $\beta \in (0,2]$, where $\beta$ is assumed to be \emph{known}. 
    
    \item \emph{(Sub)linear runtime.} The FOMT has the computational complexity of order
    $$    
    \min\left\{\bigl(\varepsilon_{\ceil{\beta}-1,\gamma_n}(f)\bigr)^{-1}{n^{\frac{2\beta}{2\beta + 1}} (\log n)^{\frac{4\beta + 3}{2\beta + 1}}},\;n (\log n)^2\right\} 
    $$
    with $\gamma_n \asymp (\log(n)/n)^{(\beta +1 -\ceil{\beta})/(2\beta+1)}$.
Under detectable alternatives, the first term in the above displayed equation determines the runtime of FOMT, as $\varepsilon_{k,\gamma_n}(f) \gtrsim h_n \asymp (\log(n)/n)^{1/(2\beta+1)}$ (cf.\ Corollary~\ref{p:>=hn} in the appendix). In particular, FOMT has a sublinear runtime of order $n^{2\beta/(2\beta + 1)}(\log n)^{(4\beta +3)/(2\beta + 1)}$ when the detectable violations occur on a subset of $[0, 1]$ with a non-vanishing (Lebesgue) measure, see Figure~\ref{F:Phase}.  It is worth noting that FOMT is significantly faster than existing minimax optimal testing procedures, which at least require a runtime of $O(n^2)$, see Table~\ref{table: complexity}.

\item \emph{Computational and statistical adaptation in concert.} When the smoothness parameter is \emph{unknown}, we introduce a variant of Lepskii principle CALM (\underline{C}omputationally \underline{A}daptive \underline{L}epskii \underline{M}ethod), with adjustments towards computational efficiency, for FOMT to automatically select the tuning parameter (i.e.\ the bandwidth)  from the data. This approach allows FOMT to achieve adaptive minimax optimality over H\"older classes of smoothness order $\beta\in(0,2]$, and meanwhile its computational complexity remains almost the same as if $\beta$ were known, under mild assumption on $f$ in form of \emph{self-similarity} \citep{gine2010confidence}. This assumption primarily {ensures} that the \emph{intrinsic} smoothness of $f$ is exactly~$\beta$. While this particular form of CALM has appeared in the literature (see e.g.\ \citealp{reiss2009pointwise,REISS2012} in the context of quantile regression), its computational advantage has not been previously recognized.
\end{enumerate}

\begin{table}
\caption{Comparison of FOMT with other existing procedures. Here DS stands for the procedure in  \citet{dumbgen2001multiscale}, GSV for \citet{ghosal2000testing}, C for \citet{chetverikov2019testing} and BHL for \citet{baraud2005testing}. For the proposed FOMT, the factor $\varepsilon_{\ceil{\beta}-1,\gamma_n}(f)$ satisfies $\bigl(\log(n)/n \bigr)^{1/(2\beta+1)} \lesssim \varepsilon_{\ceil{\beta}-1,\gamma_n}(f) \lesssim 1$. For DS, C and BHL, the factor $R$ corresponds to the number of repetitions in Monte--Carlo or bootstrap procedures. We consider the statistical optimality in terms of minimax separation rates over Hölder smooth functions of order $\beta$. In this regard, the established statistical guarantees for different methods are summarized in the last column.\label{table: complexity}}
\centering{\footnotesize
\begin{tabular}{| l| l| l| l|l|}
\toprule
\hline \xrowht[()]{10pt}
\multirow{2}{*}{\textbf{Methods}} & \multicolumn{3}{c|}{\textbf{Computational complexity} (saving multiplicative constants)}  & \multirow{2}{*}{\makecell{\textbf{Minimax}\\\textbf{optimality}}}\\[5pt]
\cline{2-4}\xrowht[()]{10pt}
& \textbf{General case}& \textbf{Best case} &\textbf{Worst case}& \\
    \hline \xrowht[()]{10pt}
    FOMT &${\varepsilon^{-1}_{\ceil{\beta}-1,\gamma_n}(f) \cdot n^{\frac{2\beta}{2\beta+1}} \left(\log n\right)^{\frac{4\beta+3}{2\beta+1}}}$  & $n^{\frac{2\beta}{2\beta+1}} \left(\log n\right)^{\frac{4\beta+3}{2\beta+1}}$ & $n\left(\log n\right)^2$ & \text{$\beta\in (0,2]$}\\[5pt]
    \hline \xrowht[()]{10pt}
    A-FOMT\tablefootnote{For achieving computational adaptability, additional conditions are required (cf.\ Theorem~\ref{t:MC Lepskii compleixity} in Section~\ref{S: Adaptivity}).} &${ \varepsilon^{-1}_{\ceil{\beta}-1,\gamma_n}(f) \cdot n^{\frac{4}{5}} \left(\log n\right)^{\frac{11}{5}}}$  & $n^{\frac{2\beta}{2\beta+1}} \left(\log n\right)^{\frac{4\beta+3}{2\beta+1}}$ & $n^{\frac{9}{5}}\left(\log n\right)^{\frac{6}{5}}$ & \text{$\beta\in (0,2]$}\\[5pt]
    \hline \xrowht[()]{10pt}
    DS & $Rn^2$ & $Rn^2$ & $Rn^2$ & \text{$\beta=1$ or $2$}\\[5pt]
    \hline \xrowht[()]{10pt}
    GSV & $n^3h_n^2$\;\; with $n^{-\frac{1}{3}}\ll h_n \ll 1$ & $n^{\frac{7}{3}}$ & $n^3$ &\text{$\beta\in (1,2]$}\\[5pt]
    \hline \xrowht[()]{10pt}
    C & $Rn^3$&$Rn^3$ &$Rn^3$ & \text{$\beta\in (1,2]$}\\[5pt]
    \hline \xrowht[()]{10pt}
    BHL & $R(l_n^3\vee n l_n)$\;\; with $1\leq l_n \leq n$ & $Rn$ & $Rn^3$ &\text{$\beta\in (0,1]$}\\[5pt]
    \hline
\bottomrule
\end{tabular}}
\end{table}

\subsection{Outline and notation} 
In Section~\ref{S: The testing procedure}, we introduce formally the testing procedure FOMT, and provide its statistical guarantees in terms of type-I error control, consistency and minimax optimal separation rates. Section~\ref{S: Compleixity Analysis} discusses the parameter choice strategy of FOMT and its computational complexity, together with a detailed comparison with existing procedures (Table~\ref{table: complexity}). Section~\ref{S: Adaptivity} studies the statistical and computational adaptivity of FOMT. Simulation studies are provided in Section~\ref{S: Simulation}. Technical details and proofs are given in the appendix.  

By $\mathbb{E}(X)$ and $\mathbb{V}(X)$, we denote the expectation and variance of the random variable $X$, respectively. For $x\in\mathbb{R}$,  we denote by $\floor{x}$ the largest integer that is smaller or equal to $x$, and similarly, by $\ceil{x}$ the smallest integer that is larger or equal to $x$. For $a,b\in \mathbb{R}$, we define $a\wedge b \coloneqq \min\{a,b\}$ and $a\vee b \coloneqq \max\{a,b\}$. For sequences $(a_n)_{n\in \mathbb{N}}$ and $(b_n)_{n\in \mathbb{N}}$ of positive numbers, we write $a_n\lesssim b_n$ or $a_n = O(b_n)$ if $a_n \leq C b_n$ for some finite constant $C > 0$. If $a_n \lesssim b_n$ and $b_n \lesssim a_n$, we write $a_n \asymp b_n$.  For readability, the multiplicative constants are suppressed in the main text via the above notation, while explicit constants are provided in the appendix.

\subsection{Code availability}
The implementation of FOMT, its adaptive variant and various existing monotonicity tests in \texttt{R} is provided, together with a documentation, on GitHub (\url{https://github.com/liuzhi1993/FOMT}).

\section{The basic testing procedure for known \texorpdfstring{$(\beta,L)$}{beta L}}
\label{S: The testing procedure}
For the construction of FOMT, we employ {local polynomial estimators} (see e.g.\ Section~\ref{s:lpe} in the appendix, and \citealp{tsybakov2008introduction,fan2018local}) for estimating $f$.
\begin{definition}\label{d:lpe}
Let $K:\mathbb{R} \rightarrow \mathbb{R}$ be a kernel, $h>0$ a bandwidth, and $k\ge 0$ an integer. Let also $U(x) = (1, x, x^2/2!, \ldots, x^k/k!)^\top$. For $x\in [0,1]$, the \emph{local polynomial estimator} (LPE) $\hat{f}_n(x)$ of order $k$ of $f(x)$ is defined as 
$\hat{f}_n(x) = \langle\hat\theta_n(x), U(0) \rangle$, where
\begin{equation*}
\hat{\theta}_n(x)\; \coloneqq\; \underset{\theta \in \mathbb{R}^{k+1}}{\mathrm{argmin}} \sum_{i=1}^{n}\left(Y_i-\left\langle\theta,\; U\left(\frac{x_i-x}{h}\right)\right\rangle\right)^{2}K\left(\frac{x_i-x}{h}\right).
\end{equation*}
\end{definition}
We make the following assumptions on the kernel $K$.
\begin{AoK}
\mbox{}
\begin{enumerate}[label=(\text{K\arabic*}),wide]
    \item \label{K1} The kernel $K:\R\to\R$ is non-negative, $\int_{-\infty}^{\infty} K(u)du = 1$ and its support belongs to $[-1,1]$. There exist constants $K_{\max}$, $K_{\min}$ and $\Delta$ such that
    $$
    0 < K_{\min}\1_{\{|u|\leq \Delta\}} \leq K(u) \leq K_{\max} < \infty, \quad \text{for all } u\in \R.$$
    \item \label{K2} The kernel $K$ is Lipschitz continuous, i.e., $K\in \Sigma(1, L_K)$ with $L_K\in (0,\infty)$.
    \item \label{K3} The kernel $K$ is symmetric, i.e.,
         $K(u)=K(-u)$ for all $u\in \R.$
\end{enumerate}
\end{AoK}
Assumptions \ref{K1} and \ref{K2} are standard conditions that guarantee the uniqueness of LPE as well as its minimax optimality in $L^p$-risk, $1\le p\le \infty$ (see e.g.\ \citealp{tsybakov2008introduction,korostelev2012minimax}). 
The symmetry Assumption~\ref{K3} is only needed to simplify technicalities (and could be relaxed). For instance, this assumption together with \hyperref[M1]{(M1)} ensures that the LPE $\hat f_n(x)$ is a non-negatively weighted sum of observations $Y_i$.

Before we come to the ({global}) testing problem in \eqref{eq:null}, we start with a \emph{local} testing problem, where the null hypothesis is
\(H_{i,j}: f(x_i)\leq f(x_j)\)    
with $x_i \equiv i/n<x_j \equiv j/n$ for some fixed $i,j \in [n]$. To this end, we will use the LPE $\hat{f}_n$ to estimate $f(x_i)$ and $f(x_j)$ by $\hat{f}_n(x_i)$ and $\hat{f}_n(x_j)$, respectively, and then to check whether the difference 
\begin{equation}
    \label{T_I,J}
    T_{i,j} \coloneqq \hat{f}_n(x_i)-\hat{f}_n(x_j),
\end{equation}
is significantly larger than zero. In particular, we employ the LPE $\hat{f}_n$ of $f$ with an optimally chosen bandwidth $h_n \asymp \bigl(\log(n)/n\bigr)^{1/(2\beta +1)}$, see its explicit formula in \eqref{e:bw} in the appendix. We reject the hypothesis $H_{i,j}$ if the statistic $T_{i,j}$ exceeds a {critical value} $q_{n,\beta,i,j} (\alpha) $, which needs to be chosen to control the type I error by $\alpha \in (0,1)$. That is, the \emph{local test} $\Phi_{i,j} \coloneqq \1_{\{T_{i,j}\geq q_{n,\beta,i,j}(\alpha)\}}$ is an $\alpha$-level test.

Towards the choice of the critical value $q_{n,\beta,i,j}(\alpha)$, we decompose $T_{i,j}$ in \eqref{T_I,J} into a sum of a deterministic part $D_{i,j}\equiv D_{i,j}(f)$ depending only on the signal $f$ and a random part $R_{i,j}\equiv R_{i,j}(\varepsilon)$ depending only on the noise $\varepsilon = (\varepsilon_i)$, and analyze each term individually. Due to the Hölder smoothness of $f\in H$, we can derive an upper bound of the deterministic part in form of  $D_{i,j}\lesssim h_n^{\beta}$, which can be improved to zero when $x_i$ and $x_j$ are away from the boundaries of the domain of $f$ (Proposition~\ref{p:upper bound Dab} in the appendix). Thus, we obtain
\begin{equation*}
    D_{i,j} \leq D_{n,\alpha,i,j} 
    \begin{cases}
        =0,&\quad \text{ if } x_i,x_j \in [h_n,1-h_n],\\
        \asymp h_n^{\beta},&\quad \text{ otherwise,} 
    \end{cases}
\end{equation*}
where the precise formula of $D_{n,\beta,i,j}$ is in \eqref{eq: D n i j} in the appendix.

{As the random term $R_{i,j}$ is normally distributed with mean zero, we can control its tail probability through an upper bound of its variance $\mathbb{V}(R_{i,j})\lesssim n^{-1}h_n^{-1}$ for all $i,j\in [n]$, which is obtained by treating separately the randomness of $\hat{f}_n(x_i)$ and $\hat{f}_n(x_j)$.} This is, however, not always optimal, particularly when $x_i$ and $x_j$ are close to each other. In this situation, $\hat{f}_n(x_i)$ and $\hat{f}_n(x_j)$ exhibit strong correlation, as the data points used by both estimates are largely overlapped.  Employing such a correlation leads to a refined upper bound (Theorem~\ref{Theorem: Bounds of variances} in the appendix) 
\begin{equation*}
    \V(R_{i,j}) \lesssim n^{-1}h_n^{-3}\left(|x_i-x_j| \wedge h_n\right)^2,\quad \text{ for all }\, i,j\in [n].
\end{equation*}
Combining this upper bound of variance with Mill's ratio, we can then find an explicit constant $C_{n,\alpha,i,j}$ (see \eqref{C,n,alpha,I,J} in the appendix), satisfying $\mathbb{P}(R_{i,j}\geq C_{n,\alpha,i,j})\leq \alpha$. Therefore, each $\Phi_{i,j}$ is an $\alpha$-level test if we  set $q_{n,\beta,i,j}(\alpha) = C_{n,\alpha,i,j}+D_{n,\beta,i,j}$.

Based on the proposed local tests, a typical approach to detect violations over the whole domain of $f$ is to \emph{scan} over all possible pairs of $i,j\in [n]$. However, this approach can be extremely time-consuming and may become computationally infeasible for large sample sizes. Alternatively, inspired by  \citet{ergun2000spot}, we propose to scan over only a (randomly generated) sparse collection of local tests, which has a sublinear cardinality and for which we can show that it is meanwhile able to locate all violations with high probability. 

\begin{figure}
\centering
\includegraphics[width=0.85\textwidth]{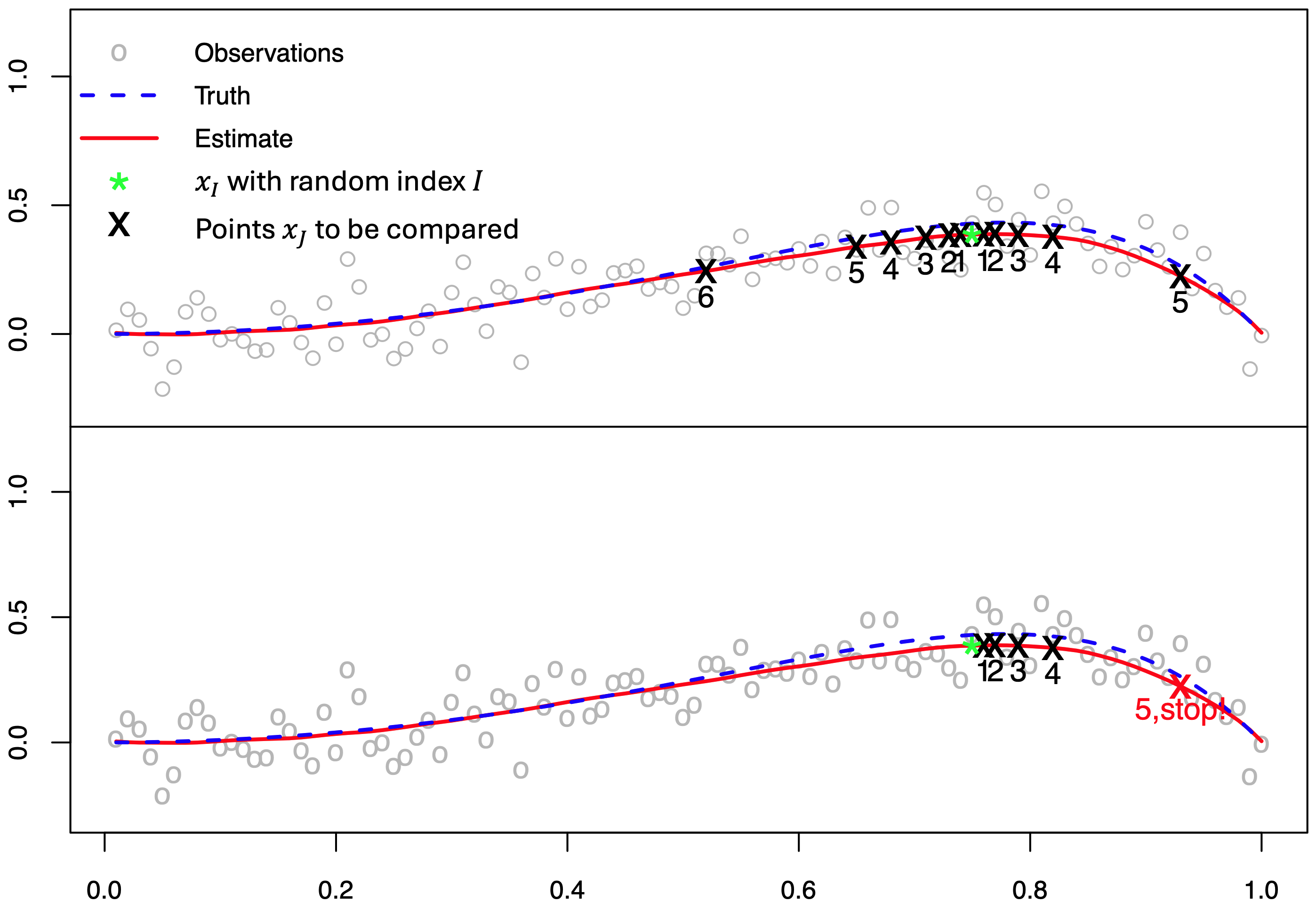}
\caption{Illustration of FOMT (\cref{alg:MC}). The true signal $f(x)= x^2-x^7$ (blue dashed line) and its LPE $\hat{f}_n$ (red solid line) together with $n = 100$ observations $Y_i$ (gray circles) are displayed in both panels with $x_I = 34/100$ and $x_I= 75/100$, respectively. For each $x_I$, locations of $x_J$ in an individual repetition of left and right searches, are highlighted by black crosses, respectively.
In the first panel, no violation is detected with $I = 34$, while in the second panel, a violation with $I = 75$ and $J=93$ (marked with a red cross) is detected.}
\label{F:MC}
\end{figure}

This sparse collection of local tests is generated as follows. We first select randomly $x_I$, with an index $I\sim \mathrm{Unif}([n])$ and then consider a few of $x_J$'s on both left and right sides of $x_I$. The locations $x_J$ are uniformly generated over intervals starting at $x_I$ with dyadically increasing lengths, which ensures efficient exploration of both \emph{nearby} and \emph{distant} indices. In order to increase the chance of finding violations associated with $f(x_I)$ and $f(x_J)$, we repeat the random generation of locations $x_I$ by around $h_n^{-1}\asymp \bigl(n/\log(n)\bigr)^{1/(2\beta+1)}$ times, (see \eqref{e: Csn} below), and for each $x_I$, we repeat the left and right searches for $x_J$ by $O(\log n)$ times. The resulting procedure is called \emph{FOMT} (\underline{F}ast and \underline{O}ptimal \underline{M}onotonicity \underline{T}est), see \cref{alg:MC}. An illustration is given in \cref{F:MC}. 

\alglanguage{pseudocode}
\begin{algorithm}[]
\begin{algorithmic}[1]
\Statex {\textbf{Input:} data $Y_1,\dots,Y_n$, and significance level $\alpha \in(0,1)$}
\Statex {\textbf{Parameters:} standard deviation $\sigma$, kernel function $K$, smoothness order $\beta$ and radius $L$}
\MRepeat{$C_n(\alpha)$}
    \State Generate $I \sim \mathrm{Unif}([n])$\;
    \If{$I \leq n-1$}
        \MRepeat{$\ceil{20 \log(n)}$}
            \For{$0 \leq k \leq \ceil{\log_2(n-I)}$}
                    \State Generate $J_k\sim \mathrm{Unif}([2^k\wedge(n-I)])$
                    \If{$\Phi_{I,I+J_k}=1$}
                        \Return $\Phi=1$\;
                    \EndIf
            \EndFor
        \EndRepeat
    \EndIf
    
    \If{$I \geq 2$}
        \MRepeat{$\ceil{20 \log(n)}$}
            \For{$0 \leq k \leq \ceil{\log_2(I-1)}$}
                    \State Generate $J_k'\sim \mathrm{Unif}([2^k\wedge(I-1)])$
                    \If{$\Phi_{I-J_k',I}=1$}
                        \Return $\Phi=1$\;
                    \EndIf
            \EndFor
        \EndRepeat
    \EndIf
\EndRepeat
\State \Return $\Phi=0$\;
\end{algorithmic}
\caption{FOMT: Fast and Optimal Monotonicity Test  $\Phi$}
\label{alg:MC}
\end{algorithm} 

By a simultaneous control on the type I error for each local test, we guarantee that the overall type I error of FOMT remains bounded by $\alpha$. Furthermore, we demonstrate that FOMT achieves a minimax optimal separation rate within the Hölder classes.

\begin{theorem}[Type I error]
\label{t:type I error}
Under the nonparametric regression model in \eqref{model}, suppose that Assumptions \hyperref[M1]{(M1)} and \ref{K1}--\ref{K3} hold. Let $H$ be defined in \eqref{eq:null} and $\alpha\in (0,1)$. Then, the FOMT $\Phi$ in \cref{alg:MC} is an $\alpha$-level test for $H$, namely, 
\begin{equation*}
    \mathbb{P}_{H}(\Phi = 1) \coloneqq \underset{f\in H}{\sup}\mathbb{P}_{f}(\Phi = 1) \leq \alpha, \quad \text{ for all } n \in \mathbb{N}.
\end{equation*}

\end{theorem}
For any $\delta>0$, define the classes $\mathcal{F}_{\beta}(\delta)$ of functions in the alternative hypothesis as 
\begin{equation*}
    \mathcal{F}_{\beta}(\delta)=
    \begin{cases}
        \left\{f\in \Sigma(\beta,L)\,|\, \max_{0\leq a<b\leq 1} f(a)-f(b)\geq \delta\right\}, &\quad \text{if } \beta \in (0, 1],\\
        \left\{f\in \Sigma(\beta,L)\,|\, \min_{x \in [0,1]}f'(x)\leq -\delta\right\}, &\quad \text{if } \beta \in (1, 2].
    \end{cases}
\end{equation*}
We show that FOMT is minimax optimal in the separation of $\mathcal{F}_{\beta}(\delta)$ from the null. 
\begin{theorem}[Minimax optimality]
\label{Theorem: Minimax rate optimality}
Under the nonparametric regression model in \eqref{model}, suppose that Assumptions~\hyperref[M1]{(M1)} and \ref{K1}--\ref{K3} hold. Let $\alpha \in (0,1)$, $h_n$ as in \eqref{e:bw}, and
\begin{equation}
    \label{e: Csn}
 C_{n}(\alpha) = -2\log\left(\frac{\alpha}{2}\right)h_n^{-1}  \quad\text{and}\quad
    \Delta_{\beta,n} = \left(\frac{\log n}{n}\right)^{\frac{\beta-\ceil{\beta}+1}{2\beta+1}}.
\end{equation}
Then the FOMT $\Phi$ is asymptotically minimax optimal in the sense that there exist two constants $0<c<C<\infty$ depending on $\beta$ such that the following statements hold:
\begin{enumerate}[label=(\roman*), wide]
    \item\label{i:opt:a}
    Upper bound:  
    \begin{equation*}
    \liminf_{n\to \infty}\underset{f\in \mathcal{F}_{\beta}(C\Delta_{\beta,n})}{\inf}\P_f (\Phi =1) \geq 1- \alpha.  
    \end{equation*}
    \item\label{i:opt:b} Lower bound: 
    For any $\alpha$-level test $\Psi$, namely, $\limsup_{n\to \infty}\mathbb{P}_{H}(\Psi = 1)\le \alpha$, 
    \begin{equation*}
        \underset{n \to \infty}{\limsup} \underset{f\in \mathcal{F}_{\beta}(c\Delta_{\beta,n})}{\inf} \P_f (\Psi=1) \leq
         \alpha.
    \end{equation*}
\end{enumerate}
The explicit formulae of $c$ and $C$ are given in the proof of \cref{Theorem: Minimax rate optimality} (cf.\ Section~\ref{ss:proof:rate} and the appendix).
\end{theorem}

The optimal separation rate $\Delta_{\beta,n}$ in \cref{Theorem: Minimax rate optimality} matches the rates reported in \citet[Theorems~3.1 and~3.2]{dumbgen2001multiscale}, who consider $\beta$ restricted to the values of $1$ and $2$, as well as the rates established in \citet[Theorem~4.2]{akakpo2014testing}, where $\beta \in (0,2]$ is considered for a Gaussian white noise model.

If $\sigma^2>0$ is unknown in the nonparametric regression model \eqref{model}, then it can be replaced by any estimators $\hat{\sigma}^2_n$ that are uniformly consistent over $\widetilde{\Sigma}(\beta,L)$ with $\beta \in (0,2]$ and $L>0$, given by 
\begin{equation*}
    \widetilde{\Sigma}(\beta,L) := \begin{cases}
        \Sigma(\beta,L), \quad& \text{ if } \beta \in (0,1],\\
        \{f\in \Sigma(\beta,L)\,|\, \norm{f'}_{\infty} \leq L\},\quad& \text{ if } \beta \in (1,2].
    \end{cases} 
\end{equation*}
\begin{corollary}
\label{c:unknown sigma}
Consider the nonparametric regression model in \eqref{model} with \emph{unknown} $\sigma^2$. Suppose that Assumptions~\hyperref[M1]{(M1)} and \ref{K1}--\ref{K3} hold. Let $\hat{\sigma}^2_n$ be a consistent estimator of $\sigma^2$ uniformly over $\widetilde{\Sigma}(\beta,L)$, that is, for any $\delta>0$, as $n \to \infty$,
\begin{equation}
    \label{ieq:ub F0}
    \underset{f\in \widetilde{\Sigma}(\beta,L)}{\sup}\mathbb{P}_f\left( \bigg|\frac{\hat{\sigma}_n^2}{\sigma^2}-1 \bigg| \geq \delta \right)\to 0.
\end{equation}
Let $\alpha \in (0,1)$, $\tilde{h}_n$ be the minimax optimal bandwidth given in \eqref{e:bw} with $\sigma$ replaced by $\hat{\sigma}_n$. We set $C_n(\alpha)$ as in \eqref{e: Csn}, and critical values $q_{n,\beta,i,j}(\alpha)$ accordingly with $\tilde{h}_n$ in place of $h_n$. Then: 
\begin{enumerate}[label=(\roman*), wide]
    \item \label{c:unknown sigma:a} FOMT is an asymptotically $\alpha$-level test for $\widetilde{H} \coloneqq \widetilde{\Sigma}(\beta,L)\cap \mathcal{M}$, i.e.,
    \begin{equation*}
    \limsup_{n\to \infty}\mathbb{P}_{\widetilde{H}}(\Phi = 1) \coloneqq \limsup_{n\to \infty}\underset{f\in \widetilde{H}}{\sup}\mathbb{P}_{f}(\Phi = 1) \leq \alpha.
\end{equation*}
\item \label{c:unknown sigma:b} For $C= 3C_{\beta}+2L$ with $C_{\beta}$ in \eqref{e:C0},
\begin{equation*}
    \liminf_{n\to \infty}\underset{f\in \mathcal{F}_{\beta}(C\Delta_{\beta,n})\cap \widetilde{\Sigma}(\beta,L)}{\inf}\P_f (\Phi =1) \geq 1- \alpha.  
\end{equation*}   
\end{enumerate}
\end{corollary}
\begin{remark}
\label{r: unknown variance}
The requirement in \eqref{ieq:ub F0} is rather weak, and can be fulfilled by many estimators (see e.g.\ \citealp{rice1984bandwidth,dette1998estimating,fan1998efficient,cai2008adaptive}). For instance, Rice, {\it et al.}~\cite{rice1984bandwidth} introduced 
\(\hat{\sigma}^2_n = \frac{1}{2(n-1)}\sum_{i=1}^{n-1}(Y_{i+1}-Y_i)^2,\)    
which satisfies,  uniformly for $f\in \widetilde{\Sigma}(\beta,L)$, 
\begin{equation*}
    \mathbb{E}(\hat{\sigma}^2_n) - \sigma^2 = \frac{1}{2(n-1)}\sum_{i=1}^{n-1} (f(x_{i+1})-f(x_i))^2 \lesssim  n^{-2(\beta \wedge 1)}
\quad \text{ and }\quad    \V(\hat{\sigma}^2_n) \lesssim \frac{1}{n}.
\end{equation*}
Thus, this estimator $\hat\sigma$ meets the requirements of \eqref{ieq:ub F0}, as 
\begin{equation*}
\mathbb{P}\left( \bigg|\frac{\hat{\sigma}_n^2}{\sigma^2}-1 \bigg| \geq \delta \right) \leq 
\frac{1}{\delta^2\sigma^4} \mathbb{E}(\abs{\hat{\sigma}_n^2-\sigma^2}^2)\lesssim \frac{1}{\delta^2\sigma^4} n^{-(4\beta \wedge 1)}\to 0.
\end{equation*}
\end{remark}    

\section{Computational complexity analysis}
\label{S: Compleixity Analysis}
\subsection{Computational complexity of FOMT}
\label{SS: CA of Phi}
\begin{theorem}[Computation]
\label{t:complexity1}
Under the nonparametric regression model in \eqref{model}, suppose that Assumptions~\hyperref[M1]{(M1)} and \ref{K1}--\ref{K3} hold. Let $\alpha\in (0,1)$, and $C_{n}(\alpha)$ as in \eqref{e: Csn}. Then:
\begin{enumerate}[label=(\roman*), wide]
    \item 
    FOMT has worst-case computational complexity $O(n (\log n)^2)$ for all $f\in \Sigma(\beta,L)$. 
     \item \label{i:comput Phi b} If further $f\in \mathcal{F}_{\beta}(C \Delta_{\beta,n})$ as in \cref{Theorem: Minimax rate optimality} with  $\Delta_{\beta,n} = \bigl(\log (n)/n\bigr)^{(\beta-\ceil{\beta}+1)/(2\beta+1)}$, then with asymptotic probability at least $1-\alpha$, FOMT detects the violation of $f$ in
    \begin{equation*}
        O\left( \bigl(\varepsilon_{\ceil{\beta}-1,\gamma_n}(f)\bigr)^{-1} \cdot n^{\frac{2\beta}{2\beta+1}} \left(\log n\right)^{\frac{4\beta+3}{2\beta+1}}\right)
    \end{equation*}
    steps. Here $\varepsilon_{\ceil{\beta}-1,\gamma_n}(f)$ is in Definition \ref{e: exceedance 0}, $\gamma_n = C_{\beta} h_n^{\beta -\ceil{\beta}+1} \asymp \bigl(\log(n)/n \bigr)^{(\beta -\ceil{\beta}+1)/(2\beta +1)}$ with the constant $C_{\beta}$ in \eqref{e:C0} and the optimal bandwidth $h_n$ in \eqref{e:bw}.
\end{enumerate}
\end{theorem}
    
        \cref{t:complexity1} \ref{i:comput Phi b} demonstrates that with asymptotic probability at least $1-\alpha$, FOMT detects a violation in the first $O\bigl(\bigl(\varepsilon_{\ceil{\beta}-1,\gamma_n}(f)\bigr)^{-1} \cdot n^{2\beta/(2\beta+1)} (\log n)^{(4\beta+3)/(2\beta+1)}\bigr)$ steps and terminates immediately.
    Note that FOMT does not require knowledge of $\varepsilon_{\ceil{\beta}-1,\gamma_n}(f)$. It has computational complexity adaptive to $\varepsilon_{\ceil{\beta}-1,\gamma_n}(f)$. If additionally $\varepsilon_{\ceil{\beta}-1,\gamma_n}(f)\gg n^{-1/(2\beta+1)}$, then its computational complexity is sublinear with high probability. Particularly,  $\varepsilon_{\ceil{\beta}-1,\gamma_n}(f)$ stays bounded from below as $n\to \infty$, for any fixed alternative $f\in \Sigma(\beta,L) \cap\mathcal{M}^c$.
    That is, FOMT has sublinear complexity $O\bigl(n^{2\beta/(2\beta+1)} (\log n)^{(4\beta+3)/(2\beta+1)}\bigr)$ for a fixed alternative. This is supported by the simulations in \cref{S: Simulation}.

{ To gain intuition on the role of the $\gamma$-exceedance fraction $\varepsilon_{\ceil{\beta}-1,\gamma_n}(f)$ for the computational complexity of FOMT, we focus on the case of $\beta \in (0,1]$ with $f \in \Sigma(\beta, L)$, noting that the case of $\beta \in (1,2]$ is analogous. Suppose that there exist indices $i,j \in [n]$ such that $f(x_i)-f(x_j)\geq \gamma_n$. If $x_i$ and $x_j$ were known, this $\gamma_n$-violation could be consistently detected using the local test $\Phi_{i,j}$. The rate $\gamma_n$ is minimax optimal for the testing problem in \eqref{eq:null}, delineating the boundary between detectable and undetectable alternatives (\cref{F:Phase} and \cref{Theorem: Minimax rate optimality}). Since $x_i$ and $x_j$ are generally unknown, we adopt a uniform sampling strategy. The probability of selecting an index associated with a $\gamma_n$-violation is given by $\varepsilon_{0,\gamma_n}(f)$ (\cref{l:epsilon_f} \ref{l:eps_g:b}). To ensure with high probability that at least one such index is sampled, we repeatedly generate $O\bigl(\bigl(\varepsilon_{0,\gamma_n}(f)\bigr)^{-1}\bigr)$ uniform indices. By Corollary \ref{p:>=hn}, we have $\varepsilon_{0,\gamma_n}(f) \geq h_n$ for all $f\in \mathcal{F}_{\beta}(C\Delta_{0,n})$, as defined in \cref{Theorem: Minimax rate optimality} (see also the green region above $h_n^{-1}$ in \cref{F:Phase}). Consequently, choosing  $C_n(\alpha) \asymp  h_n^{-1} \gtrsim \bigl(\varepsilon_{0,\gamma_n}(f)\bigr)^{-1}$, as in \eqref{e: Csn}, is sufficient to guarantee the minimax optimality (\cref{Theorem: Minimax rate optimality,t:consistency1}).}

\begin{remark}[Log-factor speed-up for $1 < \beta \le 2$]
In contrast to the case of  $\beta \in (0,1]$, where violations of monotonicity are characterized by $f(x_i)-f(x_j)$ for some $1\leq i<j\leq n$, violations under $\beta\in (1,2]$ can be measured using the first-order derivative $f'(x_i)$ at the design point $x_i$. Consequently, once an approximate match $x_I \approx x_i$ is identified, the search for $x_j$ in \cref{alg:MC} (lines 4--6 and 9--11) becomes unnecessary. Instead, it suffices to conduct local tests between adjacent design points $\Phi_{I,I+1}$ or $\Phi_{I-1,I}$, to determine whether $f'(x_I)\approx f'(x_i)$ is significantly negative. This observation motivates a simplified version of FOMT, termed \emph{S-FOMT} (see \cref{alg:MC2} in the appendix). This modification reduces its computational complexity from  $O\bigl(\bigl(\varepsilon_{1,\gamma}(f)\bigr)^{-1}\cdot nh_n \cdot \log^2 n\bigr)$ to $O\bigl(\bigl(\varepsilon_{1,{\gamma}}(f)\bigr)^{-1}\cdot nh_n\bigr)= O\bigl(\bigl(\varepsilon_{1,{\gamma}}(f)\bigr)^{-1}\cdot n^{2\beta/(2\beta+1)} (\log n)^{1/(2\beta+1)}\bigr)$ (cf. \cref{t:complexity1}), while maintaining its full detection power,  with at most a constant-factor loss, for all $f\in \Sigma(\beta,L)$ with $\beta\in (1,2]$. More precisely. \cref{t:type I error,Theorem: Minimax rate optimality} remain valid for S-FOMT. Moreover, as in \cref{t:complexity1}, one can show that with asymptotic probability at least $1-\alpha$, S-FOMT detects a violation of $f\in \mathcal{F}_{\beta}(C\Delta_{\beta,n})$ in $O(n)$ steps.
\end{remark}

\subsection{Comparison with other methods}
\label{SS: CA of other methods}
The computational complexity of existing monotonicity testing procedures may depend on three factors: estimation of $\sigma^2$, distribution of design points 
and the number of repetitions in Monte--Carlo or bootstrap procedures that are used to determine critical values. Towards a fair comparison, we consider the setup of \eqref{model},~i.e.
\begin{enumerate}
    \item[(A)] \label{A} Random errors $\varepsilon_i$ are i.i.d. Gaussian $\N(0,\sigma^2)$ with known $\sigma^2>0$.
    \item[(B)] \label{B} Design points\footnote{The statistical guarantees of  procedures by \citet{ghosal2000testing,chetverikov2019testing} are established for  random designs. Here we use regular deterministic designs only to evaluate their computational complexities.} $x_i = i/n$ for all $i = 1,\dots, n$.
\end{enumerate}
In \cref{table: complexity}, we provide comparison of the proposed FOMT with testing procedures by \citet{dumbgen2001multiscale,ghosal2000testing,chetverikov2019testing,baraud2005testing}. See \cref{DSGSVC,Details: Baraud Huet and Laurent's Procedure} for details. These results are supported by simulation studies in \cref{S: Simulation} (see particularly  \cref{F:times,F:large_times}).

\section{Adaptivity}
\label{S: Adaptivity}
If the smoothness parameter $\beta \in (0,2]$ is unknown, the optimal choice of bandwidth $h_n\asymp \bigl(\log (n)/n\bigr)^{1/(2\beta+1)}$ is not accessible. To tackle this problem, we introduce a computationally and statistically adaptive Lepskii principle for tuning $h$ and estimating $f$. It provides adaptively minimax optimal estimates, and meanwhile, its computational complexity also adapts to $\beta$. In this section, we work under the following assumption, instead of Assumption \hyperref[M1]{(M1)}.
\begin{assumption}
\begin{enumerate}
    \item[(M2)] \label{M2} $f\in \Sigma(\beta,L)$ with \emph{unknown} $\beta \in (0,2]$ and \emph{fixed} $L>0$.
\end{enumerate}
\end{assumption}
\subsection{A computationally adaptive Lepskii principle}
\label{SS: A new Lepskii principle}
Let $\hat{f}_n(x) \equiv \hat{f}_n(x;h)$ denote the LPE of order one for $f(x)$, $x\in [0,1]$ with bandwidth $h>0$, see Definition \ref{d:lpe}. 
For any $\mathcal{A}\subseteq [n]$, we define the semi-metric $d_\mathcal{A}$ by
\begin{equation*}
    d_{\mathcal{A}}\big(f,\hat{f}_n(\cdot;h)\big) \coloneqq \underset{i\in \mathcal{A}}{\max} \abs{f(x_i) - \hat{f}_n (x_i;h)}.
\end{equation*}
Clearly, it can be bounded from above by
\begin{align}
    \label{ieq:BV trade-off}
    \underset{i\in \mathcal{A}}{\max}\left| f(x_i)-\mathbb{E}\bigl(\hat{f}_{n}(x_i;h)\bigr)\right|  + \underset{i\in \mathcal{A}}{\max}\left| \hat{f}_n(x_i;h) -\mathbb{E}\bigl(\hat{f}_{n}(x_i;h)\bigr) \right| \eqqcolon B_{f,\mathcal{A}}(h) + \rho_{\mathcal{A}}(h).
\end{align}
By Hölder smoothness and \citet[Theorem~1.8]{tsybakov2008introduction}, there exist constants $c_1 = c_1(L,K)>0$ and $C_{\rho}= C_{\rho}(\sigma,K)>0$ (see \eqref{e:C rho} in the appendix) such that 
\begin{align*}
    B_{f,\mathcal{A}}(h) &\leq c_1 h^{\beta}\eqqcolon G_1(h),\qquad \qquad \qquad \quad\,\,\,\,\, \quad\text{ for all } h\in (0,1/2],\\
    \mathbb{E}\left(\rho_{\mathcal{A}}^2(h)\right) &\leq \mathbb{E}\left(\rho_{[n]}^2(h)\right)\leq C_{\rho}^2\frac{\log n}{nh}\eqqcolon G_2^2(h),\quad \text{ for all } h\in (0,1/2] \text{ and } \mathcal{A}\subseteq [n].     
\end{align*} 
Namely, $G_1$ and $G_2^2$ serve as upper bounds for the bias and variance, respectively. We consider an increasing sequence of $h_m = n^{-1} 4^{m-1}$ for $m \in [M]$ as bandwidth candidates with
$$M \equiv M_n = \frac{4}{5}\log_4 n + \frac{1}{5}\log_4 \log n.$$
The range of $(h_m)_{m\in[M]}$ covers the minimax optimal bandwidth $h_n$ 
under $\beta \in (0,2]$. For simplicity, for all $m\in[M]$, we set 
\begin{equation*}
   \hat{f}_{n,m} \coloneqq \hat{f}_{n}(\cdot;h_m),\quad \rho(m)\coloneqq \rho(h_m),\quad G_1(m) \coloneqq G_1(h_m)  \quad \text{ and }\quad G_2(m) \coloneqq G_2(h_m).
\end{equation*}
For $\kappa>1$, we introduce the \emph{Computationally Adaptive Lepskii Method} (CALM;\ \cref{Lepskii}) with the stopping index   
\begin{equation}
    \label{e: Lepskii}
    \bar{m} \coloneqq \min\{m\,|\, \text{ there exists } k\in [m-1] \text{ such that } d_{\mathcal{A}}(\hat{f}_{n,k},\hat{f}_{n,m})> 4 \kappa G_2(k)\}-1.
\end{equation}

\begin{theorem}\label{theorem: Lepskii}
Consider the nonparametric regression model in \eqref{model} with Assumptions~\hyperref[M2]{(M2)} and \ref{K1}--\ref{K3}. Let $\kappa>1$ and $\bar{m}$ be in \eqref{e: Lepskii}, $K_{\max}$ be in Assumption \ref{K1} and $\mu_2 \coloneqq \int_{-\infty}^{\infty} K^2(u)du$. Then there is a constant $C_{\max} \equiv C_{\max}(\kappa)>0$ such that 
\begin{equation*}
    \underset{\mathcal{A}\subseteq [n]}{\inf}\underset{f\in \Sigma(\beta,L)}{\inf}\mathbb{P}\left(d_\mathcal{A}(f,\hat{f}_{n,\bar{m}})  \leq C_{\max} \left(\frac{\log n}{n}\right)^{\frac{\beta}{2\beta +1}}\right) \geq 1- n^{-\frac{\mu_2 }{4K^2_{\max}}(\kappa-1)^2} \log n,
\end{equation*}
i.e., $\hat{f}_{n,\Bar{m}}$ attains the optimal estimation rates for $\Sigma(\beta,L)$ with probability tending to one.
\end{theorem}

\begin{remark}\label{r:Lepskii}
Following the standard Lepskii principle (cf. \citealp{lepskii1991problem,mathe2006regularization}), one would define the stopping index $\tilde{m}$ as 
\begin{equation}
    \label{e: standard Lepskii}
    \tilde{m} \coloneqq \max \{m\,|\, d_{\mathcal{A}}(\hat{f}_{n,k}, \hat{f}_{n,m})\leq 4 \kappa G_2(k),\, \text{ for all } k\leq m\}.
\end{equation}
By the same proof technique as in \cref{theorem: Lepskii}, we can show that $\hat{f}_{\tilde{m}}$ achieves adaptive minimax optimality as well. However, this approach is computationally suboptimal. To see this, suppose that we know that $m_*\in [M]\,$ fails to satisfy the constraint in \eqref{e: standard Lepskii}. We cannot rule out the possibility that a later index $m>m_*$ might still satisfy the condition in \eqref{e: standard Lepskii}. Thus, the worst case computational complexity of $\tilde{m}$ is $O\bigl(\abs{\mathcal{A}} \cdot n^{4/5}(\log n)^{1/5}\bigr)$, which can be attained. In contrast, our CALM in \eqref{e: Lepskii} addresses this issue by terminating at the first time that the condition in \eqref{e: Lepskii} is violated. This modification enables more efficient computation of $(\hat{f}_{n,\bar{m}}(x_i))_{i\in \mathcal{A}}$ and the corresponding bandwidth $h_{\Bar{m}}$, see \cref{Lepskii}. {In the literature, CALM has been recognized as an equivalent formulation of the standard Lepskii principle (see \citet{reiss2009pointwise,REISS2012} in quantile regression). However, the computational advantage highlighted here has remained unnoticed, to the best of our knowledge. In retrospect, CALM should be preferred in other settings, given its equivalent statistical efficiency but significantly improved computational efficiency. }
\end{remark}

\begin{figure}
\centering
\includegraphics[width=0.75\textwidth]{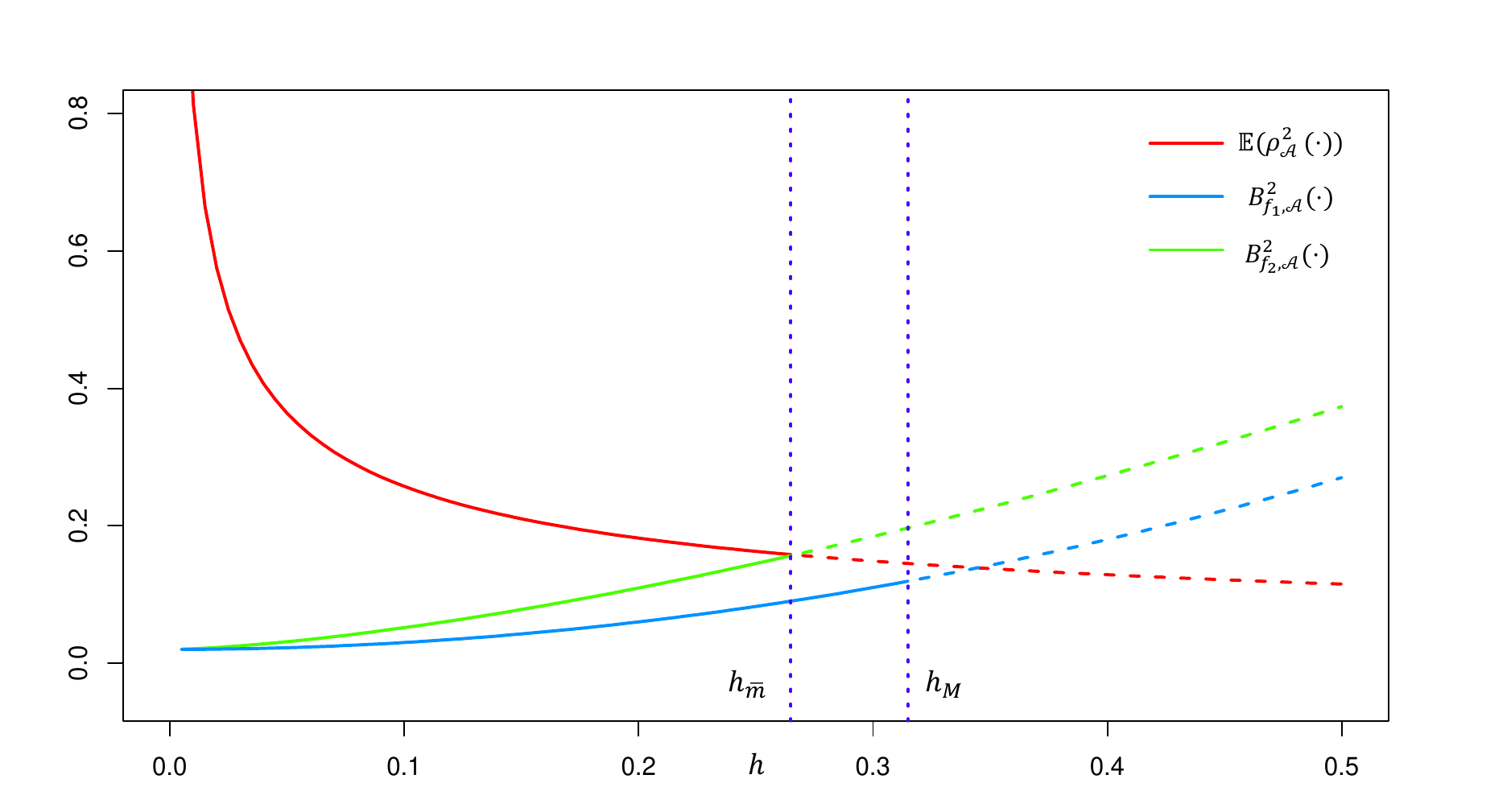}
\caption{Squared biases versus variances of LPE over various choices of bandwidth, see \eqref{ieq:BV trade-off}. We illustrate two scenarios: $f_1\in \Sigma(\beta,L)$ with $\beta>2$, and $f_2\in \mathcal{C}_{\mathcal{A},\beta}$ in \eqref{e: defn class C} with $0<\beta\le2$. In case of $f_1$, CALM (\cref{Lepskii}) returns $h = h_M$ while in the case of $f_2$, CALM returns $h_{\Bar{m}}\lesssim \bigl(\log(n)/n\bigr)^{1/(2\beta+1)}$}
\label{F:Lepskii}
\end{figure}

\alglanguage{pseudocode}
\begin{algorithm}[H]
\begin{algorithmic}[1]
\Statex {\textbf{Input:} data $Y_1,\dots,Y_n$, constant $\kappa>1$ and subset $\mathcal{A}\subseteq [n]$}
\Statex {\textbf{Parameters:} standard deviation $\sigma$, kernel function $K$}
    \For{$1 \leq m \leq M$}
        \State Compute $(\hat{f}_{n,m}(x_i))_{i\in \mathcal{A}}$ and set $\bar{m}=m$\;
        \For{$1 \leq k \leq m-1$}
            \If{$d_{\mathcal{A}}(\hat{f}_{n,k},\hat{f}_{n,m})> 4 \kappa G_2(k)$}
                \State $\bar{m} = \bar{m}-1$\;
                \State \Return $(\hat{f}_{n,\bar{m}}(x_i))_{i\in \mathcal{A}}$ and $h_{\Bar{m}}$\;
            \EndIf
        \EndFor
    \EndFor
\State \Return $(\hat{f}_{n,\bar{m}}(x_i))_{i\in \mathcal{A}}$ and $h_{\Bar{m}}$\;
\end{algorithmic}
\caption{CALM: Computationally Adaptive Lepskii Method}
\label{Lepskii}
\end{algorithm}
Similar to the standard Lepskii principle, CALM balances bias and variance and terminates when $B_{f,\mathcal{A}}$ and $\rho_{\mathcal{A}}$ in \eqref{ieq:BV trade-off} reach the same order (if possible). If $B_{f,\mathcal{A}}(h)\lesssim h^{\beta'}$ with $\beta'>2$, then $B_{f,\mathcal{A}}$ is equal to $G_2$ at some point in $(h_M,1/2]$ instead of $(0,h_M]$. Thus, CALM will terminate with $\Bar{m}=M$ and return estimates with the largest bandwidth $h_M$ (see \cref{F:Lepskii}). To quantify the intrinsic smoothness, we consider the following \emph{self-similarity} assumption:
\begin{assumption}[Self-Similarity] 
\mbox{}
\label{A:SS} For some $ \beta \in (0,2]$ and constants $c_{\mathcal{A}},c_1>0$, it holds
\begin{equation}
    \label{ieq:LB}
    c_{\mathcal{A}} h^{\beta} \leq B_{f,\mathcal{A}}(h) \equiv \underset{i\in \mathcal{A}}{\max}\left| f(x_i)-\mathbb{E}\left(\hat{f}_{n}(x_i;h)\right)\right|\leq c_1h^{\beta} ,\quad \,\text{ for all } h\in (0,1/2].
\end{equation}
\end{assumption}

\begin{remark}
A similar \emph{self-similarity condition} was introduced by \citet{gine2010confidence} for kernel density estimators $\hat{p}_h$  with bandwidth $h>0$ in the form of 
\begin{equation}
    \label{e:SS density}
    b_1 h^{\beta} \leq \norm{\mathbb{E}(\hat{p}_h)-f}_{\infty} \leq b_2 h^{\beta},\quad \text{for all } h\in (0,h_0],
\end{equation}
where $0<b_1\leq b_2$ are two constants and $h_0$ is a constant bandwidth. The condition in \eqref{e:SS density} is shown to be crucial for the existence and construction of adaptive confidence bands of density functions. It is also known that \eqref{e:SS density} is a rather weak requirement as it remains valid in Hölder classes, except for a \enquote{topologically small} subset. More precisely, the exceptional set is nowhere dense with respect to the norm topology in the Hölder class (\citealp[Proposition~4 and discussion therein]{gine2010confidence}). We refer to \citet{gine2010confidence,bull2012honest,bull2013spatially,chernozhukov2014anti,armstrong2020simple} for self-similarity conditions for \emph{wavelet projection estimators}. We emphasize that our self-similarity condition in \eqref{ieq:LB} is required only for an accurate analysis of computational complexity, instead of adaptive minimax optimality. Clearly, the self-similarity condition in \eqref{ieq:LB} may fail for some functions in $\Sigma(\beta,L)$, for instance, $f \equiv 0$ for which it holds $B_{0,\mathcal{A}} \equiv 0$ with arbitrary $\mathcal{A}\subseteq[n]$ and $h\in(0,1/2]$. It remains open whether the exceptional set of \eqref{ieq:LB} is also \enquote{topologically small}.
We leave this for further research.
\end{remark}

For $\beta \in (0,2]$, we focus on the class $\mathcal{C}_{\mathcal{A},\beta}$ given by
\begin{equation}
    \label{e: defn class C}
    \mathcal{C}_{\mathcal{A},\beta} = \{f\in \Sigma(\beta,L) \,|\, \text{self-similarity assumption \eqref{ieq:LB} holds with } \beta \}.
\end{equation}
Intuitively, $\mathcal{C}_{\mathcal{A},\beta}$ consists of all functions in $\Sigma(\beta,L)$, whose smoothness order is no larger than $\beta$.
For all $f\in \mathcal{C}_{\mathcal{A},\beta}$, CALM runs in sublinear time 
with probability tending to one. 

\begin{lemma}\label{theorem: Lepskii complexity}
Under the nonparametric regression model \eqref{model}, suppose that Assumptions~\hyperref[M2]{(M2)} and \ref{K1}--\ref{K3} hold. Let $\mu_2 = \int_{-\infty}^{\infty} K^2(u)du$. Then:
\begin{enumerate}[label=(\roman*), wide]
    \item \label{i:Lepskii compleixity:a}     
    The worst case computational complexity of CALM (\cref{Lepskii}) is
    \begin{equation*}
        O\left(\abs{\mathcal{A}} \cdot n^{\frac{4}{5}} (\log n)^{\frac{1}{5}}\right).
    \end{equation*}
    \item \label{i:Lepskii compleixity:b} If $f\in \mathcal{C}_{\mathcal{A},\beta}$ in \eqref{e: defn class C}, then with probability no less than $1-2n^{-\mu_2(\kappa-1)^2/(4K^2_{\max})}\cdot \log n$, CALM (\cref{Lepskii}) has computational complexity 
\begin{equation*}
    O\left(\abs{\mathcal{A}} \cdot n^{\frac{2\beta}{2\beta+1}}\left(\log n\right)^{\frac{1}{2\beta+1}}\right).
\end{equation*}
\end{enumerate}
\end{lemma}
\begin{remark}
If $f$ fulfills the self-similarity condition with $0 < \beta\leq 2$, i.e., $f\in \mathcal{C}_{\mathcal{A},\beta}$ in \eqref{e: defn class C}, then by Lemma \ref{theorem: Lepskii complexity} 
we have
\begin{equation*}
    h_{\Bar{m}} \asymp \left(\frac{\log n}{n}\right)^{\frac{1}{2\beta+1}}.
\end{equation*}
This reveals that CALM is indeed computationally and statistically adaptive.
\end{remark}

\subsection{Adaptive FOMT}
\label{SS: Adaptive MC}
Using CALM (\cref{Lepskii}) we now introduce an adaptive modification (A-FOMT) of FOMT (\cref{alg:MC}) that works as follows (see \cref{alg:MC Lepskii}):
\begin{enumerate}[label=\roman*.]
    \item \emph{Generation of pairs of indices}: Generate a realization $i$ from the distribution $ \mathrm{Unif}([n])$. Construct two sequences of pairs $\{(i,i+j)\}_j$ and $\{(i-j,i)\}_j$ for searching to the right and to the left starting from $i$, respectively, as described in FOMT. Let $\mathcal{P}$ represent the set of all generated pairs, and $\mathcal{A}$ the set of all indices from $\mathcal{P}$.
    \item \emph{Computation of estimates}: Compute the estimates $(\hat{f}_{n,\Bar{m}}(x_i))_{i\in \mathcal{A}}$ and $h_{\Bar{m}}$, via CALM.
    \item \emph{Testing for violations}: For each pair $(i,j)\in \mathcal{P}$, apply the local tests $\Phi_{i,j}$ with the critical value $C_{n,\alpha,i,j}$ in \eqref{C,n,alpha,I,J} based on $h_{\Bar{m}}$. A-FOMT will reject the null hypothesis in \eqref{eq:null} if and only if some violation of monotonicity $\Phi_{i,j}$  is detected.
    \item Repeat Steps 1--3 $C_{n}(\alpha) = -2\log(\alpha/2) \cdot n/\log (n)$ times.
\end{enumerate}
\algdef{SE}[REPEATN]{RepeatN}{End}[1]{\algorithmicrepeat\ #1 \textbf{times}}{\algorithmicend}
\alglanguage{pseudocode}
\begin{algorithm}[H]
\begin{algorithmic}[1]
\Statex {\textbf{Input:} data $Y_1,\dots,Y_n$, a constant $\kappa>1$ and significance level $\alpha \in(0,1)$}
\Statex {\textbf{Parameters:} standard deviation $\sigma$, kernel function $K$ and radius $L$}
\State $C_n(\alpha)= -\log(\alpha/2) \cdot (n/\log (n))$
\For{$1\leq l\leq C_{n}(\alpha)$}
    \State Generate $I\sim \mathrm{Unif}([n])$
    \State Generate $(\mathcal{P},\mathcal{A})$ with \cref{alg:IG} in the appendix with parameters $(n,I)$
    \State Compute estimates $(\hat{f}_{n,\Bar{m}}(x_i))_{i\in \mathcal{A}}$ and $h_{\Bar{m}}(\mathcal{A})$ with \cref{Lepskii}
    \For{$(i,j)\in \mathcal{P}$}
        \State Compute $C_{n,\alpha,i,j}$ in \eqref{C,n,alpha,I,J} with $h = h_{\Bar{m}}$
        \If{$\hat{f}_{n,\Bar{m}}(x_i)-\hat{f}_{n,\Bar{m}}(x_j)\geq C_{n,\alpha,i,j}$}
            \Return $\Phi_{A}=1$\;   
        \EndIf
    \EndFor
\EndFor
\State \Return $\Phi_{A}=0$\;
\end{algorithmic}
\caption{A-FOMT: Adaptive FOMT $\Phi_{A}$}
\label{alg:MC Lepskii}
\end{algorithm}

\begin{theorem}\label{theorem: Lepskii consistency}
Under the nonparametric regression model \eqref{model}, suppose that Assumptions \hyperref[M2]{(M2)} and \ref{K1}--\ref{K3} hold.
Let $\alpha \in (0,1)$, $\Delta_{\beta,n} = \bigl(\log (n)/n\bigr)^{(\beta-\ceil{\beta}+1)/(2\beta+1)}$ and $\kappa> 1+2K_{\max}/\sqrt{\mu_2}$, where $K_{\max}$ is in Assumption \ref{K1} and $\mu_2 = \int_{-\infty}^{\infty} K^2(u)du$. 
Then, for A-FOMT $\Phi_{A}$ in \cref{alg:MC Lepskii}, there exists a sufficiently large constant $C$ such that 
\begin{equation*}
    \liminf_{n \to \infty}\underset{f\in \mathcal{F}_\beta(C\Delta_{\beta,n})}{\inf}\P_f (\Phi_{A} =1) \geq 1- \alpha.
\end{equation*}
\end{theorem}
The comparison between \cref{theorem: Lepskii consistency,Theorem: Minimax rate optimality} reveals that, up to a constant factor, A-FOMT achieves adaptively minimax detection optimality. It means that even when the smoothness parameter $\beta \in(0,2]$ is unknown, A-FOMT retains its full detection power, losing at most a constant factor in performance.

Let 
\begin{equation}
\label{e: gloabl SS}
\mathcal{C}_{\beta} = \left\{f\in \Sigma(\beta,L)\,\bigg|\, f\in \mathcal{C}_{\mathcal{A},\beta} \text{ for all } \mathcal{A} \text{ with } \abs{\mathcal{A}} \gtrsim \left(\log n\right)^2\right\},  
\end{equation}
i.e., $\mathcal{C}_{\beta}$ consists of functions such that fulfilling the self-similarity condition in \eqref{e: defn class C} on all subsets $\mathcal{A}$ with not too small cardinality.
\begin{theorem} 
\label{t:MC Lepskii compleixity}
Assume that all conditions in \cref{theorem: Lepskii consistency} hold. Let $f\in\mathcal{F}_{\beta}(C\Delta_{\beta,n})$ for some sufficiently large $C$. 
Then: 
\begin{enumerate}[label=(\roman*), wide]
    \item \label{i:tMC:a}
    With probability at least $1-\alpha$,  A-FOMT (\cref{alg:MC Lepskii}) has computational complexity
\begin{equation*}
O\left(\bigl(\varepsilon_{\ceil{\beta}-1,\gamma_n}(f)\bigr)^{-1} \cdot n^{\frac{4}{5}}\left(\log n\right)^{\frac{11}{5}} \right),
\end{equation*}
where $\gamma_{n} = C_{\beta}\Delta_{\beta,n}$ for some suitable constant $C_{\beta}$ depending on $\beta$.
\item  \label{i:tMC:b}
Additionally, if $f\in \mathcal{C}_{\beta}$ in \eqref{e: gloabl SS}, then with probability at least $1-\alpha$, the computational complexity of A-FOMT (\cref{alg:MC Lepskii}) is 
\begin{equation*}
    O\left(\bigl(\varepsilon_{\ceil{\beta}-1,\gamma_n}(f)\bigr)^{-1} \cdot n^{\frac{2\beta}{2\beta+1}}\left(\log n\right)^{\frac{4\beta+3}{2\beta+1}} \right).
\end{equation*}
\item \label{i:tMC:c}
In the worst-case, the computational complexity of A-FOMT (\cref{alg:MC Lepskii}) is 
\begin{equation*}
    O\left(C_{n}(\alpha) \cdot n^{\frac{4}{5}} (\log n)^{\frac{11}{5}} \right) = O\left(n^{\frac{9}{5}} (\log n)^{\frac{6}{5}}\right).
\end{equation*}
\end{enumerate}
\end{theorem}

\section{Simulation studies}
\label{S: Simulation}
Our simulation study consists of two parts. In the first, we evaluate the empirical performance of FOMT and its adaptive variant, A-FOMT, across various scenarios with small sample sizes $n \in \{400,800,1200,1600,2000,2400,2800,3200\}$. To benchmark the performance, we include as competitors the adaptive minimax optimal testing procedures introduced by \citet{dumbgen2001multiscale,akakpo2014testing,chetverikov2019testing}, denoted by DS, ABD and C, respectively. The second part focuses on the scalability of FOMT and A-FOMT to large datasets, considering sample sizes $n \in 10^5\times [10]$. Throughout the study, we set the significance level to $\alpha = 0.05$ and noise level $\sigma = 0.3$. For DS, ABD and C, their critical values are computed via Monte--Carlo simulations with $100$ repetitions. This approach provides accurate control of type I errors but it is extremely time-consuming and becomes computationally infeasible as the sample size increases. For example, the computation time of ABD (for $n = 1600$) and C (for $n=800$) exceeds two hours (\cref{T:critical}) on a standard Laptop (Intel i5-CPU 1.40GHz with four cores and 16GB RAM). Regarding FOMT and A-FOMT, we emphasize that the constants $W$ in \eqref{eq: critical} and $C_{\rho}$ in \eqref{e:C rho} are chosen for asymptotic analysis, and may be suboptimal for finite samples. Based on simulations (not shown), we recommend replacing them by $0.58W$, $0.00175C_{\rho}$, respectively. Further, the constant $20$ of FOMT (\cref{alg:MC} in lines 4 and 9) and A-FOMT (cf. \cref{alg:IG} in lines 3 and 11) is due to technical reasons and overly conservative in practice; We instead use  $0.1$ in both algorithms. Besides, we use LPE of order one with the \emph{Epanechnikov kernel} and bandwidth $h_n = 0.3 \bigl(\log(n)/n\bigr)^{1/3}$ in FOMT. In A-FOMT, we replace the base of exponential bandwidth grid from $4$ to $1.6$ to improve accuracy.

\begin{figure}
\centering
\includegraphics[width=0.95\textwidth]{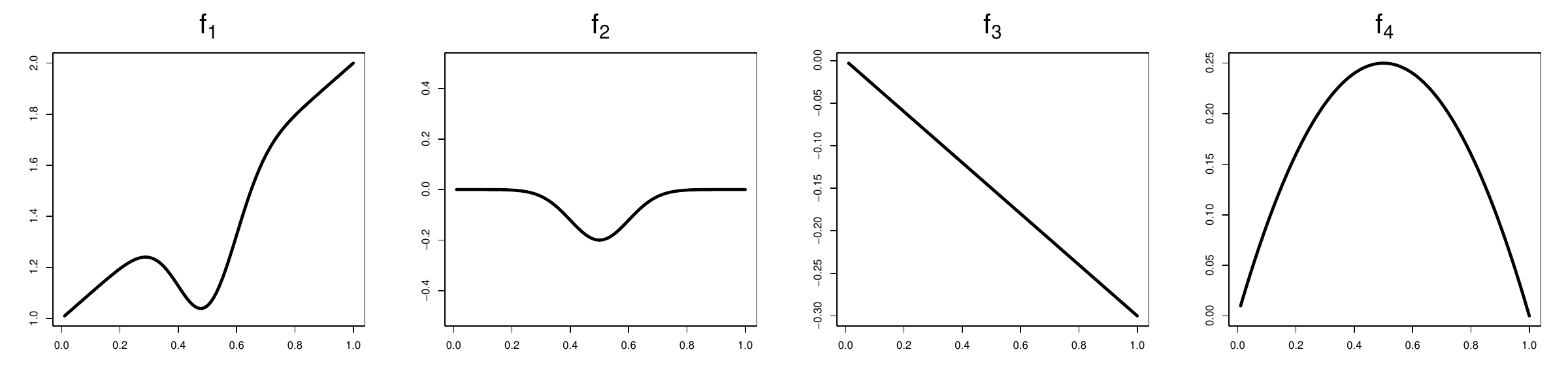}
\caption{Test signals $f_1(x) = \,1+x-0.45 \cdot \exp{(-50(x-0.5)^2)}$, $f_2(x) = \,-0.2\cdot \exp{(-50(x-0.5)^2)}$, $f_3(x) = \,-0.3 x$ and $f_4(x) = \,x(1-x)$. 
The function $f_1$ is a specific case in \citet{gijbels2000tests} with $a = 0.45$, while $f_2$ is considered by \citet{baraud2005testing}. The function $f_3$ decreases linearly with slope $-0.3$. The function $f_4$ is a \enquote{arch-like} function.} 
\label{F:signals}
\end{figure}

\begin{figure}
\centering
    \includegraphics[width=.75\textwidth]{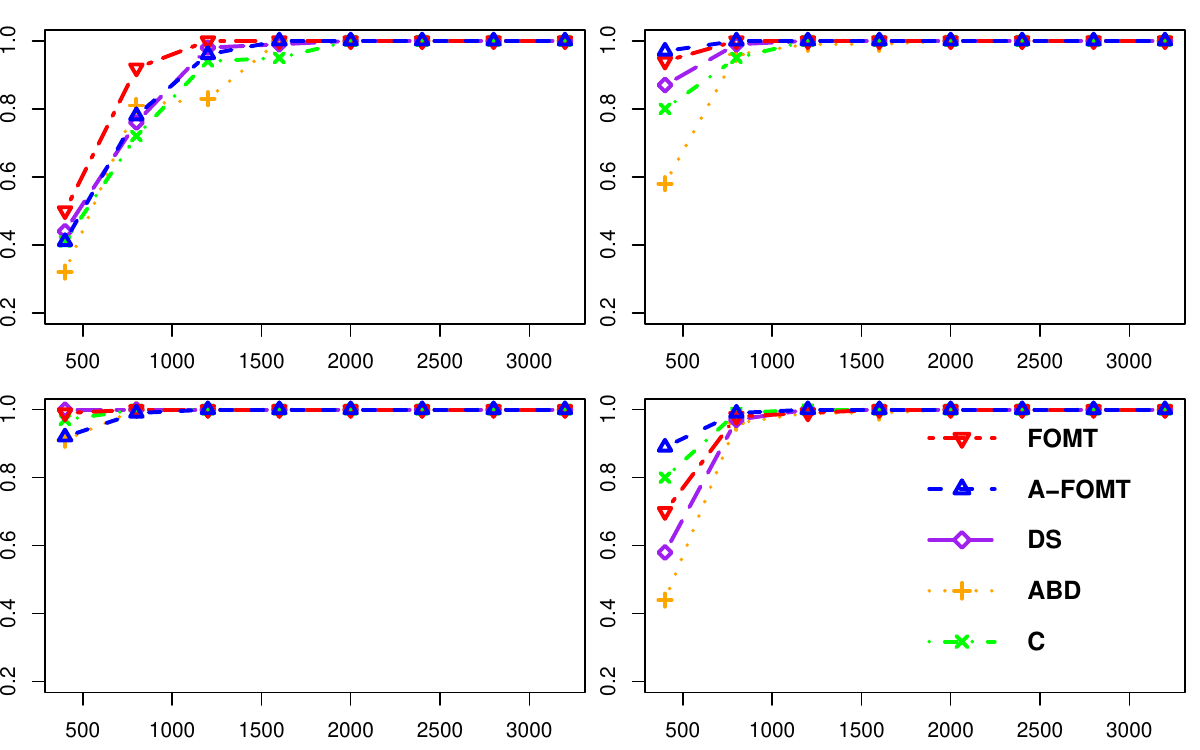}
\caption{Detection powers under the alternatives $f_i$ with $i = 1, \dots, 4$ of FOMT, A-FOMT, DS \citep{dumbgen2001multiscale}, ABD \citep{akakpo2014testing} and C \citep{chetverikov2019testing}, averaged over $100$ repetitions, for various sample sizes.} 
\label{F:powers2}
\end{figure}

\begin{figure}
\centering
\includegraphics[width=.85\textwidth]{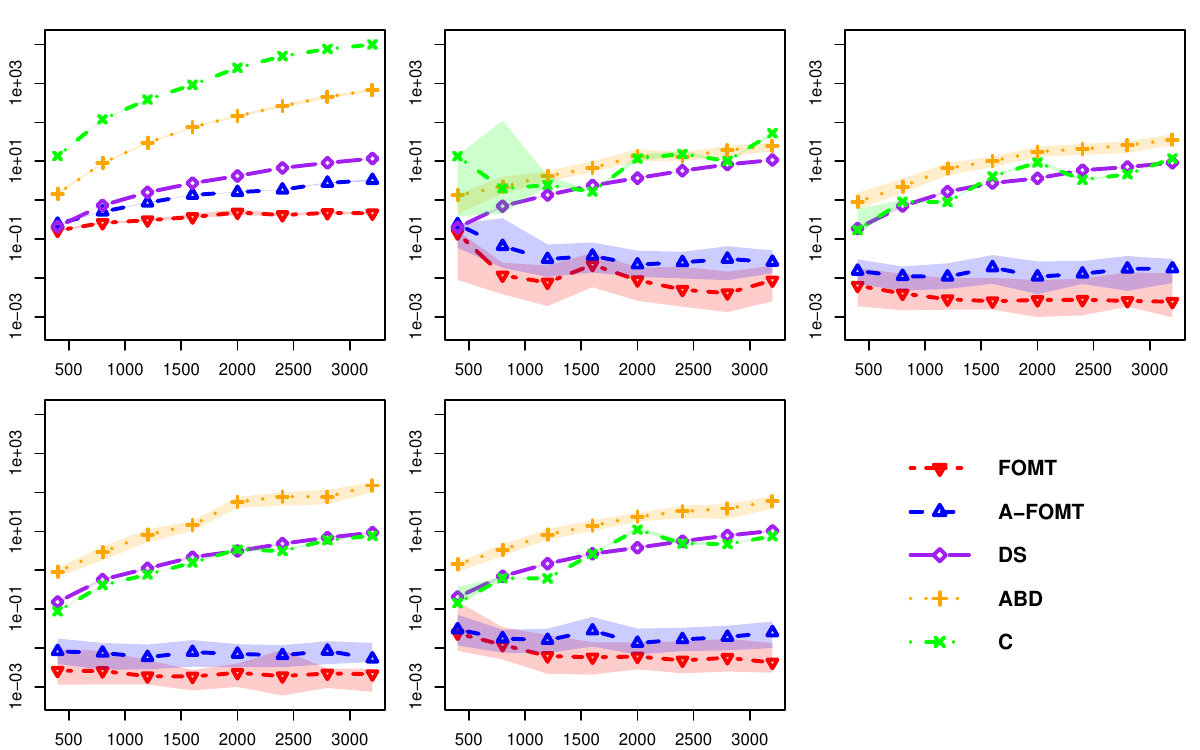}
\caption{Computation times (in seconds) of FOMT, A-FOMT, DS \citep{dumbgen2001multiscale}, ABD \citep{akakpo2014testing} and C \citep{chetverikov2019testing}. Results for $f_i$ for $i = 0,\dots,4$ are displayed in the first to fifth panels, respectively. For each signal and each method, the curve shows the median of computation times, and the shaded region marks $25\%$- and $75\%$-quantile curves, over $100$ repetitions.} 
\label{F:times}
\end{figure}

\begin{figure}
\centering
\includegraphics[width=.85\textwidth]{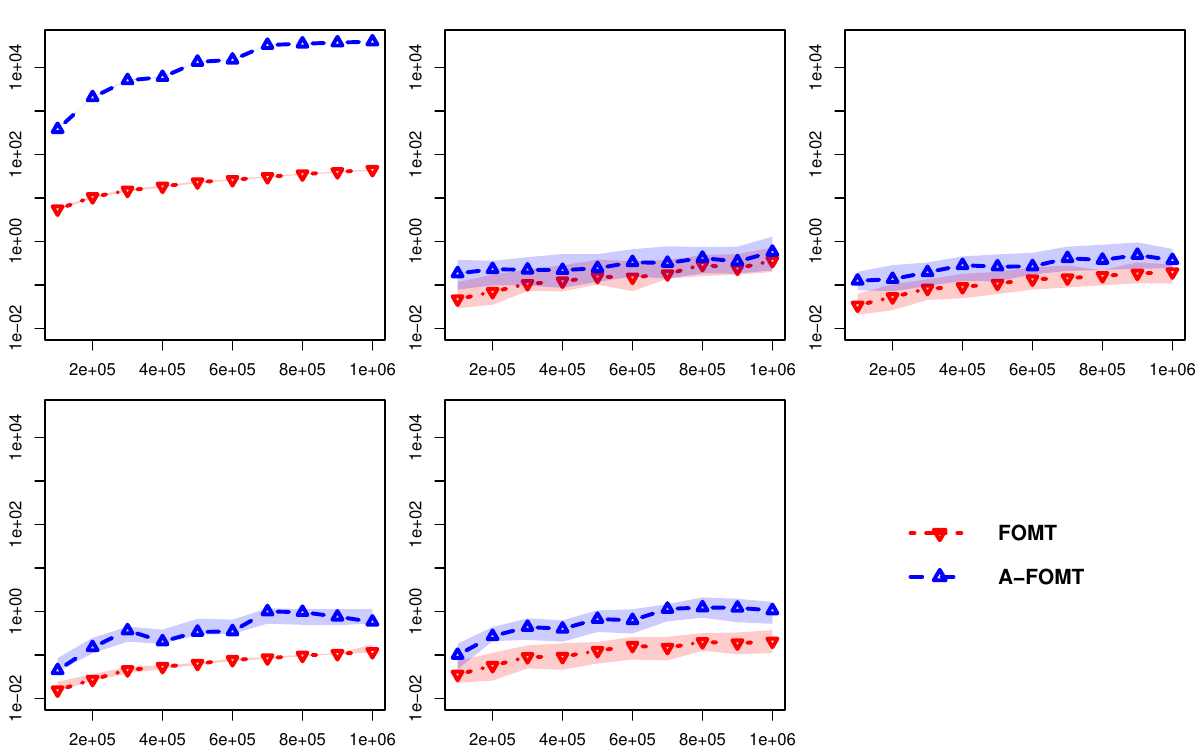}
\caption{Computation times (in seconds) of FOMT and A-FOMT for large sample sizes in the same setup as in \cref{F:times}. For each signal, the curve shows the median of computation times, and the shaded region marks $25\%$- and $75\%$-quantile curves, over $100$ repetitions for FOMT and over $10$ repetitions for A-FOMT. 
} 
\label{F:large_times}
\end{figure}

For the test signal in the null hypothesis, we consider $f_0 \equiv 0$, which corresponds to the most different signal to distinguish from the alternatives. The test signals in the alternative are shown in \cref{F:signals}. 
To ensure fair comparison in computation times, all considered approaches are implemented in the \texttt{R} language, with code available on Github (\url{https://github.com/liuzhi1993/FOMT}).
{Simulations (not displayed) reveal that all procedures exhibit almost the same control over Type I error.  We thus restrict ourselves to report the performance} of each testing procedure in terms of detection power, computation time and separately simulation time of critical values, in \cref{F:powers2,F:times,F:large_times}, with additional details in \cref{T:time,T:signals,T:critical,T:time large scale} in the appendix. 

All methods demonstrate comparable statistical performance, aligning with their minimax optimality, see \cref{F:powers2}. For $f_1$, FOMT outperforms the others slightly, with A-FOMT and DS following as two close competitors. For $f_2$, A-FOMT achieves the best detection power, with FOMT and DS trailing closely behind.  For $f_3$, DS and FOMT are slightly better than the others. For $f_4$, A-FOMT and C exhibit the best detection powers, followed by FOMT and DS. All procedures have almost perfect performance for large sample sizes. 

In \cref{F:times}, FOMT is the fastest and A-FOMT is only slightly slower.
In the scenarios with $f_0$ and $f_1$ and small sample size $n = 400$, DS has a slightly faster computational speed than A-FOMT, but its computational time becomes significant longer than that of A-FOMT for large sample sizes $n\geq 800$. In particular, A-FOMT is at least five times faster than DS, for alternatives with sample sizes $n\geq 1200$. Note that ABD and C are the slowest. 

Interestingly, as illustrated in \cref{F:times,F:large_times}, the computation times of FOMT and A-FOMT do not increase significantly with sample size and may even slightly decrease for alternatives within a medium range of sample sizes. The impact of sample size on computational time is twofold: On the one hand, increasing $n$ raises the computational cost of each local test; on the other hand, larger sample sizes enhance detection power, reducing the number of required local tests.
Consequently, for alternatives \( f_1 \), \( f_2 \), and \( f_4 \) with sample size $n\leq 3200$, increasing sample sizes lead to higher and faster detection by both methods. Further, for large scale datasets ($n\in 10^{5}\cdot \{1,\ldots,10\}$), FOMT and A-FOMT detect violations of all alternatives within one second. For functions \( f \) under the null model, corresponding to the most computationally expensive scenario, FOMT and A-FOMT exhibit a computational cost that grows linearly (up to logarithmic factors) with \( n \) (\cref{T:time}). This observation highlights the scalability of FOMT and A-FOMT for large datasets.


\section*{Acknowledgements}
This work was supported by DFG-FOR~5381 {Mathematical Statistics in the Information Age}, and in part by the DFG under Germany's Excellence Strategy, project EXC~2067 {Multiscale Bioimaging: from Molecular Machines to Networks of Excitable Cells} (MBExC).


\clearpage

\appendix
\section{Properties of $\gamma$-exceedance fraction}\label{s:gamma_exceedance}
We list the basic properties of $\gamma$-exceedance fractions $\varepsilon_{0,\gamma}(f)$ and $\varepsilon_{1,\gamma}(f)$ in Definition~\ref{e: exceedance 0} along with the proofs.
\begin{proposition}
\label{p:g*}
For any $\gamma>0$ and any continuous function $f:[0,1]\to \mathbb{R}$, there exists a function $g^*\in \mathcal{M}$ such that 
\begin{equation*}
    \varepsilon_{0,\gamma}(f)= \lambda\left\{x\in[0,1] \,| \,\abs{f(x)-g^*(x)} >\gamma \right\}.
\end{equation*}
\end{proposition}

\begin{proof}
Let $M := \max_{x\in[0,1]}\abs{f(x)}$. By definition, there is a function sequence  $(g_k)_{k\in \mathbb{N}}\subseteq \mathcal{M}$ such that $\lambda(D_k)\leq \varepsilon_{0,\gamma}(f)+1/k$, with $D_k= \{x\,|\, \abs{f(x)-g_k(x)} >\gamma\}$. Further, we define a modification $\Tilde{g}_k:[0,1]\to \mathbb{R}$ of $g_k$  as
\begin{equation*}
    \tilde{g}_k(x)=
    \begin{cases}
        \underset{t\in D_k^c\cap[0,x]}{\sup}g_k(t),\quad&  \text{ if } D_k^c\cap[0,x) \neq \emptyset,\\
        \underset{t\in D_k^c}{\inf}g_k(t), \quad& \text{ otherwise}.
    \end{cases}
\end{equation*}
Clearly, $\tilde{g}_k$ is monotone and  $g_k(x) = \tilde{g}_k(x)$ for all $x\in D_k^c$, which implies that 
\begin{equation}
    \label{inq: gk on Dkc}
    \abs{f(x)-\tilde{g}_k(x)} = \abs{f(x)-g_k(x)} \leq \gamma,\quad \text{ for all } x\in D_k^c.
\end{equation}
Thus, $\abs{\Tilde{g}_{k}(x)}\leq M+\gamma$ for all $x\in D_k$ and $k \in \mathbb{N}$. By Helly's Selection Principle (e.g.\ \citealp[Lemma 13.15]{carothers2000real}) there exists a subsequence of $(\Tilde{g}_k)_{k\in \mathbb{N}}$ that converges pointwise to $g^* \in \mathcal{M}$. With slight abuse of notation, we still write  $(\Tilde{g}_k)_{k\in \mathbb{N}}$ for this subsequence. 

Let $\Tilde{D}_k = \{x\,|\, \abs{f(x)-\Tilde{g}_k(x)}>\gamma\}$ and $D_*= \{x\,|\, \abs{f(x)-g^*(x)}>\gamma\}$. By \eqref{inq: gk on Dkc}, we have $D_k^c\subseteq \Tilde{D}_k^c$, and further $\lambda(\Tilde{D}_k)\leq \lambda(D_k)\leq \varepsilon_{0,\gamma}(f)+1/k$. By the pointwise convergence of $(\Tilde{g}_k)_{k\in \mathbb{N}}$, we have $D_*\subseteq \cup_{m = 1}^{\infty}A_m$ with $A_m \coloneqq \cap_{k\geq m}\tilde{D}_k$. Then $\lambda(A_m)\leq \lambda(\Tilde{D}_k) \leq \varepsilon_{0,\gamma}(f)+1/k$ for all $k\geq m$, which implies that $\lambda(A_m) \leq \varepsilon_{0,\gamma}(f)$ and further $\lambda(D_*)\leq \lambda(\cup_{m = 1}^{\infty}A_m) = \lim_{m\to \infty}\lambda(A_m) = \varepsilon_{0,\gamma}(f)$. By definition, $\lambda(D_*)\geq\varepsilon_{0,\gamma}(f)$. Thus, $\lambda(D_*)=\varepsilon_{0,\gamma}(f)$. 
\end{proof}

\begin{definition}\label{def:heavy}
Let $f:[0,1]\to \mathbb{R}$ be a continuous function. For $\gamma>0$, a point $a\in[0,1]$ is called $\gamma$-\emph{right-heavy} if there exists $b\in(a,1]$ such that $\lambda(A_{a,b})\geq (b-a)/2$, with 
\begin{equation*}
    A_{a,b} = \left\{x\in[a,b]\,|\,f(a)-f(x)\geq \gamma \right\}.
\end{equation*}
Similarly, we say $a$ is \emph{left-heavy} if there exists $b\in [0,a)$ such that $\lambda(B_{b,a}) \geq (a-b)/2$, with
\begin{equation*}
    B_{b,a} = \left\{x\in[b,a]\,|\,f(x)-f(a)\geq \gamma\right\}.
\end{equation*}
We call $a$ $\gamma$-\emph{heavy} if it is $\gamma$-right-heavy or $\gamma$-left-heavy. Let $H_f(\gamma)$, $H_{f,R}(\gamma)$ and $H_{f,L}(\gamma)$ be the sets of  $\gamma$-heavy points, $\gamma$-right-heavy points and $\gamma$-left-heavy points of $f$, respectively.
\end{definition}

\begin{lemma}
\label{l:epsilon_f}
Let $f:[0,1]\to\R$ be a continuous function and $\gamma>0$. Then:
\begin{enumerate}[label=(\text{\roman*}), wide]
    \item\label{l:eps_g:a} If $\varepsilon_{0,\gamma}(f) >0$, then there exist $0\leq a<b\leq 1$ such that $f(a)-f(b)> 2\gamma$.
    \item\label{l:eps_g:b} The sets $H_{f,R}(\gamma)$, $H_{f,L}(\gamma)$ and $H_f(\gamma)$ are measurable and $\lambda\bigl(H_f(\gamma)\bigr)\geq \varepsilon_{0,\gamma}(f)$.
\end{enumerate}
Assume further that $f\in \Sigma(\beta,L)$ with $\beta, L>0$. Then: 
\begin{enumerate}[label=(\text{\roman*}), wide]
\setcounter{enumi}{2}
    \item\label{l:eps_g:c} For $\beta\in (0,1]$, if there exist $0\leq a<b\leq 1$ such that $f(a)-f(b)= 2\gamma + \delta$ for some $\delta>0$, then 
    $\varepsilon_{0,\gamma}(f)\geq 2^{1-1/\beta}\bigl({\delta}/{L}\bigr)^{{1}/{\beta}}.$
    \item\label{l:eps_g:d} For $\beta\in (1,2]$, let $\delta \in (0,1/2]$ such that $D_{\gamma}(f)\cap [\delta,1-\delta] \neq  \emptyset$ with 
    \begin{equation*}
        D_{\gamma}(f) := \{x\in [0,1]\;|\;f'(x)\leq - \gamma\}.
    \end{equation*}
    Then $\lambda \bigl(\{a\in[\delta,1-\delta]\,|\, \underset{x\in[a-\delta,a+\delta]}{\max}f'(x)\leq -\gamma +3L\delta^{\beta-1} \}\bigr)\geq \varepsilon_{1,\gamma}(f) \wedge  (1-2\delta).$
    \item\label{l:eps_g:e} For $\beta\in (1,2]$, if $\min_{x\in (0,1)}f'(x)\leq -(\gamma+\delta)$, then 
        $\varepsilon_{1,\gamma}(f)\geq ({\delta}/{L})^{{1}/{(\beta-1)}}.$
\end{enumerate}
\end{lemma}
\begin{proof}
{\sc Part~}\ref{l:eps_g:a}.
Anticipating a contradiction, we suppose that
    \begin{equation}
        \label{ieq:lemma1}
        f(x)-f(y)\leq 2\gamma\qquad \text{for all }0\leq x<y\leq 1.
    \end{equation}
   By the uniform continuity of $f$ on $[0,1]$, there exist $0= x_0 <x_1<\dots<x_N = 1$ such that
    \begin{equation}
        \label{ieq:modulus}
        \sup_{x,y\in [x_{i-1},\;x_{i})}\abs{f(x)-f(y)}\;\le\;2\gamma\qquad \text{ for all } i = 1,\dots, N.
    \end{equation}
    Define \(m_1 \equiv m(x_0) \coloneqq \min\bigl\{f(x) \mid x \in [x_0, 1]\bigr\}\) and \(N_1 \equiv N(x_0, m_1) \coloneqq \max  \bigl\{i \mid f(x)\leq m_1+2\gamma\ \text{ for all } x \in [x_0, x_i) \bigr\} \).
    For $x \in [x_0, x_1)$, we have $ f(x) \le \inf_{y\in[x_1,1]} f(y) + 2\gamma$ by~\eqref{ieq:lemma1} and $f(x) \;\le\; \inf_{y\in[x_0,x_1)} f(y) + 2\gamma$ by \eqref{ieq:modulus}, and thus $f(x) \;\le\; m_1 + 2\gamma$. It follows that $N_1 \ge 1$.
    If $N_1 < N$, we further define
        \(m_2 \equiv m(x_{N_1}) \coloneqq \min\bigl\{f(x) \mid x \in [x_{N_1}, 1]\bigr\}\) and \(N_2 \equiv N(x_{N_1}, m_2) \coloneqq \max  \bigl\{i\mid f(x)\leq m_2+2\gamma\ \text{ for all } x \in [x_{N_1}, x_i) \bigr\}\).
    Similarly, we have $N_2 \ge N_1+1$. Thus, we can repeat this process and obtain a finite number of integers $N_0\coloneqq 0 < N_1 < \cdots < N_{k_0} = 1$ such that 
    $$
    m_i \le f(x) \le m_i + 2\gamma \qquad \text{for all } x \in [x_{N_{i-1}}, x_{N_i}), \text{ and } i = 1, \ldots, k_0.
    $$
    As $m_1 \le m_2 \le \cdots \le m_{k_0}$, we define a monotone increasing function $g: [0,1] \to \R$ as $$
    g(x) = m_i + \gamma \quad \text{ for }x\in[x_{N_{i-1}}, x_{N_i})\qquad\text{and}\qquad g(1) = m_{k_0} +\gamma. 
    $$
    It follows that $\vert f(x) - g(x) \vert \le \gamma$ for all $x \in[0,1]$. Then, $\varepsilon_{0,\gamma}(f) = 0$, which is a contradiction.
    
\medskip
{\sc Part~}\ref{l:eps_g:b}. By Definition~\ref{def:heavy} we have
    \begin{align*}
        H_{f,R}(\gamma) = \left\{ a\in[0,1]\,\bigg| \, \underset{b>a}{\max} \frac{\lambda([a,b]\cap \{x\,|\,f(x)-f(a)\geq \gamma\})}{b-a}\geq \frac{1}{2}\right\}, \\
      \text{and}\quad  H_{f,L}(\gamma) = \left\{ a\in[0,1]\,\bigg| \, \underset{b<a}{\max} \frac{\lambda([b,a]\cap \{x\,|\,f(a)-f(x)\geq \gamma\})}{a-b}\geq \frac{1}{2}\right\}. 
    \end{align*}
Define the function $h_1:[0,1]^2\to\R$ as $h_1(a, b) = \lambda([a,b]\cap \{x\,|\,f(x)-f(a)\geq \gamma\})/(b-a)$ if $a < b$ and equal to zero, otherwise. 
The function $h_1$ is continuous, and thus measurable. Further, the function $h_2:[0,1]\to \R$, defined as $h_2(a) = \max_{b>a}h_1(a, b)$ for $a \in[0,1]$, is also measurable. Note that $H_{f,R}(\gamma) = h_2^{-1}([1/2,1])$, we obtain the measurablity of $H_{f,R}(\gamma)$. Similarly, $H_{f,L}(\delta)$ is measurable, and consequently, $H_f(\delta)$ is measurable.

If $\varepsilon_{0,\gamma}(f)=0$, then clearly $\lambda\bigl(H_f(\gamma)\bigr)\geq \varepsilon_{0,\gamma}(f)$. Consider next $\varepsilon_{0,\gamma}(f)>0$. By Part~\ref{l:eps_g:a}, there exist $0\leq a<b \leq 1$ such that $f(a)-f(b)>2\gamma$. It implies $f(a)-\gamma> f(b)+\gamma$ and $[a,b]\subseteq A_{a,b}\cup B_{a,b}$. Then, at least one of $a$ and $b$ is heavy, as  $\lambda(A_{a,b})\vee\lambda(B_{a,b}) \ge \lambda(A_{a,b}\cup B_{a,b})/2\ge (b-a)/2$. Thus, \(f(a)-f(b)\leq 2\gamma\) for all \( a,b \in H_f^c(\gamma) \text{ with } a<b\).
From this, we can construct a monotone increasing function $g$ on $H_f^c(\gamma)$ such that $\vert f(x) - g(x)\vert\le \gamma$ for $x \in H_f^c(\gamma)$, in the same way as in Part~\ref{l:eps_g:a}. 
We further extend $g$ to $\Tilde{g}:[0,1]\to\mathbb{R}$ by
    \begin{equation*}
        \Tilde{g}(x) = \begin{cases}
            \underset{t\in [0,x]\cap H_f^c(\gamma)}{\sup}g(t) \quad &\text{ if }H_f^c(\gamma)\cap [0,x)\neq \emptyset,\\
            \underset{t\in H_f^c(\gamma)}{\inf}g(t)\quad &\text{otherwise}.
        \end{cases}
    \end{equation*}
    Then $\Tilde{g}$ is monotone on $[0,1]$ and $\Tilde{g}|_{H_f^c(\gamma)} = g$. Thus, 
\( \varepsilon_{0,\gamma}(f) \leq \lambda\bigl(\{x\,|\,\abs{f(x)-\tilde{g}(x)}> \gamma\}\bigr)= \lambda\bigl(\{x\in H_f(\gamma)\,|\,\abs{f(x)-\tilde{g}(x)}> \gamma\}\bigr)+\lambda\bigl(\{x\in H_f^c(\gamma)\,|\,\abs{f(x)-{g}(x)}> \gamma\}\bigr)\leq \lambda\bigl(H_f(\gamma)\bigr)\).

\medskip 
{\sc Part~}\ref{l:eps_g:c}. By Proposition~\ref{p:g*}, there is $g \in \mathcal{M}$ such that $\varepsilon_{0,\gamma}(f) = \lambda(D_g)$ with $D_g= \{x\,|\,\abs{f(x)-g(x)}>\gamma\}$. We consider three cases separately as follows. 
\begin{enumerate}[wide]
    \item[Case 1.] Both $a,b \not\in D_g$. By  monotonicity, we have
    $f(a)-\gamma \leq g(a) \leq g(b) \leq f(b)+ \gamma$, which contradicts to the condition $f(a)-f(b)=2\gamma+\delta>2\gamma$. 
    \item[Case 2.] Either $a \not\in D_g$ and $b \in D_g$, or $a \in D_g$ and $b \not\in D_g$. Consider first $a\not\in D_g$ and $b\in D_g$. Note that $ b-a\geq \left({2\gamma+\delta}/{L}\right)^{{1}/{\beta}}$ due to Hölder smoothness and $f(a)-f(b) = 2\gamma+\delta$. Then we have $(b-(\delta/L)^{1/\beta},\,b]\subseteq [a,\,b]$, and further, for any $x\in (b-(\delta/L)^{1/\beta},b]$, 
    \begin{equation*}
        g(x)-f(x)\geq g(a)-f(b)-L\abs{x-b}^{\beta}> g(a)-f(b)-\delta \geq g(a)-f(a) + 2\gamma \geq \gamma.
    \end{equation*}
    Thus, $\varepsilon_{0,\gamma}(f)\geq \lambda \bigl(\bigl(b-(\delta/L)^{{1}/{\beta}},\,b\bigr]\bigr)= \left(\delta/L\right)^{{1}/{\beta}}$.
    For $a\in D_g$ and $b\not\in D_g$, we consider $[a,\,a+(\delta/L)^{1/\beta})$ instead, and obtain the same lower bound.
    \item[Case 3.] Both $a,b\in D_g$. 
    If $f(a) < g(a) - \gamma$, then $g(x)-f(x)>\gamma$ for $x\in (b-((2\gamma+\delta)/L)^{1/\beta}),b]\subseteq [a,b]$. Thus,
    $\varepsilon_{0,\gamma}(f)\geq \lambda \bigl((b-((2\gamma+\delta)/L)^{\frac{1}{\beta}},b]\bigr)= \bigl((2\gamma+\delta)/L\bigr)^{\frac{1}{\beta}}$.
    Similarly, if $f(b) > g(b)+\gamma$, we have the same lower bound of $\varepsilon_{0,\gamma}(f)$.
    
    Consider now $g(a)<f(a)-\gamma$ and $g(b)>f(b)+\gamma$, and define 
    \begin{align*}
        a' = \inf\{x\in [a,b]\,|\, f(x)\leq g(x) + \gamma\} \quad \text{ and } \quad
        b' = \sup\{x\in [a,b]\,|\, f(x)\geq g(x) - \gamma\}.
    \end{align*}
    Then \(a' \le b'\), since otherwise for $x\in (b',a')$ it would lead to a contradiction that $g(x)-\gamma>f(x)>g(x)+\gamma$.  
    Further, we have $\abs{g(x)-f(x)}>\gamma$ for all $x\in [a,a')\cup (b',b]$, which implies $\varepsilon_{0,\gamma}(f)\geq a'-a + b -b'$.
    If $a'=b'$, then
    $\varepsilon_{0,\gamma}(f)\geq b-a \geq  \bigl((f(a)-f(b))/L\bigr)^{{1}/{\beta}} = \bigl((2\gamma+\delta)/L\bigr)^{{1}/{\beta}}$. 
    If $a' < b'$, we have \(
    \varepsilon_{0,\gamma}(f)\geq a'-a + b -b' \ge 
\bigl( {(f(a)-f(a'))_+}/{L}\bigr)^{{1}/{\beta}} + \bigl( {(f(b')-f(b))_+}/{L}\bigr)^{{1}/{\beta}}.
    \) Note further that \(g(a'+) \le g(b'-)\), \(f(a') \le g(a'+) + \gamma\) and \(f(b') \ge g(b'-) - \gamma\). Thus, 
    \begin{align*}
        \varepsilon_{0,\gamma}(f) & \geq \left( \frac{f(a)-g(a'_+)-\gamma}{L}\right)_+^{{1}/{\beta}} + \left( \frac{g(b'_-)-\gamma-f(b)}{L}\right)_+^{{1}/{\beta}}\\
      &  \geq \left( \frac{f(a)-g(a'_+)-\gamma}{L}\right)_+^{{1}/{\beta}} + \left( \frac{g(a'_+)-\gamma-f(b)}{L}\right)_+^{{1}/{\beta}}\ge2\left(\frac{\delta}{2L}\right)^{1/\beta}.
    \end{align*}
\end{enumerate}
Overall, we have \(\varepsilon_{0,\gamma}(f)\geq \left({\delta}/{L}\right)^{{1}/{\beta}} \wedge \left(({2\gamma+\delta})/{L}\right)^{{1}/{\beta}} \wedge 2\left({\delta}/({2L})\right)^{1/\beta} = 2\left({\delta}/{(2L)}\right)^{1/\beta}\).

\medskip

{\sc Part~}\ref{l:eps_g:d}. As $ \lambda\left( D_{\gamma}(f)\right) =\varepsilon_{1,\gamma}(f) $, we have
\(\lambda\left(D_{\gamma}(f)\cap[\delta,1-\delta]\right) \geq (\varepsilon_{1,\gamma}(f) - 2\delta)_+.\)
Since $D_{\gamma}(f) \cap [\delta,1-\delta] $ is non-empty and closed, by Proposition~\ref{p:delta expansion} its $\delta$-expansion satisfies \(
\lambda\bigl(D^{2\delta}_{\gamma}(f)\cap[\delta,1-\delta]\bigr)=\lambda\bigl(\left(D_{\gamma}(f)\cap[\delta,1-\delta]\right)^{2\delta}\cap[\delta,1-\delta]\bigr)\geq \bigl(\left(\varepsilon_{1,\gamma}(f) - 2\delta\right)_+ + 2\delta \bigr) \wedge (1-2\delta)\geq \varepsilon_{1,\gamma}(f) \wedge  (1-2\delta).\)
For $a\in D^{2\delta}_{\gamma}(f)\cap[\delta,1-\delta]$, by Hölder smoothness we have \[\max_{x\in [a-\delta,a+\delta]}f'(x)\leq f'(a)+L\delta^{\beta-1}\leq -\gamma +L(2\delta)^{\beta-1}+ L\delta^{\beta-1}  \leq -\gamma +3L\delta^{\beta-1}.\]
Thus, \(\bigl\{a\in [\delta,1-\delta]\,\big|\,\max_{x\in [a-\delta,a+\delta]}f'(x)\leq -\gamma +3L\delta^{\beta-1} \bigr\}\supseteq D^{2\delta}_{\gamma}(f)\cap[\delta,1-\delta]\), and then \(
\lambda\bigl(\{a\in [\delta,1-\delta]\,\big|\,\max_{x\in [a-\delta,a+\delta]}f'(x)\leq -\gamma +3L\delta^{\beta-1} \} \bigr) \geq \lambda\left(D^{2\delta}_{\gamma}(f)\cap[\delta,1-\delta]\right)\geq \varepsilon_{1,\gamma}(f) \wedge  (1-2\delta).
\)

\medskip

{\sc Part~}\ref{l:eps_g:e}. By assumption, $f'(x_0)\leq -(\gamma+\delta)$, for some $x_0 \in [0,1]$. The Hölder smoothness implies that $
f'(x)\leq f'(x_0) + L(\delta/L)^{\frac{\beta-1}{\beta-1}} \le  -\gamma$ for all $x\in A_0 \coloneqq \bigl[ 0\vee\bigl(x_0-(\delta/L)^{\frac{1}{\beta-1}}\bigr),\bigl(x_0+\left(\delta/L\right)^{\frac{1}{\beta-1}}\bigr)\wedge 1\bigr]$.
Thus, $\varepsilon_{1,\gamma}(f) \ge \lambda(A_0)\geq (\delta/L)^{1/(\beta-1)}$.
\end{proof}

We note that Proposition~\ref{p:g*}, Lemma~\ref{l:epsilon_f}~\ref{l:eps_g:a} and~\ref{l:eps_g:b} hold also for càdlàg functions, and the proofs follow similarly. Applying Lemma~\ref{l:epsilon_f}~\ref{l:eps_g:c} and~\ref{l:eps_g:e}, we obtain the following.

\begin{corollary}
\label{p:>=hn}
Let \(\Delta_{\beta,n} = \bigl({\log(n)}/{n}\bigr)^{\frac{\beta-\ceil{\beta}+1}{2\beta+1}}\) and $\gamma_n = C_{\beta} h_n^{\beta -\ceil{\beta}+1}$ 
with $C_{\beta}$ given in \eqref{e:C0} in the appendix. If there is a constant $C  \geq (2\gamma_n +2L h_n^{\beta-\ceil{\beta}+1})/\Delta_{\beta,n}$, for instance,  $C = (2C_{\beta}+2L)C_h^{\beta -\ceil{\beta}+1}$, then it holds that $\varepsilon_{\ceil{\beta}-1,\gamma_n}(f)\geq h_n$, for any $f\in \mathcal{F}_{\beta}(C\Delta_{\beta,n})$. 
\end{corollary}

\section{Selected proofs illustrating main ideas}\label{ss:proof:rate}
The proof strategy for establishing the statistical guarantees of FOMT involves disentangling the randomness arising from the algorithm and that inherent in the data. As an illustration, we present the proof of the separation rate for FOMT in the case of $\beta \in (0,1]$, while deferring the other proofs, based on similar principles but requiring additional technicalities, to Section~\ref{s:proof2}. For convenience, we restate this result from \cref{Theorem: Minimax rate optimality} below.
\begin{theorem}[Separation rate for $0< \beta \le 1$]
\label{t:consistency1} Under the model \eqref{model}, suppose that Assumptions~\hyperref[M1]{(M1)} and \ref{K1}--\ref{K3} hold with $\beta \in (0,1]$ and $L>0$. Let $(f_n)_{n\geq1}$ be a sequence of functions in $\Sigma(\beta,L)$ such that \(\varepsilon_{0,\gamma_n}(f_n) \geq \frac{2}{n}\) with 
\(  \gamma_n = C_{\beta} h_n^{\beta} \asymp \left(\log (n)/ n\right)^{{\beta}/{(2\beta+1)}}\), $C_{\beta} $ in \eqref{e:C0} and $h_n$ in \eqref{e:bw}. Set  \(C_{n}(\alpha) = -2\log\left({\alpha}/{2}\right)/ \varepsilon_{0,\gamma_n}(f_n)\) in FOMT $\Phi$ (\cref{alg:MC}). Then
\begin{equation*}
    \underset{n\to \infty}{\liminf}\,\mathbb{P}_{f_n}(\Phi = 1) \geq 1-\alpha.
\end{equation*}
\end{theorem}

\begin{proof}
By definition, FOMT (\cref{alg:MC}) accepts $H$ if and only if all local tests conducted in the algorithm return zero. Let $I_1,\dots, I_{C_{n}(\alpha)}$ denote the repeatedly generated indices via uniform sampling. For each $l \in [C_{n}(\alpha)]$, let $\mathcal{J}_l^+$ and $\mathcal{J}_l^-$ represent the sets of generated random indices $J_k$ and $J_k'$ corresponding to left and right searches initiated from $I_l$. Introduce 
\[
\mathcal{P} := \bigcup_{l \in [C_{n}(\alpha)]}\mathcal{P}_l\qquad\text{where} \quad\mathcal{P}_l := \left\{(I_l,I_l+J)\,\big|\, J\in \mathcal{J}_l^{+}\right\}\cup \left\{(I_l-J',I_l)\,\big|\, J'\in \mathcal{J}_l^{-}\right\}.
\] Note that random sets $\mathcal{P}_l$ are i.i.d.\ distributed and independent of both the random errors $\varepsilon_i$ and the signal $f_n$. Define 
\[
E := \bigcup_{1\leq i<j\leq n} E_{i,j}\qquad \text{where}  \quad E_{i,j} := \bigl\{R_{i,j}\leq -\sigma \sqrt{6W}{C_h^{-\beta-1/2}}h_n^{\beta}\bigr\},
\] 
with $R_{i,j}$ in \eqref{R_I,J}, $C_h$ in \eqref{defn Ch} and $W$ in \cref{Theorem: Bounds of variances} (in the appendix), and define also $\mathcal{I}:= \left\{(i,j)\,\big|\, D_{i,j}\geq \gamma_n -4Lh_n^{\beta}\right\}$ with $D_{i,j}$ in \eqref{ieq: upper bound Dab 2}. To disentangle the randomness in the algorithm from that in the data, we employ the following decomposition:
\begin{align}
    \mathbb{P}(\Phi = 0) & = \mathbb{P}(\Phi = 0,\,\,E)+\mathbb{P}(\Phi = 0,\,\,E^c)\nonumber\\
    &  \leq \mathbb{P}(\Phi = 0,\,\,E) + \mathbb{P}(\Phi_{i,j} = 0\text{ for all }(i,j)\in\mathcal{P},\, \mathcal{P}\cap \mathcal{I}\neq\emptyset,\,E^c) \nonumber \\
    & \qquad\qquad\qquad\qquad+ \mathbb{P}(\Phi_{i,j} = 0\text{ for all }(i,j)\in\mathcal{P},\, \mathcal{P}\cap \mathcal{I} = \emptyset,\,E^c).\label{e:decomp}
\end{align}
We will treat each term in \eqref{e:decomp} separately in the rest of the proof. 

\medskip

For the first term in \eqref{e:decomp},  using Mills's ratio and \eqref{variance bound} in \cref{Theorem: Bounds of variances}, we obtain
\begin{align*}
    \mathbb{P}(E_{i,j}) 
    \leq \exp\left(-\frac{1}{2}\cdot\frac{\left(\sigma \sqrt{6W}{C_h^{-\beta-1/2}}h_n^{\beta}\right)^2}{\mathbb{V}(R_{i,j})}\right) 
    \leq \exp\left(-\frac{1}{2}\cdot \frac{6\sigma^2 W C_h^{-2\beta-1}h_n^{2\beta}}{\sigma^2 W n^{-1}h_n^{-1}}\right)
    =n^{-3}.
\end{align*}
Then, for sufficiently large $n$, \(\mathbb{P}(E) \leq \sum_{1\leq i <j\leq n}\mathbb{P}\left( E_{i,j}\right)\leq n^2 \max_{i,j}\mathbb{P}\left( E_{i,j}\right) \leq {\alpha}/{2},\) which shows $\mathbb{P}(\Phi = 0,\,\,E) \le {\alpha}/{2}$.

\medskip

For the second term in \eqref{e:decomp}, we consider the case of $\mathcal{P}\cap \mathcal{I}\neq\emptyset$, and let $(i,j)$ be an arbitrary pair in $\mathcal{P}\cap \mathcal{I}$. Under the event $E^c$, it follows from \eqref{ieq: Cij} and \eqref{eq: D n i j} that
\begin{equation*}
    T_{i,j}  = D_{i,j}+R_{i,j} \geq \, \gamma_n -4Lh_n^{\beta} -\sigma \sqrt{6W}{C_h^{-\beta-1/2}}h_n^{\beta} \geq C_{n,\alpha,i,j} + D_{n,\beta,i,j}.
\end{equation*}
Thus, \(\mathbb{P}(\Phi_{i,j} = 0 \text{ for all }(i,j)\in \mathcal{P}\cap \mathcal{I}, E^c) = 0.\)

\medskip

For the last term in \eqref{e:decomp}, recalling that $\mathcal{P}_l$ are i.i.d. random sets, we have 
\begin{equation*}
    \mathbb{P}(\Phi_{i,j} = 0,\text{ for all }(i,j)\in\mathcal{P};\, \mathcal{P}\cap \mathcal{I} = \emptyset;\,E^c)
    \leq\, \mathbb{P}(\mathcal{P}\cap \mathcal{I} = \emptyset) 
    =\,\mathbb{P}(\mathcal{P}_l\cap \mathcal{I} = \emptyset)^{C_{n}(\alpha)}.    
\end{equation*}
Let $H_{f_n}(\gamma_n)$ be the set of $\gamma_n$-heavy points of $f_n$ (Definition~\ref{def:heavy}). 
By Lemma~\ref{l:epsilon_f} \ref{l:eps_g:b}, it holds
\begin{align*}
    \mathbb{P}(\mathcal{P}_l\cap \mathcal{I}= \emptyset) 
    =&\,\sum_{i=1}^n\mathbb{P}\bigl(\mathcal{P}_l\cap \mathcal{I}=  \emptyset \,\big| \, I_l = i\bigr)\cdot \mathbb{P}(I_l = i) \\ 
    \leq & \, \sum_{i:\, x_i \in H_{f_n}^{1/n}(\gamma_n)}\mathbb{P}\bigl(\mathcal{P}_l\cap \mathcal{I}= \emptyset\,\big| \, I_l = i\bigr) \cdot \mathbb{P}(I_l = i) + \sum_{i:\, x_i \notin H_{f_n}^{1/n}(\gamma_n)}\mathbb{P}(I_l = i)  \\
    \leq & \, \frac{1}{n}\sum_{i:\, x_i \in H_{f_n}^{1/n}(\gamma_n)}\mathbb{P}\bigl(\mathcal{P}_l\cap \mathcal{I}= \emptyset\,\big| \, I_l = i\bigr) + 1-\varepsilon_{0,\gamma_n}(f_n).
\end{align*}
If $x_i \in H_{f_n}^{1/n}(\gamma_n)$, we have $\abs{x_i-a}\leq 1/n$ for some $a\in H_{f_n}(\gamma_n) =H_{f_n,R}(\gamma_n)\cup H_{f_n,L}(\gamma_n)$. Consider first $a\in H_{f_n,R}(\gamma_n) $, i.e., there is $b\in (a,1]$ such that $\lambda(A)\geq (b-a)/2$ with 
$$A = \left\{x\in [a,b]\,\Big|\, f_n(a)-f_n(x)\geq \gamma_n \equiv C_{\beta} h_n^{\beta}\right\}.$$ 
Further, by Lemma~\ref{l: 1/8}, we obtain, for sufficiently large $n$,
\begin{equation*}
    \mathbb{P}\left(\text{ there is no }j\in \mathcal{J}_l^{+} \text{ satisfying }x_{i+j} \in A^{1/n}\right) \leq \left(\frac{7}{8}\right)^{\ceil{20\log n}}< \left(\frac{7}{8}\right)^{\frac{\log n}{\log\left({8}/{7}\right)}}= \frac{1}{n}.
\end{equation*}
For any $j\in \mathcal{J}_l^{+}$ satisfying $x_{i+j} \in A^{1/n}$, there exists $c\in A$ fulfilling $\abs{c-x_{i+j}}\leq 1/n$ and $f_n(a)-f_n(c)\geq \gamma_n$, by definition. Then, we have 
\begin{align*}
     D_{i,i+j} 
    =&\,\sum_{k = 1}^n W_{nk}(x_i)f_n(x_k) - \sum_{l = 1}^n W_{nk}(x_{i+j})f_n(x_k)\nonumber\\
    =& \sum_{k = 1}^n W_{nk}(x_i)\Big(f_n(x_k) -f_n(x_i)\Big)+f_n(x_i)-f_n(a)+f_n(a)-f_n(c)\nonumber\\
    &\quad + f_n(c)-f_n(x_{i+j}) +\sum_{k = 1}^n W_{nk}(x_{i+j})\Big(f_n(x_{i+j})-f_n(x_k)\Big)\nonumber\\
    \geq &\, -Lh_n^{\beta}-L\left(\frac{1}{n}\right)^{\beta}+\gamma_n - L\left(\frac{1}{n}\right)^{\beta} -L h_n^{\beta} \geq  \gamma_n -4Lh_n^{\beta},
\end{align*}
which shows $(i,i+j)\in \mathcal{I}$. Thus, given $I_l = i$ and the $\gamma_n$-right-heaviness of $x_i$, we have
\[\mathbb{P}(\mathcal{P}_l\cap \mathcal{I}=\emptyset\,\big| \, I_l = i) \leq \mathbb{P}\left(\text{ there is no }j\in \mathcal{J}_l^{+} \text{ satisfying }x_{i+j} \in A^{1/n}\right) \leq \frac{1}{n}.\] By symmetry, the above inequality remains valid for 
all $\gamma_n$-left-heavy point $x_i$. Then
\begin{align*}
     &\mathbb{P}(\Phi_{i,j} = 0\text{ for all }(i,j)\in\mathcal{P},\, \mathcal{P}\cap \mathcal{I} = \emptyset,\,E^c)      \leq\mathbb{P}(\mathcal{P}_l\cap \mathcal{I} = \emptyset)^{C_{n}(\alpha)} \\
     \leq\;& \,\Bigl( \frac{1}{n}\sum_{i:\, x_i \in H_{f_n}^{1/n}(\gamma_n)}\mathbb{P}_{f_n}(\mathcal{P}_l\cap \mathcal{I}= \emptyset\,\big| \, I_l = i) + 1-\varepsilon_{0,\gamma_n}(f_n)\Bigr)^{C_{n}(\alpha)}\\
    \leq\;& \,\Bigl(\frac{1}{n} + 1-\varepsilon_{0,\gamma_n}(f_n)\Bigr)^{C_{n}(\alpha)} \leq \,\Bigl(1-\frac{1}{2}\varepsilon_{0,\gamma_n}(f_n)\Bigr)^{C_{n}(\alpha)} \leq \frac{\alpha}{2},
\end{align*}
where the second last inequality follows from $\varepsilon_{0,\gamma_n}(f_n)\geq 2/n$.

\medskip

Therefore, summarizing the above calculations, we obtain
\begin{equation*}
    \limsup_{n\to \infty}\mathbb{P}_{f_n}(\Phi=0)\leq \frac{\alpha}{2}+ 0+ \frac{\alpha}{2} = \alpha,
\end{equation*}
which concludes the proof. 
\end{proof}

\section{Local polynomial estimators: a brief recap}\label{s:lpe}
Recall from Definition~\ref{d:lpe} that for any $x \in [0,1]$, the LPE $\hat{f}_n(x)$ of order $\ceil{\beta}-1$ corresponds to the first component of $\hat{\theta}_n(x)$, which represents a weighted least square estimator
\begin{subequations}
\begin{equation}
\label{theta}
\hat{\theta}_n(x) = \underset{\theta \in \mathbb{R}^{\ceil{\beta}}}{\arg \min} (-2\theta^{\top} a_{nx}+\theta^{\top}\B_{nx}\theta),
\end{equation}
where 
\begin{align}
\label{defn of Bnx}
\mathcal{B}_{nx} = & \frac{1}{nh}\sum_{i=1}^{n}U\left(\frac{x_i-x}{h}\right)U\left(\frac{x_i-x}{h}\right)^{\top}K\left(\frac{x_i-x}{h}\right), \\
a_{nx} = & \frac{1}{nh}\sum_{i=1}^{n}  U\left(\frac{x_i-x}{h}\right) K\left(\frac{x_i-x}{h}\right)Y_i, \nonumber\\
\label{defn of U}\text{and}\quad
U(x) =& \left(1,x,\frac{x^2}{2!},\dots, \frac{x^{\ceil{\beta}-1}}{(\ceil{\beta}-1)!} \right)^{\top}.
\end{align}
\end{subequations}

\noindent If $\mathcal{B}_{nx}$ is positive definite (which  holds for large enough $n$, see Remark~\ref{remark: assumptions} \ref{r:assump b}), the solution of \eqref{theta} is uniquely given by $\hat{\theta}_n (x) = \mathcal{B}_{nx}^{-1}a_{nx}$. In this case,
\begin{equation}
\label{f_n hat}
\hat{f}_{n}(x) = \sum_{i=1}^{n}W_{ni}(x)Y_{i},
\end{equation}
where
\begin{equation*}
W_{ni}(x) = \frac{1}{nh}U(0)^{\top}\mathcal{B}_{nx}^{-1}U\left(\frac{x_i-x}{h}\right)K\left(\frac{x_i-x}{h}\right).
\end{equation*}

\begin{remark}
\label{remark: assumptions}
\mbox{}
\begin{enumerate}[label=\text{\roman*}., wide]
    \item Many kernels satisfy Assumptions \ref{K1}--\ref{K3}, as detailed in \cref{Kernel and Parameters}.
    \begin{table}[htbp]
    \centering
    \caption{Common kernels and their parameters Assumptions \ref{K1}-\ref{K3}}
    \label{Kernel and Parameters}
    \begin{tabular}{| l| l| c|c| c|c|}
    \toprule
    \hline \xrowht[()]{10pt}
    Type & Kernel & $L_K$& $K_{\max}$& $\lambda_0\, (\beta\in (0,1])$ & $\lambda_0\, (\beta\in (1,2])$ \\[5pt]
    \hline \xrowht[()]{10pt}
    Epanechnikov &$K(u) = \frac{3}{4}(1-u^2)\cdot\1_{\{\abs{u}\leq 1\}}$ & $\frac{2}{3}$& $\frac{3}{4}$ &$0.5184$ & $0.0283$ \\[5pt]
    \hline \xrowht[()]{10pt}
    Triangular &$K(u) = (1-\abs{u})\cdot\1_{\{\abs{u}\leq 1\}}$ & $1$ & $1$ &$0.5250$ & $0.0276$  \\[5pt]
    \hline \xrowht[()]{10pt}
    Quartic &  $K(u) = \frac{15}{16}(1-u^2)^2\cdot\1_{\{\abs{u}\leq 1\}}$  & $\frac{5\sqrt{3}}{12}$& $\frac{15}{16}$ & $0.5234$ & $0.0228$ \\[5pt]
    \hline \xrowht[()]{10pt}
    Cosine  & $K(u) = \frac{\pi}{4}\cos\left(\frac{\pi}{2}u\right)\cdot\1_{\{\abs{u}\leq 1\}}$ &$\frac{\pi^2}{8}$ & $\frac{\pi}{4}$& $0.5194$ & $0.0276$ \\[5pt]
    \hline
    \bottomrule
    \end{tabular}
    \end{table}
    \item \label{r:assump b} Under Assumption \ref{K1}, if $h = h_n\searrow 0$ and $nh_n \to \infty$ as $n \to \infty$, then there exist constants $\lambda_0>0$ and $n_0\in \mathbb{N}$ such that 
    \begin{equation}
        \label{inequality: LP1}
        \lambda_{\min}(\B_{nx})\geq \lambda_0,\quad \text{ for all } n\geq n_0 \text{ and for all } x\in [0,1],
    \end{equation}
    where $\lambda_{\min}(\B_{nx})$ denotes the smallest eigenvalue of $\B_{nx}$ (cf.\ \citealp[Lemma~1.5]{tsybakov2008introduction}). This implies that $\B_{nx}$ is positive definite for all $x\in [0,1]$, and $\hat{f}_n(x)$ can be uniquely determined by \eqref{f_n hat} provided large $n$. Besides, \eqref{inequality: LP1} demonstrates that for all $n\geq n_0$,
    \begin{equation*}
        \norm{\B_{nx}^{-1}v}\leq \frac{1}{\lambda_0}\norm{v}, \quad \text{ for all } \,x\in [0,1] \text{ and } \,v\in \R^{\ceil{\beta}},
    \end{equation*}
    where $\norm{\cdot}$ denotes the Euclidean norm of $\R^{\ceil{\beta}}$. Thus, for all $n\geq n_0$
    \begin{equation}
    \label{B_nx^-1 leq lambda_0}
        \normop{\B_{nx}^{-1}}\leq \frac{1}{\lambda_0},
    \end{equation}
    where $\normop{\cdot}$ denotes the operator norm induced by Euclidean norm. 
    \item The choice of $\lambda_0$ as in \citet[Lemma~1.5]{tsybakov2008introduction} is suboptimal and we supply more accurate choices for $\beta \in (0,1]$ and $\beta \in (1,2]$. For a kernel $K$ that satisfies Assumptions \ref{K1}, and $\beta \in (0,1]$, we have $\B_{nx} \in \R^{1\times 1}$ and 
    \begin{align*}
        \min_{x\in [0,1/2]}\B_{nx} = \min_{x\in [0,1/2]} \frac{1}{nh}\sum_{i=1}^{n}K\left(\frac{x_i-x}{h}\right)
        \geq \underset{x\in [0,1/2]}{\min}\frac{1}{nh}\sum_{i\,:\,x_i\geq x}^{n}K\left(\frac{x_i-x}{h}\right)\to \int_{0}^{1}K(u)du = \frac{1}{2}.
    \end{align*}
    A similar result for $x\in [1/2,1]$ follows in the same way. Therefore, we derive an asymptotic uniform lower bound $\lambda_{\min}(\B_{nx})$, given by $1/2$ for $\beta \in (0,1]$. This bound is sharp by choosing $x =0$ or $1$. From this perspective, we set $\lambda_0 = 1/2$. 
    In the case of $\beta \in (1,2]$, an explicit (asymptotic) uniform lower bound of $\lambda_{\min}(\B_{nx})$ is not directly accessible. Alternatively, we consider a numerical approximation as 
    \(
    \lambda_0 \approx  \min_{x\in \mathcal{X},\; h\in H}\lambda_{\min}(\B{n_0,x}),
\)
with $n_0 = 100$, $\mathcal{X} = \{i/n_0\,|\, i \in [n_0]\}$ and $H = \{0.05, 0.10,\dots, 0.4\}$.
    See \cref{Kernel and Parameters} for details. We set
     \begin{equation}
        \label{e: lambda_0}
         \lambda_0 = \begin{cases}
             0.5,\quad& \text{ if } \beta \in (0,1],\\
             0.0228,\quad& \text{ if } \beta \in (1,2].
         \end{cases}
     \end{equation}
\end{enumerate}
\end{remark}

\begin{lemma}[{\citealp[Lemma~1.3 and Proposition~1.12]{tsybakov2008introduction}}]
\label{lemma of weights}
Suppose that $x_k \equiv k/n$ for all $k=1,\dots,n$ and Assumption \ref{K1} holds. Let $x \in [0,1]$ and $h\geq1/(2n)$. Then the weights $(W_{nk})_{k=1}^n$ of the LPE($\ceil{\beta}-1$) satisfy, for sufficiently large $n$,
\begin{enumerate}[label=(\text{\roman*}), wide]
\item $\max_{k,x} |W_{nk}(x)| \leq C_{*}/(nh)$;
\item $\sum_{k=1}^{n}|W_{nk}(x)| \leq C_{*} $;
\item $W_{nk}(x)=0$, if $|x_k-x|>h$;
\item $\sum_{k=1}^{n}W_{nk}(x)=1$ and $\sum_{k=1}^n W_{nk}(x)(x_i-x)^l = 0$ for all $l = 1,\dots, \ceil{\beta}-1$.
\end{enumerate}
where $C_{*} = 8K_{\max}/ \lambda_0$, and $\lambda_0$ in \eqref{inequality: LP1}. 
\end{lemma}

\begin{theorem}[{\citealp[Theorem~1.8]{tsybakov2008introduction}}]
\label{theorem: Minimax optimality LPE}
 Under the nonparametric Gaussian model in \eqref{model}, suppose that Assumptions~\hyperref[M1]{(M1)} and \ref{K1}--\ref{K3} hold. Let $\hat{f}_n$ be the $LPE(\ceil{\beta}-1)$ estimator of $f$ with bandwidth
\begin{subequations}\label{e:bw}
\begin{equation}
    \label{defn h_n}
    h_n = C_h\left( \frac{\log n}{n} \right)^{\frac{1}{2\beta+1}},
\end{equation}
where 
\begin{equation}
    \label{defn Ch}
    C_h = \left(\frac{q_2}{4q_1^2\beta}\right)^{\frac{1}{2\beta+1}}, \quad 
    q_1 = \frac{C_*L}{(\ceil{\beta}-1)!}, \quad
    q_2 = \sigma^2\frac{16\ceil{\beta}}{\lambda_0^2}\int_{-\infty}^{\infty}K(u)^2du,
\end{equation}
\end{subequations}
are constants independent of $n$, with $C_* = 8K_{\max}/\lambda_0$ and $\lambda_0$ in \eqref{inequality: LP1}. Then
\begin{equation*}
    \underset{n \rightarrow \infty}{\limsup} \underset{f \in \Sigma(\beta, L)}{\sup} \E_f \left( \left(\frac{\log n}{n}\right)^{-\frac{2\beta}{2\beta+1}} \norm{\hat{f}_n-f}_{\infty}^2 \right) \leq \frac{q_2}{C_h}+2q_1^2C_h^{2\beta}.
\end{equation*}
\end{theorem}

Let \( h > 0 \) and \( x \in [0,1] \). We introduce the index set of grid points that are in an $h$-neighborhood of $x$ as:
\begin{equation}
    \label{support of x}
    I_x \equiv I_x(h) = \{i \in [n] \mid |x_i - x| < h \}.
\end{equation}
For simplicity, if \( x = x_i \) for some \( i \in [n] \), we write \( I_i \) instead of \( I_x \).
Under Assumptions \ref{K1} and \ref{K2} we have $K(-1)= K(1)= 0$. Then by Lemma \ref{lemma of weights}, 
\begin{equation*}
    \hat{f}_n(x) = \sum_{i \in I_x} W_{ni}(x)Y_i,\quad \text{ for all }x\in [0,1].
\end{equation*}
\begin{proposition}[Translation]
\label{translation}
Suppose that Assumptions \ref{K1} and \ref{K2} hold. Let \( h \leq x_i < x_j \leq 1 - h \) for some \( i, j \in [n] \). Then 
\[
    W_{nk}(x_i) = W_{n(k+j-i)}(x_j), \quad \text{for all } k \in I_i.
\]
\end{proposition}

\begin{proof}
By straightforward computation, we have
\begin{equation}
    \label{equality: translation of supports}
    I_j = I_i + (j - i) \coloneqq \{k + (j - i) \mid k \in I_i\}. 
\end{equation}
Let \( \mathcal{B}_{ni} \) and \( \mathcal{B}_{nj} \) denote \( \mathcal{B}_{n x_i} \) and \( \mathcal{B}_{n x_j} \), respectively. We verify that \( \mathcal{B}_{ni}^{-1} = \mathcal{B}_{nj}^{-1} \). Using \eqref{defn of Bnx} and \eqref{equality: translation of supports}, we have:
\begin{align*}
    \mathcal{B}_{nj} &= \frac{1}{nh} \sum_{k=1}^{n} U\left(\frac{x_k - x_j}{h}\right) U\left(\frac{x_k - x_j}{h}\right)^{\top} K\left(\frac{x_k - x_j}{h}\right) \\
    &= \frac{1}{nh} \sum_{k \in I_j} U\left(\frac{x_k - x_j}{h}\right) U\left(\frac{x_k - x_j}{h}\right)^{\top} K\left(\frac{x_k - x_j}{h}\right) \\
    &= \frac{1}{nh} \sum_{k \in I_j} U\left(\frac{x_k - (x_j - x_i) - x_i}{h}\right) U\left(\frac{x_k - (x_j - x_i) - x_i}{h}\right)^{\top}  K\left(\frac{x_k - (x_j - x_i) - x_i}{h}\right) \\ 
    &= \frac{1}{nh} \sum_{k \in I_j} U\left(\frac{x_{k - (j - i)} - x_i}{h}\right) U\left(\frac{x_{k - (j - i)} - x_i}{h}\right)^{\top} K\left(\frac{x_{k - (j - i)} - x_i}{h}\right) \\
    &= \frac{1}{nh} \sum_{k \in I_i} U\left(\frac{x_k - x_i}{h}\right) U\left(\frac{x_k - x_i}{h}\right)^{\top} K\left(\frac{x_k - x_i}{h}\right) \\
    &= \frac{1}{nh} \sum_{k=1}^{n} U\left(\frac{x_k - x_i}{h}\right) U\left(\frac{x_k - x_i}{h}\right)^{\top} K\left(\frac{x_k - x_i}{h}\right) = \mathcal{B}_{ni}.
\end{align*}
By \citet[Lemma~1.5]{tsybakov2008introduction}, both \( \mathcal{B}_{ni} \) and \( \mathcal{B}_{nj} \) are non-singular under Assumption \ref{K1}, which implies \( \mathcal{B}_{ni}^{-1} = \mathcal{B}_{nj}^{-1} \). Furthermore, using \eqref{equality: translation of supports}, we obtain
\begin{align*}
    W_{nk}(x_i) &= \frac{1}{nh} U(0)^{\top} \mathcal{B}_{ni}^{-1} U\left(\frac{x_k - x_i}{h}\right) K\left(\frac{x_k - x_i}{h}\right) \\
    &= \frac{1}{nh} U(0)^{\top} \mathcal{B}_{nj}^{-1} U\left(\frac{x_k - x_i}{h}\right) K\left(\frac{x_k - x_i}{h}\right) \\
    &= \frac{1}{nh} U(0)^{\top} \mathcal{B}_{nj}^{-1} U\left(\frac{x_{k + (j - i)} - x_j}{h}\right) K\left(\frac{x_{k + (j - i)} - x_j}{h}\right) = W_{n(k + j - i)}(x_j),
\end{align*}
which finishes the proof.
\end{proof}

\begin{lemma}[Nonnegativity]
\label{l: W grt 0}
Suppose that Assumptions \hyperref[M1]{(M1)}, \ref{K1} and \ref{K3} are valid. Assume that either of the following two conditions is fulfilled:
\begin{enumerate}[label=(\text{\roman*}), wide]
    \item $\beta\in (0,1]$ and $x\in [0,1]$, or
    \item $\beta\in (1,2]$ and $x= x_i\in [h,1-h]$ for some $i\in \{1,\dots,n\}$. 
\end{enumerate}
Then it holds $$W_{nk}(x) = \frac{K\left(\frac{x_k-x}{h}\right)}{\sum_{j=1}^n K\left(\frac{x_j-x}{h}\right) }\geq 0,\quad \text{for all }k\in I_x.$$
\end{lemma}

\begin{proof}
For $\beta \in (0,1]$, the local polynomial estimator reduces to the \emph{Nadaraya–Watson} estimator, where the weight is given by
\begin{equation*}
    W_{nk}(x) = \frac{K\left(\frac{x_k-x}{h}\right)}{\sum_{j=1}^n K\left(\frac{x_j-x}{h}\right) }\geq 0,
\end{equation*}
and the inequality follows from Assumption \ref{K1}.

Now assume that $\beta \in (1,2]$. By straightforward computation we obtain that 
\begin{align*}
     W_{nk}(x_i) = \frac{1}{\alpha_{n,0}\alpha_{n,2}-\alpha_{n,1}^2}\left(\alpha_{n,2} -\alpha_{n,1}\cdot\frac{x_k-x_i}{h}\right)\cdot K\left(\frac{x_k-x_i}{h}\right),
\end{align*}
where   
\begin{align*}
    \alpha_{n,0} &\coloneqq \sum_{j=1}^n K\left(\frac{x_j-x_i}{h}\right), \\
    \alpha_{n,1} &\coloneqq \sum_{j=1}^n \frac{x_j-x_i}{h}K\left(\frac{x_j-x_i}{h}\right),\\
    \text{and}\quad\alpha_{n,2} &\coloneqq \sum_{j=1}^n \left(\frac{x_j-x_i}{h}\right)^2 K\left(\frac{x_j-x_i}{h}\right).
\end{align*}
Note that when $x_i \in [h,1-h]$, the terms $\alpha_{n,0}$, $\alpha_{n,1}$ and $\alpha_{n,2}$ are independent of $i$. From Assumption \ref{K3} and the condition $x_i \in [h,1-h]$, it follows that $\alpha_{n,1}= 0$. Therefore, applying Assumption \ref{K1}, we have
\begin{equation*}
    W_{nk}(x_i) = \frac{K\left(\frac{x_k-x_i}{h}\right)}{\sum_{j=1}^n K\left(\frac{x_j-x_i}{h}\right) }\geq 0.
\end{equation*}
\end{proof}
\begin{remark}
For $x = x_i \in [h,1-h]$ with $i\in [n]$, the \emph{local linear estimator} (LPE($1$)), coincides with the \emph{Nadaraya–Watson} estimator (LPE($0$)).
\end{remark}

\begin{proposition}\label{p:upper bound Dab}
Suppose that Assumptions \hyperref[M1]{(M1)}, \ref{K1} and \ref{K3} hold. Let $f \in H \equiv \mathcal{M} \cap \Sigma(\beta, L)$ and $0 \leq a < b \leq 1$ and $n$ be sufficiently large. Then 
\begin{equation}
    \label{ieq: upper bound Dab 1}
    D_{a,b} \coloneqq \sum_{k=1}^{n}W_{nk}(a)f(x_k) - \sum_{k=1}^{n} W_{nk}(b)f(x_k)\leq \left(\frac{16K_{\max}}{\lambda_0}\vee 2\right)Lh^{\beta}.
\end{equation}
If additionally $a = x_i$ and $b = x_j$ for some $i, j \in [n]$ with $x_i, x_j \in [h, 1-h]$, then
\begin{equation*}
    \label{ieq: upper bound Dab 2}
    D_{i,j} \coloneqq D_{x_i,x_j} \leq 0.
\end{equation*}
More precisely, we have
\begin{equation}
    \label{eq: D n i j}
        D_{i,j}\leq D_{n,\beta,i,j} \equiv D_{n,\beta,i,j}(h) \coloneqq \begin{cases}
            0,&\quad \text{ if } x_i, x_j \in [h, 1-h],\\
            \left(\frac{16K_{\max}}{\lambda_0}\vee 2\right)Lh^{\beta},&\quad \text{otherwise,}
        \end{cases}
\end{equation}
i.e., $D_{n,\beta,i,j}$ serves as an upper bound for $D_{i,j}$ for all $f \in H$.
\end{proposition}
\begin{proof}
Let $0 \leq a < b \leq 1$. Note that $\sum_{k=1}^{n} W_{nk}(x) = 1$ for all $x \in [0, 1]$, so we have, for any $f\in H$,
\begin{align}
    \label{ieq: upper bound Dab 3}
    D_{a,b} =& \sum_{k=1}^{n}W_{nk}(a)f(x_k) - \sum_{k=1}^{n} W_{nk}(b)f(x_k) \nonumber\\
    =& \sum_{k=1}^{n}W_{nk}(a)(f(x_k) - f(a)) + f(a) - f(b) + \sum_{k=1}^{n} W_{nk}(b)(f(b) - f(x_k)) \nonumber\\
    \leq& \sum_{k=1}^{n}W_{nk}(a)(f(x_k) - f(a)) + \sum_{k=1}^{n} W_{nk}(b)(f(b) - f(x_k)).
\end{align}
If $\beta \in (0,1]$, then by Lemmas \ref{lemma of weights} and \ref{l: W grt 0} and the Hölder smoothness of $f$,
\begin{equation*}
    \sum_{k=1}^{n}W_{nk}(a)(f(x_k) - f(a)) \leq \sum_{k=1}^{n}W_{nk}(a) L\abs{x_k-a}^{\beta}\leq Lh^{\beta}.
\end{equation*}
If $\beta \in (1,2]$, then by Lemma \ref{lemma of weights}, 
\begin{align}
    \label{ieq: beta <= 2}
    \sum_{k=1}^{n}W_{nk}(a)(f(x_k) - f(a)) =& \sum_{k=1}^{n}W_{nk}(a)(f'(\xi_k) - f'(a))(x_k-a) \nonumber \\
    \leq& \sum_{k=1}^{n}\abs{W_{nk}(a)} \cdot \abs{f'(\xi_k) - f'(a)} \cdot \abs{x_k-a} 
    \leq \frac{8K_{\max}}{\lambda_0}Lh^{\beta}.
\end{align}
Combining \eqref{ieq: upper bound Dab 3}--\eqref{ieq: beta <= 2} gives \eqref{ieq: upper bound Dab 1}.
If additionally $a = x_i$ and $b = x_j$ for some $i, j \in [n]$ with $h \leq x_i < x_j \leq 1-h$, then by Proposition \ref{translation} and Lemma \ref{l: W grt 0}, we have
\begin{align*}
    D_{i,j} &= \sum_{k=1}^{n}W_{nk}(x_i)f(x_k) - \sum_{k=1}^{n} W_{nk}(x_j)f(x_k) \notag \\
    &= \sum_{k \in I_i} W_{nk}(x_i)f(x_k) - \sum_{k \in I_j} W_{n(k-(j-i))}(x_i)f(x_k) \notag \\
    &= \sum_{k \in I_i} W_{nk}(x_i)(f(x_k) - f(x_{k+j-i})) \leq 0. 
\end{align*}
\end{proof}

\begin{theorem}
\label{Theorem: Bounds of variances}
Suppose that $K$ is a kernel satisfying Assumptions \ref{K1} and \ref{K2}. For any $a,b\in [0,1]$, the following inequality holds:
\begin{equation}
    \label{inequality: Bounds of variances}
    \sum_{k=1}^{n}\Big(W_{nk}(a)-W_{nk}(b)\Big)^2 \leq W n^{-1}h^{-3}\left(\frac{8L_2}{\lambda_0}|a-b| \wedge h\right)^2,
\end{equation}
where $W= 4C_*^2 \vee 4\widetilde{C}_*^2$, $\widetilde{C}_* = L_1/(8L_2)+\sqrt{e}K_{\max}/\lambda_0$, $L_1 = \sqrt{e}L_K + \sqrt{e} K_{\max}$ and $L_2 =~ 2\ceil{\beta}K_{\max} + \ceil{\beta}L_K$,
with $\lambda_0$ in \eqref{inequality: LP1}, $C_*$ in Lemma \ref{lemma of weights}, and $K_{\max}$ and $L_K$ in Assumptions \ref{K1} and \ref{K2}, respectively.
In particular, under the assumption of normality, the variance of $R_{a,b}\coloneqq \sum_{k=1}^{n}W_{nk}(a)\varepsilon_k - \sum_{k=1}^{n} W_{nk}(b)\varepsilon_k$ satisfies
\begin{equation}
    \label{variance bound}
    \V(R_{a,b}) \leq \sigma^2 W n^{-1}h^{-3}\left(\frac{8L_2}{\lambda_0}|a-b| \wedge h\right)^2.
\end{equation}
\end{theorem}

\begin{proof}
For any $a,b\in [0,1]$, we have, by Lemma \ref{lemma of weights},
\begin{align*}
    \sum_{k=1}^{n}(W_{nk}(a)-W_{nk}(b))^2 \leq 2\sum_{k=1}^n W_{nk}^2(a) +2\sum_{k=1}^n W_{nk}^2(b) 
    \leq 2\underset{k,x}{\max} \abs{W_{nk}(x)}\left(\sum_{k=1}^n \abs{W_{nk}(a)}+\sum_{k=1}^n \abs{W_{nk}(b)} \right)\leq \frac{4C_*^2}{nh}.
\end{align*}
Now we prove the statement for $a,b$ with $\abs{a-b}/h \leq \lambda_0 / (8L_2)$ by verifying that $\B_{n\cdot}$ in \eqref{defn of Bnx}, $\B_{n\cdot}^{-1}$, and $W_{nk}(\cdot)$ are all Lipschitz continuous, in the following steps:
\begin{enumerate}
    \item[Step 1:]  For $\B_{n\cdot}$, define $G:\mathbb{R}\to \mathbb{R}^{\ceil{\beta}\times \ceil{\beta}}$ by $G(u) = U(u)U(u)^{\top}K(u)$ for all $u \in \R$, where $U$ is given by \eqref{defn of U} and $K$ is a kernel.
Then,
\begin{align}
    \normop{\B_{na}-\B_{nb}} 
    &\leq \frac{1}{nh} \sum_{k=1}^{n}\normop{G\left(\frac{x_k-a}{h}\right) - G\left(\frac{x_k-b}{h}\right)  } \nonumber\\
    &= \frac{1}{nh} \sum_{k \in I_a \cup I_b} \normop{G\left(\frac{x_k-a}{h}\right)-G\left(\frac{x_k-b}{h}\right)}
    \label{inequality: B_nx Lipschitz}
    \leq \frac{1}{nh} \sum_{i \in I_a \cup I_b}L_2 \frac{|a-b|}{h} \leq 4 L_2 \frac{|a-b|}{h},
\end{align}
where $I_a$ and $I_b$ are given in \eqref{support of x}, and the second last inequality follows from the Lipschitz continuity of $G$ with parameter $L_2 = 2\ceil{\beta}K_{\max}+\ceil{\beta}L_{K}$ (cf.\ Lemma~\ref{Lipschitz Lemma}).
\item[Step 2:] For $\B_{n\cdot}^{-1}$, let $r =  \normop{\B_{nb}^{-1}(\B_{na}-\B_{nb})}$. Note that $\abs{a-b}/h\leq \lambda_0 /(8L_2)$ and it follows from \eqref{B_nx^-1 leq lambda_0} and \eqref{inequality: B_nx Lipschitz} that 
\begin{align*}
    r \leq  \normop{\B_{nb}^{-1}} \normop{\B_{na}-\B_{nb}}  
    \leq \normop{\B_{nb}^{-1}} 4L_2 \frac{|a-b|}{h} 
    \leq \frac{4L_2}{\lambda_0} \frac{|a-b|}{h} \leq \frac{1}{2}.
\end{align*}
Applying \citet[Theorem~2.3.4]{golub1983matrix} 
with $r \in (0,1/2]$ we have
\begin{equation}
    \label{inequality: B_nx^-1 Lipschitz}
    \normop{\B_{na}^{-1}-\B_{nb}^{-1}}\leq \normop{\B_{na}-\B_{nb}}\norm{\B_{nb}^{-1}}^2\frac{1}{1-r} \leq \frac{4L_2}{\lambda_0^2}\frac{|a-b|}{(1-r)h}  \leq \frac{8L_2}{\lambda_0^2}\frac{|a-b|}{h}.
\end{equation}
\item[Step 3:] For $W_{nk}(\cdot)$, define the function $F$ as $F(u) = U(u)K(u)$ for all $u \in \R$.
Note that $\norm{U(0)^{\top}}=1$. We obtain
\begin{align*}
     &\abs{W_{nk}(a)-W_{nk}(b)} \\
    \leq&\frac{1}{nh} \norm{U(0)} \cdot \norm{\B_{na}^{-1}F\left(\frac{x_k-a}{h}\right)-\B_{nb}^{-1}F\left(\frac{x_k-b}{h}\right)}\\
    \leq& \frac{1}{nh}\left(\norm{\B_{na}^{-1}F\left(\frac{x_k-a}{h}\right)-\B_{na}^{-1}F\left(\frac{x_k-b}{h}\right)} +\norm{\B_{na}^{-1}F\left(\frac{x_k-b}{h}\right)-\B_{nb}^{-1}F\left(\frac{x_k-b}{h}\right)}  \right) \\
    \leq& \frac{1}{nh}\left(\normop{\B_{na}^{-1}}\cdot \norm{F\left(\frac{x_k-a}{h}\right)-F\left(\frac{x_k-b}{h}\right)}+ \normop{\B_{na}^{-1}-\B_{nb}^{-1}} \cdot\norm{F\left(\frac{x_k-b}{h}\right)}\right).
\end{align*}
It follows from \eqref{B_nx^-1 leq lambda_0} and Lemma \ref{Lipschitz Lemma} that
\begin{equation*}
     \normop{\B_{na}^{-1}}\cdot \norm{F\left(\frac{x_k-a}{h} \right)-F\left( \frac{x_k-b}{h}\right)} \leq \frac{L_1}{\lambda_0}  \frac{|a-b|}{h}.
\end{equation*}
Further, since $\norm{U(u)}\leq \sqrt{e}$ and $|K(u)|\leq K_{\max}$ for all $u\in [-1,1]$, we obtain 
$$\norm{F(u)} \leq \sqrt{e}K_{\max},\quad \text{for all } u\in \R.$$
Together with \eqref{inequality: B_nx^-1 Lipschitz} we have the upper bound
\begin{equation*}
    \normop{\B_{na}^{-1}-\B_{nb}^{-1}} \norm{F\left( \frac{x_k-b}{h} \right)} \leq \sqrt{e}K_{\max}\frac{8L_2}{\lambda_0^2}\frac{|a-b|}{h}. 
\end{equation*}
Thus, $\abs{W_{nk}(a)-W_{nk}(b)} \leq \left(L_1/ \lambda_0  +8L_2\sqrt{e}K_{\max}/\lambda_0^2\right) \cdot |a-b| /(nh^2)$.
\item[Step 4:] Since $W_{nk}$ is Lipschitz continuous, we have
\begin{align*}
    \sum_{k=1}^{n} (W_{nk}(a)-W_{nk}(b))^2 = \sum_{k\in I_a \cup I_b} |W_{nk}(a)-W_{nk}(b)|^2 \leq  4\left(\frac{L_1}{\lambda_0}  +\sqrt{e}K_{\max}\frac{8L_2}{\lambda_0^2}\right)^2\frac{|a-b|^2}{nh^3}.
\end{align*}
\end{enumerate}
This finishes the proof of \eqref{inequality: Bounds of variances}, which implies \eqref{variance bound} directly.
\end{proof}

By Mill's ratio and Theorem~\ref{Theorem: Bounds of variances}, we obtain the following result.
\begin{theorem}
\label{general tail distribution of R_ab}
Under the nonparametric regression model in \eqref{model}, we assume that the conditions in Theorem \ref{Theorem: Bounds of variances} hold. Then for any $\alpha \in (0,1)$,
\begin{equation*}
    \P\Big(R_{a,b} \geq C_{n,\alpha}(a,b)\Big) \leq \alpha,
\end{equation*}
where 
\begin{equation*}
    C_{n,\alpha}(a,b) \coloneqq \sigma\sqrt{-2\log(\alpha)W}n^{-\frac{1}{2}}h^{-\frac{3}{2}}\left(\frac{8L_2}{\lambda_0}|a-b|\wedge h\right).
\end{equation*}
\end{theorem}

Recall that the product of Lipschitz functions is Lipschitz continuous. In particular, if $K:[-1,1]\to \R$ is Lipschitz continuous with $L_K>0$ and $U:[-1,1]\to \R^{m\times n}$ is Lipschitz continuous with $L_U>0$, i.e., $\normop{U(u)-U(u')}\leq L_U\abs{u-u'}$ for all $u,u'\in[-1,1]$. Then $K\cdot U$ is Lipschitz continuous, such that 
\begin{equation}
    \label{ieq:KU Lip}
    \normop{K(u)U(u)-K(u')U(u')}\leq (U_{\max}L_K + K_{\max} L_U)\abs{u-u'}
\end{equation}
for all $u,u'\in [-1,1]$, where $K_{\max} = \underset{u\in [-1,1]}{\max}\abs{K(u)}$, $U_{\max} = \underset{u\in [-1,1]}{\max}\norm{U(u)}$, and $\normop{\cdot}$ stands for the operator norm.

\begin{lemma}
\label{Lipschitz Lemma}
Let $\beta>0$ and $K$ be a kernel satisfying Assumptions \ref{K1} and \ref{K2}. Define the functions $F$ and $G$ as $F(u) = U(u)K(u)$ and $G(u) = U(u)U(u)^{\top}K(u)$ for all $u\in \R$ with
$U$ in \eqref{defn of U}. Then $F$ and $G$ are Lipschitz continuous with parameters 
\begin{equation*}
    L_1 = \sqrt{e}(K_{\max}+L_K) \quad \text{ and } \quad L_2 = 2\ceil{\beta}K_{\max}+\ceil{\beta}L_{K}, \text{ respectively.}
\end{equation*}
\end{lemma}

\begin{proof}
Let $V(u) = U(u)U(u)^{\top}$. For $u,v\in [-1,1]$ with $u<v$, by Lagrange's mean value theorem we have 
\begin{align*}
    \norm{U(u)-U(v)} 
    =& \norm{\left(0,u-v, (u-v)\xi_1,\dots, \frac{(u-v)\xi_{\ceil{\beta}-2}^{\ceil{\beta}-2}}{(\ceil{\beta}-2)!}\right)}\\
    \leq & \abs{u-v} \sqrt{1+1+\frac{1}{2!}+\cdots+\frac{1}{(\ceil{\beta}-2)!}}
    <\abs{u-v}\sqrt{e}, 
\end{align*}
for some $\xi_i \in(u,v)$ for all $i =1,\dots, \ceil{\beta}-2$, i.e., $U$ is Lipschitz continuous with $L_U = \sqrt{e}$.
Similarly, for any $u\in[-1,1]$, 
\begin{align*}
    \norm{U(u)}
    \leq \sqrt{1+1+\frac{1}{2!}+\cdots+ \frac{1}{\floor{\beta}!}}
    < \sqrt{e}\eqqcolon U_{\max}. 
\end{align*}
By \eqref{ieq:KU Lip}, we obtain the Lipschitz continuity of $F$ with $L_1 = \sqrt{e}(K_{\max}+L_K)$.

Applying again Lagrange's mean value theorem,  we have 
\begin{align*}
    \abs{(V(u)-V(v))_{i,j}} 
    = \abs{u-v} \left|\frac{(i+j-2)\xi_{i,j}^{i+j-3}}{(i-1)!(j-1)!}\right| \leq 2\abs{u-v},
\end{align*}
for some $\xi_{i,j}\in[-1,1]$. It follows 
\begin{equation*}
    \normop{V(u)-V(v)} \leq \normf{V(u)-V(v)} 
    = \sqrt{\sum_{i=1}^{\ceil{\beta}}\sum_{j=1}^{\ceil{\beta}} (V(u)-V(v))^2_{i,j}}
    \leq 2\ceil{\beta}\cdot \abs{u-v},
\end{equation*}
which indicates that $V$ is Lipschitz continuous 
with $L_V = 2\ceil{\beta}$.
Note that $\abs{V(u)_{i,j}}\leq 1$ for all $i,j\in [\ceil{\beta}]$ and all $u \in [-1,1]$. It implies that 
\begin{align*}
    \normop{V(u)}\leq  \normf{V(u)} 
    \leq \sqrt{\sum_{i=1}^{\ceil{\beta}} \sum_{j=1}^{\ceil{\beta}} 1} = \ceil{\beta}.
\end{align*}
Thus, by \eqref{ieq:KU Lip}, we obtain that $G$ is Lipschitz continuous with $L_2 = 2\ceil{\beta}K_{\max}+\ceil{\beta}L_{K}$.
\end{proof}

\section{An auxiliary result from measure theory}
\begin{definition}
\label{e: delta expansion}
Let $A\subseteq\R$, $A \neq \emptyset$  and $\delta>0$. The $\delta$-expansion of $A$ is defined as 
\begin{equation*}
    A^{\delta} = \{x\in\R\,|\,\text{ there exists } a\in A \text{ with } \abs{x-a}\leq \delta\}.
\end{equation*}
\end{definition}

\begin{proposition}
\label{p:delta expansion}
Suppose that $-\infty < m \leq M < \infty$. Let $A$ be a non-empty and closed subset of $[m,M]$ and $\delta>0$. Then 
\begin{equation*}
    \lambda\left(A^{\delta}\cap[m,M]\right)\geq \left(\lambda\left(A\right)+\delta\right)\wedge (M-m).
\end{equation*}
\end{proposition}
\begin{proof}
If $\delta > M-m-\lambda(A)$, then $A^{\delta}\cap[m,M] = [m,M]$ and $\lambda\left(A^{\delta}\cap[m,M]\right) =M-m$. Next we consider $\delta\leq M-m-\lambda(A)$, and define $a = \min A$ and $b = \max A$. If $a\geq m+\delta$, then $B \coloneqq [a-\delta,a) \subseteq  \left(A^{\delta} \cap [m,M]\right) \backslash A$ with $\lambda(B)=\delta$. If $b\leq M-\delta$, we have similarly $B \coloneqq (b,b+\delta]\subseteq  \left(A^{\delta} \cap [m,M]\right) \backslash A$ with $\lambda(B)=\delta$. Now we suppose that $a<m+\delta$ and $b>M-\delta$. Then
\begin{equation*}
    C \coloneqq [a,b] \backslash A = (a,b)\cap A^c
\end{equation*}
is a open subset of $(a,b)$.  Thus, by \citet[Theorem~4.6]{carothers2000real}, the set $C$ can be written as a countable union of disjoint open intervals, i.e., $C=\bigcup_{i=1}^{\infty}(a_i,b_i)$ and $(a_i,b_i)\cap (a_j,b_j) = \emptyset$ for $i\neq j$. Clearly, for all $i\in \mathbb{N}$, $a\leq a_i\leq b_i \leq b$,  $a_i,b_i \in A$, and thus
\begin{equation}
    \label{eq: lambda A}
    \lambda(A) = b-a -\lambda(C) = b-a -\sum_{i=1}^{\infty} (b_i-a_i).
\end{equation}
If there exists $i$ with $b_i-a_i \geq \delta$, then $B \coloneqq (a_i,a_i+\delta)\subseteq \left( A^{\delta} \cap [m,M] \right) \backslash A$ with $\lambda(B)= \delta$. If $b_i-a_i<\delta$ for all $i\in \mathbb{N}$, we obtain $(a_i,b_i)\subseteq(a_i,a_i+\delta)\subseteq (A^{\delta}\cap [m.M])\backslash A$, and $[m,a)\cup (b,M]\subseteq \left(A^{\delta}\cap[m,M]\right)\backslash A$. Let $B = \cup_{i = 1}^{\infty}(a_i,b_i) \cup [m,a)\cup (b,M]$, then by \eqref{eq: lambda A},
\begin{align*}
    \lambda(B) = \sum_{i = 1}^{\infty}(b_i-a_i) + a-m + M-b 
    = M-m-\lambda(A) \geq \delta.
\end{align*}
Summerizing all cases, we can always find a set $B$ with Lebesgue measure $\lambda(B)\geq \delta$ and $B \subseteq \left(A^{\delta}\cap[m,M]\right)\backslash A$. Therefore, 
\begin{equation*}
    \lambda \left( A^{\delta} \cap [m,M] \right) \geq \lambda \left( A \cup B \right) = \lambda\left( A\right) + \lambda \left( B \right) \geq \lambda\left( A \right) + \delta.
\end{equation*}
\end{proof}


\section{On the average deviation of $\chi^2$-exponentials}

\begin{lemma}
\label{Lemma: 6.2}
Let $\Gamma_1,\Gamma_2,\dots$ be independent, standard Gaussian distributed random variables. If $\nu_m = \sqrt{2\log m}(1-\varepsilon_m)$ with $\lim_{m\to \infty} \varepsilon_m = 0$ and $\lim_{m \to \infty} \sqrt{\log m}\varepsilon_m = \infty$, then 
\begin{equation*}
    \lim_{m\to \infty} \E \left( \left( \frac{1}{m}\sum_{i=1}^{m} \left|\exp{\left(\nu_m \Gamma_i -\frac{\nu_m^2}{2}\right)}-1\right| \right) \wedge 1\right) = 0.
\end{equation*}
\end{lemma}
\begin{proof}
Let $Z_m = \exp{(\nu_m \Gamma_i - \nu_m^2/2)}$. Then 
\[
\E(Z_m) = \int \exp{\left(\nu_m x - \frac{\nu_m^2}{2} \right)} \frac{1}{\sqrt{2\pi}} \exp{\left(-\frac{x^2}{2}\right)} dx = 1.
\]
For any $\eta \geq 1/m$ and any $\delta>0$ we have
\[
   \E\left(Z_m \1_{\{Z_m \geq \eta m\}}\right) \leq \E \left(Z_m^{1+\delta}(\eta m)^{-\delta}\right)
    = \exp{\left(\frac{\delta(1+\delta)\nu_m^2}{2}-\delta  \log(\eta m)\right)}.
\]
By choosing $\delta = \varepsilon_m$, the latter bound becomes
\begin{equation*}
    \exp{\left( -(\varepsilon_m^2+O(\varepsilon_m^3)\log m) +o(1)\right)} \to 0 \quad \text{ as } m \to \infty.
\end{equation*}
Then by \citet[Theorem~4.12]{petrov1995limit} we have
\begin{equation*}
    \frac{1}{m}\sum_{i=1}^{m} \left|\exp{\left(\nu_m \Gamma_i - \frac{\nu_m^2}{2}\right)}-1\right| \overset{P}{\rightarrow} 0,
\end{equation*}
which is equivalent to 
\begin{equation*}
    \lim_{m\to \infty} \E \left( \left( \frac{1}{m}\sum_{i=1}^{m} \left|\exp{\left(\nu_m \Gamma_i - \frac{\nu_m^2}{2}\right)}-1\right| \right) \wedge 1\right) = 0.
\end{equation*}
\end{proof}
\begin{remark}
Lemma~\ref{Lemma: 6.2} is a modification of \citet[Lemma~6.2]{dumbgen2001multiscale}.
\end{remark}


\section{Parameter choices of FOMT}
In this section, we provide precise expressions of quantities introduced in of \cref{S: The testing procedure}. Recall from \eqref{e:bw} that the minimax optimal bandwidth $h_n$ is defined as
\begin{equation*}
    h_n = C_h \left( \frac{\log n}{n}\right)^{\frac{1}{2\beta+1}}, 
\end{equation*}
with $C_h$ given in \eqref{defn Ch}. For $1\leq i<j\leq n$, the test statistic $T_{i,j}$ in \eqref{T_I,J} is decomposed into a deterministic term $D_{i,j}$ and a random term $R_{i,j}$ as follows:
\begin{equation}
    D_{i,j} = \sum_{k=1}^{n} \bigl(W_{nk}(x_i)-W_{nk}(x_{j})\bigr)f(x_k),\quad
    \label{R_I,J}
    R_{i,j} = \sum_{k=1}^{n} \bigl(W_{nk}(x_i)-W_{nk}(x_{j})\bigr)\varepsilon_k.    
\end{equation} 
For each local test $\Phi_{i,j}$, the corresponding critical value $q_{n,\beta,i,j}(\alpha)$ is defined as
\begin{subequations}
\begin{equation}
        \label{eq: critical}
        q_{n,\beta,i,j}(\alpha) = C_{n,\alpha,i,j} + D_{n,\beta,i,j},
\end{equation}
where
    \begin{align}
        \label{C,n,alpha,I,J}
        C_{n,\alpha,i,j} &= \sigma\sqrt{-2\log\left(\frac{\alpha}{N_{\max}}\right)W}\cdot n^{-\frac{1}{2}}h_n^{-\frac{3}{2}}\left(\frac{8L_2}{\lambda_0}|x_i-x_j|\wedge h_n\right),\\
        D_{n,\beta,i,j} &=\begin{cases}
            0,&\quad \text{ if } x_i,x_j\in [h_n,1-h_n],\\
            \left(\frac{16K_{\max}}{\lambda_0}\vee2\right)Lh_n^{\beta},&\quad \text{otherwise}, 
        \end{cases}\nonumber \\
        \label{inequality: N< N_max}
      \text{and}\quad  N_{\max}&\coloneqq  \frac{40}{\log(2)} C_{n}(\alpha)\log^2 n,
    \end{align}
    with a to-be-determined parameter $C_{n}(\alpha)$. The quantity $N_{\max}$ is an upper bound of the number of local tests conducted in FOMT.
\end{subequations}
Finally, we define constant $C_{\beta}$ by
\begin{equation}
    \label{e:C0}
    C_{\beta}= 
    \begin{cases}
        \left(\left(\frac{16K_{\max}}{\lambda_0}\vee 2\right)+4\right)L+ \sigma \left(\sqrt{\frac{4W}{C_h^{2\beta+1}}} + \sqrt{\frac{6W}{C_h^{2\beta+1}}}\right),&\quad \beta \in(0,1],\\
        4L+\frac{16L_2\sigma}{\lambda_0}\sqrt{\frac{4W}{C_h^{2\beta+1}}},&\quad \beta \in(1,2],
    \end{cases}\\
\end{equation}
with $W$ in \cref{Theorem: Bounds of variances} and $C_h$ \eqref{defn Ch}, respectively.

\algdef{SE}[REPEATN]{RepeatN}{End}[1]{\algorithmicrepeat\ #1 \textbf{times}}{\algorithmicend}
\alglanguage{pseudocode}
\begin{algorithm}[h]
\begin{algorithmic}[1]
\Statex {\textbf{Input:} data $Y_1,\dots,Y_n$, and significance level $\alpha \in(0,1)$}
\Statex {\textbf{Parameters:} standard deviation $\sigma$, kernel function $K$, smoothness order $\beta$ and radius $L$}
\For{$1 \leq m \leq C_n(\alpha)$}
    \State Generate $I \sim \mathrm{Unif}([n])$\;
    \If{$I \leq n-1$}
        \If{$\Phi_{I,I+1}=1$}
            \Return $\Phi_S=1$\;
        \EndIf 
    \EndIf
    \If{$I \geq 2$}
        \If{$\Phi_{I-1,I}=1$}
            \Return $\Phi_S=1$\;
        \EndIf
    \EndIf
\EndFor
\State \Return $\Phi_S=0$\;
\end{algorithmic}
\caption{S-FOMT: Simplified-FOMT \texorpdfstring{$\Phi_S$}{tilde Phi}}
\label{alg:MC2}
\end{algorithm}

\section{Further proofs}\label{s:proof2}
\subsection{Remaining proofs for \texorpdfstring{\cref{S: The testing procedure}}{S:testing}}
\begin{proof}[Proof of Theorem \ref{t:type I error}]
For any $1\leq i<j\leq n$ and $f\in H$, it follows from Proposition~\ref{p:upper bound Dab} and Theorem~\ref{Theorem: Bounds of variances} that
\begin{align*}
    \P_{f}(\Phi_{i,j}=1) =\, \P_{f}(D_{i,j}+R_{i,j}\geq C_{n,\alpha,i,j}+D_{n,\beta,i,j}) 
    \leq \,\P_{f}(R_{i,j}\geq C_{n,\alpha,i,j})
    \leq \, \exp\left(-\frac{C_{n,\alpha,i,j}^2}{2\mathbb{V}(R_{i,j})}\right)
    =\, \frac{\alpha}{N_{\max}}.
\end{align*}Namely, each local test $\Phi_{i,j}$ is of $(\alpha/N_{\max})$-level.
Further,
\begin{align*}
    \P_{f}(\Phi=1) 
    =&  \P_{f}(\text{ there exists } l \in [N_{\max}] \text{ such that } \Phi_{I_l,J_l}= 1) \\
    \leq & N_{\max} \max_{l \in [N_{\max}]} \P_{f}(\Phi_{I_l,J_l}=1) \\
    =&   N_{\max}\max_{l \in [N_{\max}]} \sum_{1\leq i,j \leq n}  \P_{f}(\Phi_{I_l,J_l}=1\,|\, I_l = i, J_l= j)\cdot \P_{f}(I_l=i,J_l=j) \nonumber \\ 
    \leq&  N_{\max} \underset{1\leq i, j \leq n}{\max}\P_{f}(\Phi_{i,j}=1) \nonumber \leq \alpha.
\end{align*}
\end{proof}

\begin{lemma}
\label{l:lambda A}
Suppose that $I\sim \mathrm{Unif}(\{i,(i+1),\dots,j\})/n$ for some $1\leq i \leq j\leq n$ and let $A$ be a non-empty measurable subset of $[0,1]$. Then 
\begin{equation*}
    \mathbb{P}\left(I\in A^{1/n}\right) \geq \frac{n}{j-i+1}\lambda\left(A\cap \left[\frac{i-1}{n},\frac{j}{n}\right]\right),
\end{equation*}
where $A^{\delta}$ denotes the $\delta$-expansion of $A$ in Definition \ref{e: delta expansion}.
\end{lemma}

\begin{proof}
Let $J$ be a uniformly distributed random variable on $((i-1)/n,j/n]$ and set $I'= \ceil{nJ}/n$. Then $I$ and $I'$ share the same distribution. Furthermore, it always holds that $\abs{J-I'}\leq 1/n$, which implies that $\{J \in A\} \subseteq \{I'\in A^{1/n}\}$. Thus,
\begin{equation*}
    \mathbb{P}\left(I\in A^{1/n}\right) = \mathbb{P}\left(I'\in A^{1/n}\right) \geq  \mathbb{P}(J\in A) = \frac{n}{j-i+1}\lambda\left(A\cap \left[\frac{i-1}{n},\frac{j}{n}\right]\right).
\end{equation*}
\end{proof}

\begin{lemma}
\label{l: 1/8}
Let $A$ be a measurable subset of $[a,b]\subseteq [0,1]$. Assume that $x_i = i/n$ for some $i\in[n]$ satisfying $\abs{x_i-a}\leq 1/n$. Suppose that $J_k$'s are independent and each $J_k\sim \mathrm{Unif}(\{1,2,3,\dots, 2^k\wedge(n-i)\})$ for $k = 0,\dots, \ceil{\log_2(n-i)}$. Then, for sufficiently large $n$,
\begin{equation*}
    \mathbb{P}\left(\text{ there exists } k \text{ such that } x_{i+J_k} \in A^{1/n}\right)\geq \frac{\lambda(A)}{4(b-a)}.
\end{equation*}
Symmetrically, if $\abs{x_i-b}\leq 1/n$, then, for sufficiently large $n$,
\begin{equation*}
    \mathbb{P}\left(\text{ there exists } k \text{ such that } x_{i-J_k'} \in A^{1/n}\right)\geq \frac{\lambda(A)}{4(b-a)},
\end{equation*}
where $J_k'\sim \mathrm{Unif}(\{1,2,3,\dots, 2^k\wedge(i-1)\}$, $k = 0,\dots,\ceil{\log_2(i-1)})$, are independent random variables.
\end{lemma}

\begin{proof}
Let $r = \floor{\log_2(\floor{nb}-i)}$. For $J_{r+1}\sim \mathrm{Unif}(\{1,2,3,\dots,2^{r+1} \wedge(n-i)\})$, the random variable $x_{i+J_{r+1}}\sim \mathrm{Unif}(\{i+1,i+2,i+3,\dots, ((i+2^{r+1})\wedge n)\})/n$. Then by Lemma~\ref{l:lambda A} we obtain
\begin{align*}
    \mathbb{P}\left(x_{i+J_{r+1}}\in A^{1/n}\right) &\geq \frac{n}{2^{r+1}\wedge (n-i)}\lambda\left(A\cap\left[\frac{i}{n},\frac{(i+2^{r+1})\wedge n}{n}\right]\right) \nonumber\\
        &\geq \frac{n}{2^{r+1}}\left(\lambda\left(A\right)-\lambda\left(A\cap \left[a,\frac{i}{n}\right)\right)\right)\nonumber \\
        &\geq \frac{1}{2(b-a)+\frac{2}{n}\cdot \mathbbm{1}_{\{n a\geq i\}} }\left(\lambda\left(A\right)-\left(\frac{i}{n}-a\right)_+\right) 
        \geq\frac{\lambda(A)}{4(b-a)},
\end{align*}
for all $n\geq \max\{2/\lambda(A), 1/(b-a)\}$.
\end{proof}

\begin{lemma}\label{l: upper bound of C n,alpha}
Let $\alpha\in(0,1)$, $\beta \in (0,2]$, $L>0$, $\gamma_n = C_{\beta} h_n^{\beta -\ceil{\beta}+1}$ and $C_h' = \sigma^{-2/(2\beta+1)}C_h$ with $C_{\beta}$ in \eqref{e:C0}, $h_n$ and $C_h$ in \eqref{e:bw}. Suppose that $f\in \Sigma(\beta,L)$ with $\varepsilon_{\ceil{\beta}-1,\gamma_n}(f) \gtrsim n^{-1}$. Set 
\begin{equation*}
    C_n(\alpha) = -2\log\left(\frac{\alpha}{2}\right) \varepsilon_{\ceil{\beta}-1,\gamma_n}^{-1}(f).
\end{equation*}
Then, for any $1\leq i <j \leq n$, we have, 
\begin{subequations}
\begin{align}
    \label{ieq: Cij}
    C_{n,\alpha,i,j}&\leq \sigma  \sqrt{\frac{4W}{C_h^{2\beta+1}}} h_n^{\beta} = \sqrt{\frac{4W}{(C_h')^{2\beta+1}}}h_n^{\beta},\\
    \label{ieq: Cii+1}
    C_{n,\alpha,i,i+1}&\leq \sigma  \frac{8L_2}{\lambda_0}\sqrt{\frac{4W}{C_h^{2\beta+1}}} n^{-1} h_n^{\beta-1} = \frac{8L_2}{\lambda_0}\sqrt{\frac{4W}{(C_h')^{2\beta+1}}}n^{-1} h_n^{\beta-1},
\end{align}    
\end{subequations}
for all sufficiently large $n$ with $C_{n,\alpha,i,j}$ in \eqref{C,n,alpha,I,J} and $W$ in \cref{Theorem: Bounds of variances}, respectively.
\end{lemma}

\begin{proof}
It follows from the definitions of $C_{n,\alpha,i,j}$ in \eqref{C,n,alpha,I,J} and $N_{\max}$ in \eqref{inequality: N< N_max} that for sufficiently large $n$,
\begin{align*}
    C_{n,\alpha,i,j} \leq 
    &\sigma \sqrt{2W} \sqrt{\log\left(\frac{80(\log(2)-\log(\alpha))}{\log(2)}  \varepsilon_{\ceil{\beta}-1,\gamma_n}^{-1}(f)\log^2 n \right) -\log(\alpha)} n^{-\frac{1}{2}}h_n^{-\frac{1}{2}} \\
    \leq& \sigma \sqrt{4W} (\log n)^{\frac{1}{2}}  n^{-\frac{1}{2}}h_n^{-\frac{1}{2}}
    = \sigma  \sqrt{\frac{4W}{C_h^{2\beta+1}}} h_n^{\beta}
    = \sqrt{\frac{4W}{(C_h')^{2\beta+1}}}h_n^{\beta},
\end{align*} 
where the second inequality is due to $\varepsilon_{\ceil{\beta}-1,\gamma_n}(f) \gtrsim n^{-1}$. The upper bound of $C_{n,\alpha,i,i+1}$ can be obtained in a similar way.
\end{proof}

Note that the following result for $ \beta \in(0,1]$ is proven in Section~\ref{ss:proof:rate}.

\begin{theorem}[Separation rate for $1< \beta \le 2$]
\label{t:consistency1b} Under the nonparameteric regression model in \eqref{model}, suppose that Assumptions~\hyperref[M1]{(M1)} and \ref{K1}--\ref{K3} hold with $\beta \in ( 1,2]$ and $L > 0$. Let $(f_n)_{n\geq1}$ be a sequence of function in $\Sigma(\beta,L)$ such that \(\varepsilon_{1,\gamma_n}(f_n) \geq \frac{2}{n}\) with 
\(  \gamma_n = C_{\beta} h_n^{\beta -1} \asymp \left(\log (n)/ n\right)^{\frac{\beta-1}{2\beta+1}}\), $C_{\beta} $ in \eqref{e:C0} and $h_n$ in \eqref{e:bw}. Assume additionally that $D_{\gamma_n}(f_n)\cap [h_n +1/n,1-h_n-1/n]\neq \emptyset$ with $D_{\gamma_n}(f_n)$ in Lemma~\ref{l:epsilon_f}~\ref{l:eps_g:d}. Set
\begin{equation*}
    C_{n}(\alpha) = -2\log\left(\frac{\alpha}{2}\right) \varepsilon_{1,\gamma_n}^{-1}(f_n),
\end{equation*}
in FOMT (\cref{alg:MC}). Then
\begin{equation*}
    \underset{n\to \infty}{\liminf}\,\mathbb{P}_{f_n}(\Phi = 1) \geq 1-\alpha.
\end{equation*}
\end{theorem}

\begin{proof}
The proof follows similarly as in Theorem~\ref{t:consistency1}. We consider the same decomposition~\eqref{e:decomp} but with following definitions of $E$ and $\mathcal{I}$:
\begin{align*}
    E &= \underset{1\leq i\leq n-1}{\bigcup} E_{i},\quad \text{ with } \quad E_{i} = \left\{R_{i,i+1}\leq -\frac{8L_2\sigma}{\lambda_0}\sqrt{\frac{4W}{C_h^{2\beta+1}}}n^{-1}h_n^{\beta-1}\right\}
\end{align*}
with $R_{i,i+1}$ in \eqref{R_I,J}, $C_h$ in \eqref{defn Ch} and $W$ in \cref{Theorem: Bounds of variances}, respectively, and 
\begin{equation*}
    \mathcal{I} = \left\{(i,i+1)\,\Big|\, h_n\leq x_i < x_{i+1}\leq 1-h_n,\,D_{i,i+1}\geq \frac{1}{n}\left(\gamma_n -4Lh_n^{\beta-1} \right)\right\},
\end{equation*}
where $D_{i,j}$ is given in \eqref{eq: D n i j} with bandwidth $h_n$.

\medskip

For the first term in \eqref{e:decomp}, by Mill's ratio and \eqref{variance bound} in \cref{Theorem: Bounds of variances}, we have $\mathbb{P}(E_{i}) \leq \exp\left(-2\log n\right) = n^{-2}$.
Thus, for all large $n$,
\begin{equation*}
    \mathbb{P}(E)\leq  \sum_{i = 1}^{n-1}\mathbb{P}\left( E_{i}\right)\leq n \max_{i}\mathbb{P}\left( E_{i}\right) \leq \frac{1}{n}\leq \frac{\alpha}{2}.
\end{equation*}

\medskip

For the second term in \eqref{e:decomp}, consider an arbitrary pair $(i,i+1)\in \mathcal{I}$. On event $E^c$, we have
\begin{equation*}
    R_{i,i+1} > -\frac{8L_2\sigma}{\lambda_0}\sqrt{\frac{4W}{C_h^{2\beta+1}}}n^{-1}h_n^{\beta-1}.
\end{equation*}
Consequently,
\begin{align*}
    T_{i,i+1} = & D_{i,i+1}+R_{i,i+1} \nonumber\\
    >&\left(C_{\beta}-4L\right)\frac{h_n^{\beta-1}}{n} -\frac{8L_2\sigma}{\lambda_0}\sqrt{\frac{4W}{C_h^{2\beta+1}}}\frac{h_n^{\beta-1}}{n} \nonumber\\
    =&  \frac{8L_2\sigma}{\lambda_0}\sqrt{\frac{4W}{C_h^{2\beta+1}}}\frac{h_n^{\beta-1}}{n} \nonumber\\
    \geq& C_{n,\alpha,i,i+1} = q_{n,\beta,i,i+1}(\alpha),
\end{align*}
where the last inequality follows from \eqref{ieq: Cii+1} and the last equality is implied by $D_{n,\beta,i,i+1}=0$ for $x_i,x_{i+1}\in [h_n,1-h_n]$. Thus, $\Phi_{i,i+1}=1$ for all $(i,i+1)\in \mathcal{P}\cap \mathcal{I}\neq \emptyset$ under $E^c$, which shows that the second term is zero.

\medskip

For the last term in \eqref{e:decomp}, we proceed similarly as in the proof of Theorem~\ref{t:consistency1b}. Applying Lemma~\ref{l:epsilon_f} with $\delta_n = h_n+1/n$ and $\gamma_n = C_1 h_n^{\beta-1}$, we have 
\begin{equation*}
    \lambda\left(\Bar{D}_{\theta_n}(f_n)\right) \geq \varepsilon_{1,\gamma_n}(f_n) \wedge (1-2\delta_n),
\end{equation*}
where 
\begin{equation*}
    \Bar{D}_{\theta_n}(f_n) = \left\{a\in [\delta_n,1-\delta_n]\,\bigg|\, \underset{x\in [a-\delta_n,a+\delta_n]}{\max}f'(x)\leq -\theta_n\right\}\quad \text{ and } \quad \theta_n = \gamma_n -3L\delta_n^{\beta-1}.
\end{equation*}
Here, we may assume that $\varepsilon_{1,\gamma_n}(f_n) \leq 1-2\delta_n$. Otherwise, we have $\varepsilon_{1,\gamma_n}(f_n)\geq 1-2\delta_n \to 1$, as $n \to \infty$, i.e., the function $f_n$ decreases sufficiently fast \enquote{almost} everywhere and its violation is much easier to detect. We then use the same technique as in the case of $\beta \in (0,1]$ to the set $\Bar{D}_{\theta_n}(f_n)$. More precisely, for each $l\in [C_{n}(\alpha)]$ by decomposing $x_i$ into $\Bar{D}^{1/n}_{\theta_n}(f_n)$ and its complement, we obtain
\begin{align*}
    \mathbb{P}_{f_n}(\mathcal{P}_l\cap \mathcal{I}= \emptyset)  \leq \frac{1}{n}\sum_{i:\, x_i \in \Bar{D}^{1/n}_{\theta_n}(f_n)}\mathbb{P}_{f_n}(\mathcal{P}_l\cap \mathcal{I}= \emptyset\,\big| \, I_l = i) + 1-\varepsilon_{1,\gamma_n}(f_n).
\end{align*}
We focus on $I_l = i $ with $x_i \in \Bar{D}_{\theta_n}^{1/n}(f_n)$. Since $x_i \in \bar{D}^{1/n}_{\theta_n}(f_n)$, there exists $a\in \bar{D}_{\theta_n}(f_n)$ such that $\abs{x_i-a}\leq 1/n$, which implies that at least one of $x_{i-1}$ and $x_{i+1}$ is contained in the $1/n$-neighborhood of $a$. With loss of generality, we suppose that $\abs{x_{i+1}-a}\leq 1/n$. By definition $a\in [\delta_n,1-\delta_n] = [h_n+1/n,1-h_n-1/n]$, so both $x_i$ and $x_{i+1}$ are contained in $[h_n,1-h_n]$. Applying Lagrange's mean value theorem and the same proof technique as in Proposition~\ref{p:upper bound Dab}, we have, for large enough $n$,
\begin{align*}
    D_{i,i+1} 
    =& \sum_{k \in I_i} W_{nk}(x_i)\left(f_n(x_k)-f_n(x_{k+1})\right)\nonumber= -\frac{1}{n}\sum_{k \in I_i} W_{nk}(x_i)f_n'(\eta_k)\nonumber\\
    \geq& -\frac{1}{n}\max_{x\in [a-\delta_n,a+\delta_n]}f_n'(x)\nonumber
    \geq \frac{1}{n}(C_{\beta}-4L)h_n^{\beta-1},\nonumber
\end{align*}
for some $\eta_k \in (x_k,x_{k+1})$ for all $k\in I_i$. This implies that $(i,i+1)\in \mathcal{I}$. As $(i,i+1)\in \mathcal{P}$ when $I_l = i\leq n-1$, we have  
\begin{equation*}
    \mathbb{P}(\mathcal{P}_l\cap \mathcal{I}= \emptyset\,\big| \, I_l = i) = 0, \quad \text{ for all } i \text{ with } x_i \in \Bar{D}_{\theta_n}^{1/n}(f_n).
\end{equation*}
It follows, for sufficiently large $n$,
\begin{equation*}
    \mathbb{P}(\Phi_{i,j} = 0,\text{ for all }(i,j)\in\mathcal{P};\, \mathcal{P}\cap \mathcal{I} = \emptyset;\,E^c)
    \leq \Big(1-\varepsilon_{1,\gamma_n}(f_n)\Big)^{C_n(\alpha)}\leq \frac{\alpha}{2}.
\end{equation*}

\medskip
Combining the above calculations, we have proven the assertion for $\beta \in (1,2]$.
\end{proof}

\begin{proof}[Proof of \cref{Theorem: Minimax rate optimality}]
\label{proof: minimax rate optimality} We set $C =(2C_{\beta}+2L)\cdot C_h^{\beta-\ceil{\beta}+1}$ and $\gamma_n = C_{\beta}h_n^{\beta-\ceil{\beta}+1}$, where $C_{\beta}$ and $h_n$ are given in \eqref{e:C0} and \eqref{e:bw}, respectively. We prove the statements under $\beta \in(0,1]$ and $\beta \in (1,2]$ separately.
\medskip

{\sc Part~I.} $\beta \in (0,1]$.
\begin{enumerate}[label=(\text{\roman*}), wide]
    \item[\ref{i:opt:a}] Upper bound. Suppose that $(f_n)_{n\in \mathbb{N}}$ is a sequence of functions such that for each $n\in \mathbb{N}$, $f_n \in \mathcal{F}(C\Delta_n)$ and 
    \begin{equation*}
        \lim_{n\to \infty}\mathbb{P}_{f_n}(\Phi=1) = \liminf_{n\to \infty}\underset{f\in {\mathcal{F}}_{\beta}(C\Delta_n)}{\inf}\P (\Phi =1).
    \end{equation*}
    Then by Corollary \ref{p:>=hn} we have $\varepsilon_{0,\gamma_n}(f_n)\geq h_n\gg 1/n$ for all large $n$ , which implies 
    \begin{equation*}
        C_n(\alpha) = -2\log\left(\frac{\alpha}{2}\right)h_n^{-1}\geq -2\log\left(\frac{\alpha}{2}\right)\varepsilon^{-1}_{0,\gamma_n}(f_n).
    \end{equation*}
    Applying \cref{t:consistency1}, we have proven the consistency of FOMT under $\beta\in (0,1]$.
    \item[\ref{i:opt:b}] Lower bound. We prove the lower bound via constructing a series of special functions. Define
\begin{equation*}
    \psi(x) =
    \begin{cases}
    (x+1)^{\beta} & x\in [-1,0],\\
    (1-x)^{\beta} & x\in (0,1],\\
    0   & \text{otherwise}.
    \end{cases}
\end{equation*}
Let $h = h_n \asymp \bigl(\log (n)/n\bigr)^{1/(2\beta+1)}$ and its precise choice will be given later. Set $$g_{n,j}(x) =  Lh_n^{\beta} \psi_{n,j}(x),\quad \text{ with } \quad \psi_{n,j}(x) = \psi\left( \frac{x-2jh_n}{h_n}\right),$$ 
for $j\in I_n \coloneqq \{j\,|\, ((2j-1)h_n, (2j+1)h_n)\subseteq [0,1]\}$.
It is easy to see the following facts:
\begin{itemize}
    \item The support of $g_{n,j}$ is $[(2j-1)h_n, (2j+1)h_n]$.
    \item The cardinality of $I_n$ is $|I_n| \approx 1/(2h_n) = O(h_n^{-1})$.
\end{itemize}
Furthermore, we have $g_{n,j}\in \Sigma(\beta,L)$ for all $j\in I_n$, as
\begin{equation*}
    \abs{g_{n,j}(x)-g_{n,j}(y)} = Lh_n^{\beta}\abs{\psi_{n,j}(x)-\psi_{n,j}(y)}=Lh_n^{\beta} h_n^{-\beta}\abs{\psi(x)-\psi(y)} \leq L \abs{x-y}^{\beta}.
\end{equation*}
Moreover, $g_{n,j}(2jh_n)-g_{n,j}((2j+1)h_n)= Lh_n^{\beta}$, which implies  that $\{g_{n,j} \,|\, j\in I_n \}\subseteq \mathcal{F}_{\beta}(Lh_n^{\beta})$. Let $\Psi$ be any arbitrary $\alpha$-level test, i.e., $\Psi: Y=(Y_1,\dots,Y_n)\mapsto \{0,1\}$ with $\P_{0}(\Psi(Y)=1) = \E_0 (\Psi(Y))\leq \alpha$.
Then
\begin{align*}
    \underset{g\in \mathcal{F}_{\beta}(Lh_n^{\beta})}{\inf}\P_g(\Psi(Y)=1) -\alpha &= \underset{g\in \mathcal{F}_{\beta}(Lh_n^{\beta})}{\inf}\E_g(\Psi(Y))-\alpha \\
    &\leq \underset{j\in I_n}{\min}\, \E_{g_{n,j}}(\Psi(Y)) -\E_0 (\Psi(Y)) \\
    &\leq \frac{1}{|I_n|}\sum_{j\in I_n} \Big( \E_{g_{n,j}}\bigl(\Psi(Y)\bigr)-\E_0 \bigl(\Psi(Y)\bigr)\Big). 
\end{align*}
Note that 
$$\E_{g_{n,j}}(\Psi(Y)) = \int \Psi(Y) d\P_{g_{n,j}} = \int \Psi(Y) \frac{d\P_{g_{n,j}}}{d\P_0}(Y) d\P_0 = \E_0\left( \Psi(Y) \frac{d\P_{g_{n,j}}}{d\P_0}(Y) \right).$$
Thus,
\begin{align*}
    \underset{g\in \mathcal{F}_{\beta}(Lh_n^{\beta})}{\inf}\P_g(\Psi(Y)=1) -\alpha 
        &\leq \left(\frac{1}{|I_n|}\sum_{j\in I_n}\E_0\left( |\Psi(Y)| \left|\frac{d\P_{g_{n,j}}}{d\P_0}(Y)-1\right| \right)\right) \wedge 1\\
    &\leq 
    \E_0\left( \left(\frac{1}{|I_n|}\sum_{j\in I_n} \left|\frac{d\P_{g_{n,j}}}{d\P_0}(Y)-1\right|\right)\wedge 1 \right).
\end{align*}
For simplicity of notation, we assume that $\sigma= 1$. 
Then, under $H_0$, the Radon--Nikodym derivative $d\P_{g_{n,j}}/d\P_0(Y)$ satisfies
\begin{align*}
    \frac{d\P_{g_{n,j}}}{d\P_0}(Y) 
    = \exp \left( \sum_{i=1}^{n}g_{n,j}(x_i)\varepsilon_i -\frac{1}{2} \sum_{i=1}^{n}g_{n,j}^2(x_i)\right). 
\end{align*}
The random variable $$ \sum_{i=1}^{n}g_{n,j}(x_i)\varepsilon_i -\frac{1}{2} \sum_{i=1}^{n}g_{n,j}^2(x_i)$$
has the same distribution as $\omega_n Z-\omega_n^2/2$ with $Z\sim \mathcal{N}(0,1)$. Then, we have
$$\sum_{i=1}^{n}g_{n,j}^2(x_i) = n\normm{g_{n,j}}^2(1+o(1))= nL^2h_n^{2\beta} \normm{\psi_{n,j}}^2\bigl(1+o(1)\bigr)$$
and $$\normm{\psi_{n,j}}^2 = \int_{-\infty}^{\infty} \psi^2\left(\frac{x-2jh_n}{h_n}\right)dx = h_n \int_{-\infty}^{\infty} \psi^2(x)dx = \frac{2}{2\beta+1}h_n.$$
Thus,
$$w_n \coloneqq \sqrt{\sum_{i=1}^{n}g_{n,j}^2(x_i)} =\sqrt{\frac{2}{2\beta+1}}Ln^{\frac{1}{2}}h_n^{\beta+\frac{1}{2}}(1+o(1)).$$
Suppose that $h = h_n = L^{-2/(2\beta+1)}(1-\eta_n) \bigl(\log (n) / n\bigr)^{1/(2\beta+1)}$, where $\eta_n$ satisfies $\eta_n \to 0$ and $\sqrt{\log n }\eta_n \to \infty$. 
Let
\begin{equation*}
    \rho_n = 1-\sqrt{\frac{w_n^2}{2\log(h_n^{-1})}}.
\end{equation*}
Consequently, $w_n = \sqrt{2\log h_n^{-1}}(1-\rho_n)$. Furthermore, we have 
\begin{align*}
    \lim_{n\to \infty} \rho_n 
    = 1- \lim_{n\to \infty} \left(\frac{\frac{2}{2\beta+1}(1-\eta_n)^{2\beta+1}\log n }{\frac{2}{2\beta+1}\left(\log n -\log\log n - (2\beta+1)\log(1-\eta_n)+2\log(L)\right)}\right)^{{1}/{2}} =0.
\end{align*}
Next, 
\begin{align*}
    \lim_{n\to \infty} \sqrt{\log (h_n^{-1})}\rho_n 
    &= \lim_{n \to \infty} \sqrt{\log(h_n^{-1})} \left(1-\frac{w_n^2}{2\log(h_n^{-1})}\right) \left(1+\sqrt{\frac{w_n^2}{2\log(h_n^{-1})}}\right)^{-1}\\
    &= \lim_{n\to \infty} \sqrt{\log(h_n^{-1})} (1-(1-\eta_n)^{2\beta+1})\cdot \frac{1}{2}=\infty.
\end{align*}
It follows from Lemma~\ref{Lemma: 6.2} that 
\begin{equation*}
    \E_0\left( \left(\frac{1}{|I_n|}\sum_{j\in I_n} \left|\frac{d\P_{g_{n,j}}}{d\P_0}(Y)-1\right|\right)\wedge 1 \right) \to 0, \quad \text{ as } n\to \infty.
\end{equation*}
This finishes the proof of lower bound statement with $\beta \in (0,1]$.
\end{enumerate}
{\sc Part~II.} $\beta \in (1,2]$.
\begin{enumerate}[label=(\text{\roman*}), wide]
    \item[\ref{i:opt:a}] Upper bound. Suppose that $(f_n)_{n\in \mathbb{N}}$ is a sequence of functions such that for each $n\in \mathbb{N}$, $f_n \in \mathcal{F}_{\beta}(C\Delta_{\beta,n})$ and 
    \begin{equation*}
    \lim_{n\to \infty}\mathbb{P}_{f_n}(\Phi=1) = \liminf_{n\to \infty}\underset{f\in \mathcal{F}_{\beta}(C \Delta_{\beta,n})}{\inf}\P (\Phi =1),
    \end{equation*}
    with $C =(2C_{\beta}+2L)\cdot C_h^{\beta-1}$ and $\gamma_n = C_{\beta}h_n^{\beta-1}$. Since $f_n \in \mathcal{F}_{\beta}(C\Delta_{\beta,n})$, then there exists $x_n\in (0,1)$ with $f_n'(x_n)\leq - (2C_{\beta}+2L)h_n^{\beta-1}$. Using Corollary~\ref{p:>=hn}, we see that $\varepsilon_{1,n}(f_n)\geq h_n\gg 1/n$ for all $n$ large.
    
    Let $D_{\gamma_n}(f_n)= \{x\in [0,1]\,|\, f_n'(x)\leq -\gamma_n\}$. Clearly, if $x_n \in [h_n+1/n,1-h_n-1/n]$, then $D_{\gamma_n}(f_n)\cap [h_n+1/n,1-h_n-1/n] \neq \emptyset$. For $x_n \not \in [h_n+1/n,1-h_n-1/n]$, without loss of generality, we assume that $x_n \in [0,h_n+1/n]$. Then 
    \begin{equation*}
        f'\left(h_n+\frac{1}{n}\right)\leq f'(x_n) + L\left(h_n+\frac{1}{n}\right)^{\beta-1}\leq - (2C_{\beta}+2L)h_n^{\beta-1}+ 2Lh_n^{\beta-1} < -\gamma_n.
    \end{equation*}
    This ensures that $D_{\gamma_n}(f_n)\cap [h_n+1/n,1-h_n-1/n] \neq \emptyset$.  The assertion follows from \cref{t:consistency1b}.
    \item[\ref{i:opt:b}] Lower bound. We prove the second inequality via constructing a series of special functions. For $\beta \in (1,2]$, define
\begin{equation*}
    \psi(x) = 
    \begin{cases}
    \frac{1}{\beta}(x+2)^{\beta}, & x\in [-2,-1],\\
    -\frac{1}{\beta}(-x)^{\beta}+\frac{2}{\beta}, & x\in (-1,0],\\
    -\frac{1}{\beta}x^{\beta}+\frac{2}{\beta}, & x\in (0,1],\\
    \frac{1}{\beta}(2-x)^{\beta}, & x\in (1,2],\\
    0,   & \text{otherwise}.
    \end{cases}, \text{ then } 
    \psi'(x) =
    \begin{cases}
    (x+2)^{\beta-1}, & x\in [-2,-1],\\
    (-x)^{\beta-1}, & x\in (-1,0],\\
    -x^{\beta-1}, & x\in (0,1],\\
    -(2-x)^{\beta-1}, & x\in (1,2],\\
    0,   & \text{otherwise}.
    \end{cases}
\end{equation*}
Clearly, $\psi\in \Sigma(\beta,1)$. Let $C_{\psi} =C_{\psi}(\beta) = \int_{-\infty}^{\infty} \psi^2(x)dx$ and 
$$h_n = \sqrt{\frac{2}{(2\beta+1)L^2 C_{\psi}}}(1-\eta_n) \left(\frac{\log n}{n}\right)^{1/(2\beta+1)}\asymp\left(\frac{\log n}{n}\right)^{1/(2\beta+1)}, $$
where $\eta_n$ is a sequence satisfying $\lim_{n\to \infty}\eta_n = 0$ and $\lim_{n\to \infty}\sqrt{\log n}\eta_n = \infty$. Set $$
g_{n,j}(x) =  Lh_n^{\beta} \psi_{n,j}(x)\qquad \text{with}\quad \psi_{n,j}(x) = \psi\left( \frac{x-4jh_n}{h_n}\right),$$
for $j\in I_n = \{j\,| \,((4j-2)h_n, (4j+2)h_n)\subseteq [0,1]\}$.
It is easy to see the following facts:
\begin{enumerate}
    \item[1.] The support of $g_{n,j}$ is $[(4j-2)h_n, (4j+2)h_n]$.
    \item[2.] The cardinality of $I_n$ is $|I_n| \approx 1/(4h_n) = O(h_n^{-1})$.
\end{enumerate}
Furthermore, 
\begin{align*}
    \abs{g'_{n,j}(x)-g'_{n,j}(y)} 
=Lh_n^{\beta} h_n^{-1}\left|\psi'\left(\frac{x-4jh_n}{h_n}\right)-\psi'\left(\frac{y-4jh_n}{h_n}\right)\right| 
\leq Lh_n^{\beta} h_n^{-1}\left|\frac{x-y}{h_n}\right|^{\beta-1} = L\abs{x-y}^{\beta-1}.
\end{align*}

Namely, $g_{n,j}\in \Sigma(\beta,L)$ for all $j\in I_n$. Moreover, 
$$g_{n,j}'((4j+1)h_n)= Lh_n^{\beta}h_n^{-1}\psi'(1)=-Lh_n^{\beta-1},$$ 
which implies that $\{g_{n,j} \,|\, j\in I_n \}\subseteq \mathcal{F}_{\beta}(Lh_n^{\beta-1})$. Let $\Psi$ be any arbitrary $\alpha$-level test, i.e., $\Psi: Y=(Y_1,\dots,Y_n)\mapsto \{0,1\}$ with $\P_{0}(\Psi(Y)=1) = \E_0 (\Psi(Y))\leq \alpha$.
Then applying exactly the same technique as in Part I, we have
\begin{align*}
    \underset{g\in \mathcal{F}_{\beta}(Lh_n^{\beta-1})}{\inf}\P_g(\Psi(Y)=1) -\alpha \to 0,\quad\text{ as } n\to \infty. 
\end{align*}
\end{enumerate}
This finishes the proof.
\end{proof}

\begin{proof}[Proof of Corollary \ref{c:unknown sigma}]
Let $f\in \widetilde{\Sigma}(\beta,L)$. Recall \eqref{e:bw} and define, with $S\coloneqq \hat{\sigma}^2_n/\sigma^2$,  
\begin{align*}
    \Tilde{h}_n(S) &:= h_n S^{\frac{1}{2\beta+1}},\\
    \widetilde{C}_{n}(S) &:= -2\log \left(\frac{\alpha}{2}\right)h_n^{-1}S^{-\frac{1}{2\beta+1}},\\
    \widetilde{N}_{\max}(S) &:= -\frac{80}{\log (2)}\log \left(\frac{\alpha}{2}\right) (\log n)^2 h_n^{-1} S^{-\frac{1}{2\beta+1}} ,\\
    \widetilde{C}_{n,\alpha,i,j}(S) &:= \sigma \sqrt{-2\log\left(\frac{\alpha}{\widetilde{N}_{\max}(S)}\right)W}\cdot n^{-\frac{1}{2}}h_n^{-\frac{3}{2}} \left(\frac{8L_2}{\lambda_0}|x_i-x_j|\wedge h_n S^{\frac{1}{2\beta+1}} \right)S^{\frac{\beta-1}{2\beta+1}},\\
    \text{and }\widetilde{D}_{n,\beta,i,j}(S) &:= \begin{cases}
        0,&\quad \text{ if } x_i,x_j\in \left[h_n S^{\frac{1}{2\beta+1}},1-h_n S^{\frac{1}{2\beta+1}}\right],\\
        \left(\frac{16K_{\max}}{\lambda_0}\vee2\right)L h_n^{\beta} S^{\frac{\beta}{2\beta+1}},&\quad \text{ otherwise},
    \end{cases}
\end{align*}
where $K_{\max}$, $W$, $L_2$ and $\lambda_0$ are constants given in Assumption \ref{K1} and \cref{Theorem: Bounds of variances}.
For any fixed $\delta>0$, let  $F_{\delta} = F_{\delta}(f) = \left\{S(f) \notin \left[1-\delta,1+\delta\right]\right\}$.
Then it follows from \eqref{ieq:ub F0} that
\begin{align}
    \label{ieq:ub F} 
    \underset{f\in \widetilde{\Sigma}(\beta,L)}{\sup} \mathbb{P}_{f}(F_{\delta})\to 0,\quad \text{ as } n \to \infty.
\end{align}
{\sc Part~} \ref{c:unknown sigma:a}. For $f\in \widetilde{H}$, we have 
\begin{equation}
\label{ieq:decom Phi F} 
    \mathbb{P}_f(\Phi = 1) = \mathbb{P}_f(\Phi = 1,F_{\delta})+\mathbb{P}_f(\Phi = 1,\,F_{\delta}^c)\leq \mathbb{P}_f(F_{\delta})+ \mathbb{P}_f(\Phi = 1,\,F_{\delta}^c).
\end{equation}
Under $F_{\delta}^c = \{S\in [1-\delta,1+\delta]\}$, we have $\widetilde{N}_{\max}(S)\leq \widetilde{N}_{\max}(1-\delta)$ and
\begin{align}
    \mathbb{P}_f(\Phi = 1,\,F_{\delta}^c) &\leq \mathbb{P}_f\left(\text{ there exists} \,l\in \left\{1,\dots, \widetilde{N}_{\max}(S)\right\} \text{ such that }\Phi_{I_l,J_l} = 1;\,F_{\delta}^c\right) \nonumber\\
    &\leq \sum_{l = 1}^{\widetilde{N}_{\max}(1-\delta)}\mathbb{P}_f(\Phi_{I_l,J_l} = 1;\,F_{\delta}^c)\nonumber\\
    &= \sum_{l = 1}^{\widetilde{N}_{\max}(1-\delta)}\sum_{1\leq i<j\leq n}\mathbb{P}_f(\Phi_{I_l,J_l} = 1,\,I_l = i,\, J_l = j;\, F_{\delta}^c)\nonumber\\
    &\leq \sum_{l = 1}^{\widetilde{N}_{\max}(1-\delta)}\underset{1\leq i<j\leq n}{\max}\mathbb{P}_f(\Phi_{i,j} = 1;\,F_{\delta}^c).\nonumber
\end{align}
For $1\leq i<j\leq n$, 
it follows from Proposition~\ref{p:upper bound Dab} that $D_{i,j}\left(S\right) \leq  \widetilde{D}_{n,\beta,i,j}\left(S\right)$.
Therefore,
\begin{align}
    \,&\mathbb{P}_f(\Phi_{i,j} = 1,\,F^c) \nonumber\\
    \leq\,& \mathbb{P}_f \left(R_{i,j}(S)\geq \widetilde{C}_{n,\alpha,i,j}(S);\,F^c\right) \nonumber \\
    \leq\,&\underset{s\in [1-\delta,1+\delta]}{\max} \mathbb{P}_f \left(R_{i,j}(s)\geq \widetilde{C}_{n,\alpha,i,j}(s);\,F^c\right) \nonumber \\
    \leq& \underset{s\in [1-\delta,1+\delta]}{\max} \exp \left(-\frac{1}{2}\cdot \frac{\widetilde{C}^2_{n,\alpha,i,j}(s)}{\mathbb{V}(R_{i,j}(s))}\right) \nonumber \\
    \leq& \exp \left(-\frac{1}{2}\cdot \underset{s\in [1-\delta,1+\delta]}{\min}\frac{-2\sigma^2 \log\left(\frac{\alpha}{\widetilde{N}_{\max}(s)}\right)Wn^{-1}h_n^{-3}\left(\frac{8L_2}{\lambda_0}\abs{x_i-x_j}\wedge h_n s^{-\frac{1}{2\beta+1}}\right)^2 s^{\frac{2\beta-2}{2\beta+1}}}{\sigma^2 W n^{-1}h_n^{-3}\left(\frac{8L_2}{\lambda_0}\abs{x_i-x_j}\wedge h_n s^{-\frac{1}{2\beta+1}}\right)^2 s^{\frac{2\beta-2}{2\beta+1}}}\right)\nonumber\\
    =&\underset{s\in [1-\delta,1+\delta]}{\max} \frac{\alpha}{\widetilde{N}_{\max}(s)}= \frac{\alpha}{\widetilde{N}_{\max}(1+\delta)},\nonumber
\end{align}
where the last second inequality and the last inequality follow from Mill's ratio and \eqref{variance bound} in \cref{Theorem: Bounds of variances}, respectively. Therefore, as $n\to \infty$,
\begin{equation}
    \label{ieq:Fc4}
    \mathbb{P}_f(\Phi = 1;\,F_{\delta}^c) \leq \widetilde{N}_{\max}(1-\delta) \frac{\alpha}{\widetilde{N}_{\max}(1+\delta)} = \left(\frac{1+\delta}{1-\delta}\right)^{\frac{1}{2\beta+1}} \cdot \alpha,
\end{equation}
which gives an  upper bound uniformly for $f\in \widetilde{H}$.
Combining \eqref{ieq:ub F}, \eqref{ieq:decom Phi F} and \eqref{ieq:Fc4} we see that 
\begin{equation*}
    \underset{n\to \infty}{\limsup}\underset{f\in \widetilde{H}}{\sup}\mathbb{P}_f(\Phi = 1)\leq \underset{n\to \infty}{\limsup}\underset{f\in \widetilde{H}}{\sup}\mathbb{P}_f(\Phi = 1,F_{\delta}^c)\leq \left(\frac{1+\delta}{1-\delta}\right)^{\frac{1}{2\beta+1}} \cdot \alpha.
\end{equation*}
Since this inequality holds for arbitrary $\delta>0$, we have 
\begin{equation*}
    \underset{n\to \infty}{\limsup}\underset{f\in \widetilde{H}}{\sup}\mathbb{P}_f(\Phi = 1) \leq \alpha.
\end{equation*}

\medskip
{\sc Part~} \ref{c:unknown sigma:b}. Consider $\beta \in (0,1]$. 
Let $1\leq i<j\leq n$ and we define the events $E_{i,j}$ and $E$ as follows:
\begin{align*}
    E = \underset{1\leq i<j\leq n}{\bigcup} E_{i,j}, \quad \text{ with } \quad 
    E_{i,j} = \left\{R_{i,j}(S)\leq -\sigma \sqrt{\frac{6W}{C_h^{2\beta+1}}}h_n^{\beta}\right\}
\end{align*}
with $R_{i,j}$ in \eqref{R_I,J}, $C_h$ in \eqref{defn Ch} and $W$ in \cref{Theorem: Bounds of variances}, respectively. 
Suppose that $(f_n)_{n\in \mathbb{N}}$ is a sequence of functions such that for each $n\in \mathbb{N}$, $f_n \in \mathcal{F}_{\beta}(C\Delta_{\beta,n})\cap \widetilde{\Sigma}(\beta,L)$. We have 
\begin{align}
    \label{e: p decom}
    \mathbb{P}_{f_n}(\Phi = 0) =& \, \mathbb{P}_{f_n}(\Phi = 0,F_{\delta}) + \mathbb{P}_{f_n}(\Phi = 0,F_{\delta}^c) \nonumber\\
    \leq&\, \mathbb{P}_{f_n}(F_{\delta}) + \mathbb{P}_{f_n}(F_{\delta}^c\cap E) + \mathbb{P}_{f_n}(\Phi = 0,E^c \cap F_{\delta}^c).
\end{align}
By \eqref{ieq:ub F},  we have $\mathbb{P}_{f_n}(F_{\delta})\to 0$ as $n\to \infty.$
Furthermore, for any fixed $1\leq i<j\leq n$, it follows from Mill's ratio that 
\begin{align*}
    \mathbb{P}_{f_n}(E_{i,j}\cap F_{\delta}^c) =& \, \mathbb{P}\left(R_{i,j}(S)\leq -\sigma \sqrt{\frac{6W}{C_h^{2\beta+1}}}h_n^{\beta},\, S\in [1-\delta,1+\delta]\right)\\
    \leq& \, \underset{s\in [1-\delta,1+\delta]}{\max}\mathbb{P}\left(R_{i,j}(s)\leq -\sigma \sqrt{\frac{6W}{C_h^{2\beta+1}}}h_n^{\beta}\right)\\
    \leq&\, \exp \left(- \frac{1}{2}\frac{\left(\sigma \sqrt{\frac{6W}{C_h^{2\beta+1}}}h_n^{\beta}\right)^2}{\underset{s\in [1-\delta,1+\delta]}{\max}\V(R_{i,j}(s))}\right)\\
    \leq&\, \exp \left(- \frac{1}{2}\frac{\left(\sigma \sqrt{\frac{6W}{C_h^{2\beta+1}}}h_n^{\beta}\right)^2}{\sigma^2W n^{-1}h_n^{-1}(1-\delta)^{-\frac{1}{2\beta+1}}}\right)
    =\, n^{-3(1-\delta)^{\frac{1}{2\beta+1}}}
\end{align*}
where the last inequality follows from \eqref{variance bound}. Therefore, for sufficiently small $\delta>0$,
\begin{equation*}
    \mathbb{P}_{f_n}(E\cap \,F_{\delta}^c) \leq \sum_{1\leq i<j\leq n}\mathbb{P}_{f_n}(E_{i,j} \cap F_{\delta}^c) \leq n^{2-3(1-\delta)^{\frac{1}{2\beta+1}}}\to 0,\quad \text{ as } n \to \infty.
\end{equation*}
Regrading the last term in \eqref{e: p decom}, let 
\begin{equation*}
    \mathcal{I}=\left\{(i,j)\,\Bigg|\, \underset{s\in [1-\delta,1+\delta]}{\min}D_{i,j}(s)-\widetilde{C}_{n,\alpha,i,j}(s)-\widetilde{D}_{n,\beta,i,j}(s) \geq \sigma \sqrt{\frac{6W}{C_n^{2\beta+1}}}h_n^{\beta}\right\},
\end{equation*}
which 
is deterministic and depends only on $f_n$. Furthermore, by \cref{alg:MC}, FOMT accepts $H$ if and only if all conducted local tests return zero. Let $I_1,\dots, I_{C_{n}(S)}$ denote the repeatedly generated indices via uniform sampling. For each $l \in [C_{n}(S)]$, let $\mathcal{J}_l^+$ and $\mathcal{J}_l^-$ denote the sets of generated random indices $J_k$ and $J_k'$ for left and right searches starting from $I_l$. Further, we define
\begin{align*}
    \mathcal{P}(S) = \bigcup_{l \in [C_{n}(S)]}\mathcal{P}_l\qquad \text{with}\quad\mathcal{P}_l &= \left\{(I_l,I_l+J)\,\big|\, J\in \mathcal{J}_l^{+}\right\}\cup \left\{(I_l-J',I_l)\,\big|\, J'\in \mathcal{J}_l^{-}\right\}.
\end{align*}
Although $\mathcal{P}(S)$ depends on $S$, random sets $\mathcal{P}_l$ are independent of $S$ and are i.i.d. distributed. We decompose $\mathbb{P}_{f_n}(\Phi = 0,E^c \cap F_{\delta}^c)$ as follows:
\begin{align}
    \label{eq:p Phi=0 EcFc}
    \mathbb{P}_{f_n}(\Phi = 0,E^c \cap F_{\delta}^c) =& \,\mathbb{P}_{f_n}(\Phi_{i,j} = 0,\,\text{ for all } (i,j)\in \mathcal{P}(S),\,E^c \cap F_{\delta}^c) \nonumber \\
    =&\,\mathbb{P}_{f_n}(\Phi_{i,j} = 0,\,\text{ for all } (i,j)\in \mathcal{P}(S),\, \mathcal{P}(S)\cap \mathcal{I} \neq \emptyset, \,E^c \cap F_{\delta}^c) \nonumber\\
    &+ \,\mathbb{P}_{f_n}(\Phi_{i,j} = 0,\,\text{ for all } (i,j)\in \mathcal{P}(S),\, \mathcal{P}(S)\cap \mathcal{I} = \emptyset,\,E^c \cap F_{\delta}^c)
\end{align}
If $\mathcal{P}(S)\cap \mathcal{I}\neq \emptyset$, then with the event $E^c\cap F_{\delta}^c$ we obtain, for any $(i,j)\in \mathcal{P}(S)\cap \mathcal{I}$, 
\begin{align}
    &\underset{s\in [1-\delta,1+\delta]}{\min} \Big(T_{i,j}(s)-\widetilde{C}_{n,\alpha,i,j}(s)-\widetilde{D}_{n,\beta,i,j}(s) \Big)\nonumber\\
    =& \underset{s\in [1-\delta,1+\delta]}{\min}\Big( D_{i,j}(s)+R_{i,j}(s)-\widetilde{C}_{n,\alpha,i,j}(s)-\widetilde{D}_{n,\beta,i,j}(s) \Big)\nonumber\\
    \geq &\,\sigma \sqrt{\frac{6W}{C_n^{2\beta+1}}}h_n^{\beta} + \underset{s\in [1-\delta,1+\delta]}{\min} R_{i,j}(s)\geq 0, \nonumber
\end{align}
where the second last inequality follows from the definition of $\mathcal{I}$ and the last one is implied by the event $E^c\cap F_{\delta}^c$. This means that $\Phi_{i,j}= 1$ for all $(i,j)\in \mathcal{P}(S)\cap \mathcal{I}$ under the event $E^c\cap F_{\delta}^c$, i.e.,
\begin{equation*}
    \mathbb{P}_{f_n}(\Phi_{i,j} = 0,\,\text{ for all } (i,j)\in \mathcal{P}(S);\, \mathcal{P}(S)\cap \mathcal{I} \neq \emptyset; \,E^c \cap F_{\delta}^c) = 0.
\end{equation*}
For the second probability in \eqref{eq:p Phi=0 EcFc}, we see that
\begin{align*}
    &\,\mathbb{P}_{f_n}(\Phi_{i,j} = 0,\,\text{ for all } (i,j)\in \mathcal{P}(S);\, \mathcal{P}(S)\cap \mathcal{I} = \emptyset;\,E^c \cap F_{\delta}^c) \nonumber \\
    \leq&\,\mathbb{P}_{f_n}(\mathcal{P}(S)\cap \mathcal{I} = \emptyset,\,F_{\delta}^c) \nonumber\\
    =&\,\mathbb{P}_{f_n}\left(\mathcal{P}_l\cap \mathcal{I} = \emptyset,\,\text{ for all } l\in \left[\widetilde{C}_{n}(S)\right];\,F_{\delta}^c\right)\nonumber\\
    \leq&\,\underset{s\in [1-\delta,1+\delta]}{\max}\mathbb{P}_{f_n}\left(\mathcal{P}_l\cap \mathcal{I} = \emptyset,\,\text{ for all } l\in \left[\widetilde{C}_{n}(s)\right]\right)\nonumber
    =\,\mathbb{P}_{f_n}(\mathcal{P}_l\cap \mathcal{I} = \emptyset)^{\widetilde{C}_{n}(1+\delta)}.
\end{align*}
Let $$\gamma_n = \underset{s\in [1-\delta,1+\delta]}{\max}C_{\beta} \tilde{h}_n^{\beta}(s) = C_{\beta} (1+\delta)^{\frac{\beta}{2\beta+1}} h_n^{\beta}$$
with $C_{\beta}$ in \eqref{e:C0}. Let $H_{f_n}(\gamma_n)$, $H_{f_n,R}(\gamma_n)$ and $H_{f_n,L}(\gamma_n)$ denote the sets of $\gamma_n$-heavy points, $\gamma_n$-right-heavy points and $\gamma_n$-left-heavy points of $f_n$, respectively. For any fixed $l$, we have
\begin{align*}
    \,\mathbb{P}_{f_n}(\mathcal{P}_l\cap \mathcal{I}= \emptyset) 
    \leq \, \frac{1}{n}\sum_{i:\, x_i \in H_{f_n}^{1/n}(\gamma_n)}\mathbb{P}_{f_n}(\mathcal{P}_l\cap \mathcal{I}= \emptyset\,\big| \, I_l = i) + 1-\varepsilon_{0,\gamma_n}(f_n).
\end{align*}
Without loss of generality, we assume that $a\in H_{f_n,R} $, i.e., there exists $b\in (a,1]$ such that $\lambda(A)\geq (b-a)/2$, where 
$$A = \left\{x\in [a,b]\,\Big|\, f_n(a)-f_n(x)\geq \gamma_n \equiv C_{\beta}(1+\delta)^{\frac{\beta}{2\beta+1}} h_n^{\beta}\right\}.$$ 
Again as in the proof of \cref{t:consistency1}, we have 
\begin{equation*}
    \mathbb{P}_{f_n}\left(\text{ there exists }j\in \mathcal{J}_l^{+} \text{ satisfying }x_{i+j} \in A^{1/n}\right) \geq 1- \frac{1}{n},
\end{equation*}
and for any $j\in \mathcal{J}_l^{+}$ satisfying $x_{i+j} \in A^{1/n}$,
\begin{align*}
     D_{i,i+j}(s) 
    \geq \, \gamma_n -4L\tilde{h}_n^{\beta}(s)= \gamma_n -4Lh_n^{\beta}s^{\frac{\beta}{2\beta+1}}, \text{ for all } s\in [1-\delta,1+\delta].
\end{align*}
Further, for such a pair $(i,i+j)$, we have
\begin{align*}
    &\underset{s\in [1-\delta,1+\delta]}{\min}D_{i,i+j}(s)-\widetilde{C}_{n,\alpha,i,i+j}(s)-\widetilde{D}_{n,\beta,i,i+j}(s) \\
    \geq& \underset{s\in [1-\delta,1+\delta]}{\min}D_{i,i+j}(s)-\sigma \sqrt{-2\log\left(\frac{\alpha}{\widetilde{N}_{\max}(s)}\right)W}\cdot n^{-\frac{1}{2}}h_n^{-\frac{1}{2}}  s^{\frac{\beta}{2\beta+1}}-\left(\frac{16K_{\max}}{\lambda_0}\vee2\right)L h_n^{\beta} s^{\frac{\beta}{2\beta+1}}\\
    \geq& \underset{s\in [1-\delta,1+\delta]}{\min}\gamma_n -4Lh_n^{\beta}s^{\frac{\beta}{2\beta+1}} - \sigma  \sqrt{\frac{4W}{C_h^{2\beta+1}}} h_n^{\beta} s^{\frac{\beta}{2\beta+1}} -\left(\frac{16K_{\max}}{\lambda_0}\vee2\right)L h_n^{\beta} s^{\frac{\beta}{2\beta+1}}\\ 
    =& \left(C_{\beta} -4L-\sigma  \sqrt{\frac{4W}{C_h^{2\beta+1}}} - \left(\frac{16K_{\max}}{\lambda_0}\vee2\right)L\right)h_n^{\beta} (1+\delta)^{\frac{\beta}{2\beta+1}}\\
    >& \left(C_{\beta} -4L-\sigma  \sqrt{\frac{4W}{C_h^{2\beta+1}}} - \left(\frac{16K_{\max}}{\lambda_0}\vee2\right)L\right)h_n^{\beta} = \sigma \sqrt{\frac{6W}{C_n^{2\beta+1}}}h_n^{\beta},
\end{align*}
which implies $(i,i+j) \in \mathcal{I}$. Thus, given $I_l = i$ and the $\gamma_n$-right-heaviness of $x_i$, we have
\begin{align*}
    &\,\mathbb{P}_{f_n}(\mathcal{P}_l\cap \mathcal{I}\neq \emptyset\,\big| \, I_l = i,\, x_i\in H_{f_n,R}(\gamma_n)) \\
    \geq&\,\mathbb{P}_{f_n}\left(\text{ there exists }j\in \mathcal{J}_l^{+} \text{ satisfying }x_{i+j} \in A^{1/n}\,\big|\,I_l = i,\, x_i\in H_{f_n,R}(\gamma_n)\right) \geq 1-\frac{1}{n}. 
\end{align*}
Analogously, we can achieve the same upper bound with any $\gamma_n$-left-heavy point $x_i$. Thus, 
\begin{equation}
    \label{ieq:ub PlIeE Il2}
    \mathbb{P}_{f_n}\left(\mathcal{P}_l\cap \mathcal{I}=  \emptyset\,\big| \, I_l = i,\, x_i\in H_{f_n,R}(\gamma_n)\right)\leq \frac{1}{n}.
\end{equation}
Finally, note that $f_n\in \mathcal{F}_{\beta}(C\Delta_{\beta,n})$ with $C = (3C_{\beta}+2L)\cdot C_h^{\beta -\ceil{\beta}+1}$, then by Corollary \ref{p:>=hn} we see that all large $n$ and under $F_{\delta}^c$ with $\delta$ small, 
$$\varepsilon_{0,\gamma_n}(f_n) \geq \max_{s\in [1-\delta,1+\delta]}\tilde{h}_n(s)= (1+\delta)^{\frac{1}{2\beta+1}}h_n.$$ Consequently, 
$$\widetilde{C}_n(1+\delta)=\underset{s\in [1-\delta,1+\delta]}{\min}\widetilde{C}_{n}(s) = -2\log \left(\frac{\alpha}{2}\right) \underset{s\in [1-\delta,1+\delta]}{\min} \tilde{h}^{-1}_n(s) \geq -2\log \left(\frac{\alpha}{2}\right) \varepsilon^{-1}_{0,\gamma_n}(f_n).$$
Combining \eqref{ieq:ub F} and \eqref{e: p decom}--\eqref{ieq:ub PlIeE Il2}
we obtain 
\begin{align*}
    &\,\lim_{n\to \infty}\mathbb{P}_{f_n}(\Phi = 0) \\ 
    \leq & \, \lim_{n\to \infty}\mathbb{P}_{f_n}(F_{\delta})+\lim_{n\to \infty}\mathbb{P}_{f_n}(E\cap F_{\delta}^c) + \lim_{n\to \infty}\left(\frac{1}{n}\sum_{i: x_i \in H_{f_n}^{1/n}(\gamma_n)} \frac{1}{n} + 1-\varepsilon_{0,\gamma_n}(f_n)\right)^{\widetilde{C}_n(1+\delta)} \\
    \leq & \,0+0 + \lim_{n\to \infty}\left(\frac{1}{n} + 1-\varepsilon_{0,\gamma_n}(f_n)\right)^{ -2\log \left(\frac{\alpha}{2}\right) \varepsilon^{-1}_{0,\gamma_n}(f)} \\
    \leq&\, \lim_{n\to \infty} \left(1-\frac{1}{2}\varepsilon_{0,\gamma_n}(f_n)\right)^{ -2\log \left(\frac{\alpha}{2}\right) \varepsilon^{-1}_{0,\gamma_n}(f)}< \alpha.
\end{align*}
This proves the second assertion for $\beta \in (0,1]$. The statement with $\beta \in (1,2]$ follows in a similar way.
\end{proof}

\subsection{Proofs for \texorpdfstring{\cref{S: Compleixity Analysis}}{S: Complexity Analysis}}
\subsubsection*{The proposed FOMT}
\begin{proof}[Proof of \cref{t:complexity1}]
Recall from \eqref{f_n hat} that 
\begin{equation*}
    \hat{f}_n(x) = \sum_{k= 1}^{n}W_{nk}(x)Y_k = \sum_{k \in I_x}W_{nk}(x)Y_k, \text{ for all } x\in [0,1]
\end{equation*}
with $\abs{I_x} \leq 2n h_n$. All $W_{nk}(x)$ with $k\in I_x$ can obtained in $ O(nh_n)$ time. Thus, $\hat{f}_n(x)$ can be computed in $O(nh_n)$ steps, which ensures that the computational cost of each local test $\Phi_{i,j}$ is $O(nh_n)$. Since there are at most $C_{n}(\alpha) = -2\log(\alpha/2) h_n^{-1}\asymp h_n^{-1}$ indices of $I$, and for each $I$ there are at most $O(\log^2 n)$ indices of $J$, we see that $\Phi$ computes its results in $O(h_n^{-1}\cdot nh_n \cdot \log^2 n) = O(n \log^2 n)$ steps. Furthermore, for detectable alternative $f\in \mathcal{F}(C\Delta_{\beta,n})$, it follows from \cref{Theorem: Minimax rate optimality} that with probability $1-\alpha$, FOMT is able to detect violation by conducting local tests involving $O(\varepsilon_{\ceil{\beta}-1,\gamma_n}^{-1}(f))$ uniformly distributed $I$'s for all large $n$. Therefore, with probability $1-\alpha$, FOMT has computational complexity $O(\varepsilon_{\ceil{\beta}-1,\gamma_n}^{-1}(f)\cdot n h_n \cdot \log^2 n)$. Plugging in $h_n\asymp \bigl(\log (n)/n\bigr)^{1/(2\beta+1)}$ we then prove the second statement.
\end{proof}

\subsubsection*{The Procedures of \citet{dumbgen2001multiscale,ghosal2000testing,chetverikov2019testing}}
\label{DSGSVC}
The test statistics proposed by \citet{dumbgen2001multiscale,ghosal2000testing,chetverikov2019testing} can be expressed as
\begin{equation*} 
T = \underset{h \in H_n}{\max} \, \underset{l \in L_n(h)}{\max} \frac{S(l, h, K, Y)}{\sigma(l, h, K)} - C(h), 
\end{equation*}
where $h$ represents the \emph{scale} or \emph{bandwidth} used in local estimators, $H_n$  the set of all possible bandwidths candidates, and $L_n(h)$  the set of all possible \emph{locations} associated with a given $h$. The kernel function $K$ is selected based on the context of the test.

Given observations $Y= (Y_1,\dots,Y_n)$, the test statistic is constructed as follows. For a fixed $h \in H_n$ and location $l \in L_n(h)$:
\begin{enumerate}
    \item $S(l, h, K, Y)$ quantifies the discrepancy of the observations from the null hypothesis.
    \item $\sigma(l, h, K)$ serves as a scaling factor to standardize $S(l, h, K, Y)$.
    \item The function $C(h)$ is introduced by \citet{dumbgen2001multiscale} as a \emph{scale calibration} term to ensure that the multiscale test statistic remains finite almost surely. For other test statistics, $C(h)\equiv 0$.
\end{enumerate}
The null hypothesis $f \in H \equiv \Sigma(\beta, L) \cap \mathcal{M}$ will be rejected if $T$ exceeds a critical value. This threshold can be approximated using the limiting distribution of $T$ (see \citealp{ghosal2000testing}) or estimated via Monte--Carlo simulations (see \citealp{dumbgen2001multiscale,chetverikov2019testing}).
Dümbgen and Spokoiny \citet{dumbgen2001multiscale} investigated the Gaussian white noise model and introduce two multiscale test statistics, denoted by $T_{DS,1}$ and $T_{DS,2}$, to examine the monotonicity of regression in $\Sigma(1,L)$ and $\Sigma(2,L)$, respectively. The first test employing $T_{DS,1}$ detects violations of monotonicity by comparing multiscale regression estimates at different locations and scales. A significant deviation between two local estimates suggests a monotonicity violation. More precisely, $T_{DS,1}$ takes the following form
\begin{align*}
    H_n =& \left\{\frac{1}{n},\frac{2}{n},\dots, \frac{\floor{n/2}}{n}\right\},\\
    L_n(h) =& \bigl\{(s,t)\,\big|\, s<t,\, s,t\in \{kh\,|\,k\in \mathbb{N}\}\cap [h,1-h]\bigr\},\\
    S((s,t),h,K,Y) =&  \frac{\sum_{k= 1}^n K\left(\frac{x_k-s}{h}\right)Y_k}{\sigma\sqrt{\sum_{k=1}^n K^2\left(\frac{x_k-s}{h}\right)}}- \frac{\sum_{k= 1}^n K\left(\frac{x_k-t}{h}\right)Y_k}{\sigma\sqrt{\sum_{k=1}^n K^2\left(\frac{x_k-t}{h}\right)}},\\
    \sigma((s,t),h,K) \equiv& 1,\\
    C(h) =& 2\sqrt{\log \left(\frac{1}{2h}\right)},\\
    K(u) =& \1_{[-1,1]}(u)\cdot(1-\abs{u}).
\end{align*}
In contrast, the second testing procedure utilizing $T_{DS,2}$ can be interpreted as estimating the first order derivative and it rejects $H$ if the estimated derivative far below zero. More precisely, $T_{DS,2}$ uses the same $H_n$ and $C(h)$ as in $T_{DS,1}$ but 
\begin{align*}
    L_n(h) =& \{kh\,|\,k\in \mathbb{N}\}\cap [h,1-h],\\
    S(t,h,K,Y) =& \sum_{k= 1}^n K\left(\frac{x_k-t}{h}\right)Y_k,\\
    \sigma(t,h,K) =& \sigma \sqrt{\sum_{k=1}^n K^2\left(\frac{x_k-t}{h}\right)},\\
    K(u) =& \1_{[-1,1]}(u)\cdot u (1-\abs{u}).
\end{align*}

Towards robustness, Ghosal, {\it et al.}~\citet{ghosal2000testing} and Chetverikov~\citet{chetverikov2019testing} proposed a testing procedure based on local versions of Kendall's Tau statistic. The test statistic proposed in \citet{ghosal2000testing}, denoted as $T_{GSV}$, employs a fixed bandwidth $h = h_n$ satisfying $n^{-1/3}\ll h= h_n \ll 1$, corresponding to $H_n = \{h_n\}$ as a singleton. The test statistic $T_{GSV}$ is constructed as follows: 
\begin{align*}
    H_n =& \left\{h_n\right\},\\
    L_n(h) =& \left\{\frac{k}{n}\,\bigg|\,,k\in \mathbb{N}\right\}\cap [0,1],\\
    S(t,h,K,Y) =& -\frac{2}{n(n-1)}\sum_{1\leq i < j\leq n} \sign{(Y_j-Y_i)} K\left(\frac{x_i-t}{h_n}\right)K\left(\frac{x_j-t}{h_n}\right),\\
    \sigma(t,h,K) =& \frac{4}{3n(n-1)(n-2)\sqrt{n}}\sum_{1\leq i,j,k\leq n,\; i\neq j \neq k} \sign{(x_i-x_j)}\sign{(x_i-x_k)} \\
        &\quad\quad\quad \times K\left(\frac{x_j-t}{h_n}\right)K\left(\frac{x_k-t}{h_n}\right)K^2\left(\frac{x_i-t}{h_n}\right),\\
        \approx& \frac{4}{3}\int \int\int \sign{(\omega_1-\omega_2)}\sign{(\omega_1-\omega_3)} 
        K\left(\frac{\omega_2-t}{h_n}\right)K\left(\frac{\omega_3-t}{h_n}\right)\\
        & \quad \quad \quad \times K^2\left(\frac{\omega_1-t}{h_n}\right) \1_{[0,1]}(\omega_1)
        \1_{[0,1]}(\omega_2)\1_{[0,1]}(\omega_3) d\omega_1 d\omega_2 d\omega_3,\\
    \sign{(x)} \coloneqq &
    \begin{cases}
    1,\quad & \quad \text{if }x>0,\\
    0,\quad & \quad \text{if }x=0,\\
    -1, \quad & \quad \text{if }x<0.
    \end{cases}        
\end{align*}
Here, $K$ is assumed to be twice continuously differentiable and to have compact support contained in $[-1,1]$. A key limitation of this method is its lack of adaptivity: it requires prior knowledge of the smoothness parameter $\beta$, which determines the minimax-optimal bandwidth $h_n$. The test statistic $T_{GSV}$ converges in distribution to a known limiting distribution, enabling the direct selection of the corresponding quantiles as critical values. Consequently, the critical value for the procedure proposed by \citet{ghosal2000testing} can be obtained in $O(1)$ time. In contrast, the critical values for other testing procedures are determined through Monte--Carlo simulations with $R$ repetitions.

To overcome the adaptivity issue, an alternative test statistic, $T_C$, proposed in \citet{chetverikov2019testing}, employs a sequence of exponentially decreasing bandwidths and a generalized Kendall’s Tau statistic. The test statistic $T_C$ is defined as:
\begin{align*}
    H_n =& \left\{\frac{1}{2^k}\,\bigg|\,k\in \mathbb{N}\right\}\cap \left[\left(\frac{\log n}{n}\right)^{\frac{1}{3}} \cdot \frac{1}{5},\frac{1}{2}\right],\\
    L_n(h) =& \left\{\frac{k}{n}\,\bigg|\,k\in \mathbb{N}\right\}\cap [0,1],\\
    S(t,h,K,Y) =& \frac{1}{2}\sum_{1\leq i, j\leq n} (Y_j-Y_i)\sign{(x_j-x_i)} \abs{x_i-x_j}^p  K\left(\frac{x_i-t}{h_n}\right)K\left(\frac{x_j-t}{h_n}\right),\\
    \sigma(t,h,K) =& \sigma \sqrt{\sum_{i=1}^{n}\left(\sum_{j=1}^{n}\sign{(x_j-x_i)}\abs{x_i-x_j}^p  K\left(\frac{x_i-t}{h_n}\right)K\left(\frac{x_j-t}{h_n}\right) \right)^2},      
\end{align*}
where $p \in \{0,1\}$ is a tuning parameter.

\begin{theorem}
\label{t: computational complexity kernel}
Suppose that Assumptions \hyperref[A]{(A)} and \hyperref[B]{(B)} hold in the nonparametric regression model~\eqref{model}. Then: 
\begin{enumerate}[label=(\roman*), wide]
    \item The computational complexities of two testing procedures in \citet{dumbgen2001multiscale} are of the same order, more precisely, $O(Rn^2)$.
    \item The test in \citet{ghosal2000testing} has computational complexity $O(n^3h_n^2)$ with bandwidth $h_n$ satisfying $n^{-1/3}\ll h_n \ll 1$.
    \item The computational complexity of test proposed by \citet{chetverikov2019testing} is $O(Rn^3)$. 
\end{enumerate}
\end{theorem}
\begin{proof}
The computational cost of all aforementioned procedures is independent of the choice of $K$, so we use $S(l,h,Y)$ and $\sigma(l,h)$ instead. 

For any fixed $l$ and $h$, the exact forms of all aforementioned $\sigma(l,h)$ in all procedures can be computed explicitly in $O(1)$ time. Thus, the dominant computational cost arises from evaluating $S(l,h,Y)$. In $T_{GSV}$ and $T_{C}$, where local versions of Kendall's Tau statistic are used, computing $S(l,h,Y)$ requires $O(n^2h^2)$ operations for each $l$ and $h$. In $T_{DS,2}$, computing $S(l,h,Y)$ takes $O(nh)$ steps. Summing over all computational costs over all $h\in H_n$ and $l\in L_n(h)$ we obtain computational complexities of each test statistic,
\begin{align*}
    T_{DS,2}: \sum_{h\in H_n} \sum_{l\in L_n(h)}O(nh) =&\sum_{h\in H_n} O(h^{-1})\cdot O(nh)=\sum_{h\in H_n}O(n) = O(n^2),\\
    T_{GSV}: \sum_{h\in H_n} \sum_{l\in L_n(h)}O(nh) =& O(n)\cdot O(n^2h_n^2) = O(n^3h_n^2),\\
    T_{C}:\quad \sum_{h\in H_n} \sum_{l\in L_n(h)}O(nh) =& \sum_{h\in H_n} O(n)\cdot O(n^2h^2)\leq O(n^3)\sum_{k=1}^{\infty} \cdot\frac{1}{2^k}= O(n^3).
\end{align*}
Next we consider $T_{DS,1}$. For each $t \in \left\{k/n\,|\,,k\in \mathbb{N}\right\}\cap [0,1]$, it requires $O(nh)$ steps to computes
\begin{equation*}
    \widehat{\Psi}_{t,h} \coloneqq\frac{\sum_{k= 1}^n K\left(\frac{x_k-t}{h}\right)Y_k}{\sigma\sqrt{\sum_{k=1}^n K^2\left(\frac{x_k-t}{h}\right)}},
\end{equation*}
then the computational cost of all $\widehat{\Psi}_{t,h}$ is $O(nh\cdot h)$. Then we compare $\widehat{\Psi}_{s,h}$ and $\widehat{\Psi}_{t,h}$ over all pairs $(s,t)\in L_{n}(h)$, which requires $O(h^{-2})$ steps. Sum up these computational cost over all $H_n$ we then obtain the computational complexity 
\begin{equation*}
    T_{DS,1}: \sum_{h\in H_n} O(n\vee h^{-2})= O(n^2) \vee O\left(n^2 \cdot\sum_{k=1}^{n/2}\frac{1}{k^2}\right)=O(n^2).
\end{equation*}
Note that the critical values of $T_{DS,1}$, $T_{DS,2}$, $T_{GSV}$ and $T_C$ are simulated via Monte--Carlo method by replacing observations by i.i.d. Gaussian noises. We then obtain their computational complexity in the first and third statements. For the testing procedure in \citet{ghosal2000testing}, its critical is directly accessible, then it has computational complexity $O(n^3h_n^2)$.
\end{proof}

\begin{table}[H]
\centering{\scriptsize
\begin{tabular}{| l| l| l| l|}
\toprule
\hline \xrowht[()]{10pt}
\textbf{Methods} & \textbf{$H_n$}  &\textbf{$L_{n}(h)$} & \textbf{$S(l,h,K,Y)$} \\[5pt]
    \hline \xrowht[()]{10pt}
    DS$_1$ & $\{k/n\,|\,k\in \mathbb{N}\}\cap[0,1/2]$ & $\{(s,t)\,|\, s<t,\, s,t\in \{kh\,|\,k\in \mathbb{N} \}\cap[0,1]\}$ & \makecell{Difference between \\multiscale estimates \\at $s$ and $t$} \\[5pt]
    \hline \xrowht[()]{10pt}
    DS$_2$ & $\{k/n\,|\,k\in \mathbb{N}\}\cap[0,1/2]$ & $\{kh\,|\,k\in \mathbb{N} \}\cap[0,1]$ & Multiscale estimators \\[5pt]
    \hline \xrowht[()]{10pt}
    GSV & $\{h_n\}$, with $n^{-1/3}\ll h_n \ll 1$ & $\{k/n\,|\,k\in\mathbb{N}\}\cap[0,1]$ & Local Kendall's Tau  \\[5pt]
    \hline \xrowht[()]{10pt}
    C & $\left\{h = (1/2)^k\,|\, k\in \mathbb{N} \text{ with } h\geq 0.2\cdot \left(\frac{\log(n)}{n}\right)^{\frac{1}{3}} \right\}$ & $\{k/n\,|\,k\in\mathbb{N}\}\cap[0,1]$ & \makecell{Modified local \\Kendall's Tau} \\[5pt]
    \hline
\bottomrule
\end{tabular}
\caption{Comparison of current existing procedures. We employ the abbreviations DS for \citet{dumbgen2001multiscale}, GSV for \citet{ghosal2000testing}, C for \citet{chetverikov2019testing} and BHL for \citet{baraud2005testing}. }
\label{table: details of kernel type methods}}
\end{table}

\subsubsection*{Baraud, {\it et al.}~\citet{baraud2005testing}'s Procedure}
\label{Details: Baraud Huet and Laurent's Procedure}
\citet{baraud2005testing} proposed a testing procedure based on partitioning the interval $[0,1]$
uniformly into $l_n$ subintervals. Their test statistic utilizes regression line slopes computed over various scales and locations. The critical values are determined via Monte--Carlo simulation. To ensure statistical consistency, they recommend choosing $l_n = n/2$.

Before formally defining the test statistic, we introduce the following notation:
\begin{enumerate}
    \item We define an almost regular partitions of the set of indices $[n] \equiv \{1,\dots,n\}$ into $l_n$ sets as follows: for each $k$ in $[n]$, let
    \begin{equation*}
        J_k\coloneqq \left\{i\in [n] \,\bigg|\, \frac{k-1}{l_n}<\frac{i}{n}\leq \frac{k}{l_n}\right\}
    \end{equation*}
    and define the partition as
    \begin{equation*}
        \mathcal{J}^{l_n} \coloneqq \{J_k\;|\; k \in [n]\}.
    \end{equation*}
    Further, for each $l \in [n]$, we gather consecutive sets $J_k$ to obtain a partition of $[n]$ with $l$ sets. This new partition is given by
    \begin{equation*}
        \mathcal{J}^{l} \coloneqq \left\{ J_j^l= \underset{(j-1)/l<k/l_n\leq j/l}{\bigcup}J_k \;\bigg|\; j\in[l]\right\}.
    \end{equation*}
    \item For any vector $v \in \R^n $ and any subset $J$ of $[n]$, let $v_J$ denote the vector in $\R^n$ whose coordinates coincides with those of $v$ on $J$ and vanish elsewhere. That is,
    \begin{align*}
        v_{J} =& (u_j)_{j=1,\dots,n},\\
        u_j =&
        \begin{cases}
        v_j & \text{ if } j \in J,\\
        0   & \; \text{else}.
        \end{cases}
    \end{align*}
    Let $\Bar{v}_J$ denote the quantity $\sum_{j\in J}v_i/|J|$, where $|J|$ denote the cardinality of $J$.
    \item Let $\mathbf{1}$ denote the $\R^n$-vector $(1,\dots,n)^{\top}$. Moreover, we define $V_{n}$ as the linear span of $\{\mathbf{1}_J\,|\; J\in \mathcal{J}^{l_n}\}$. Note that the dimension of $V_n$ is $l_n$.
\end{enumerate}
Let $I$ and $J$ be two disjoint subsets of $[n]$ such that $I$ is on the left of $J$ in the sense that every element in $I$ is smaller than every element in $J$. For any $l\in \{2,\dots,l_n\}$, let
\begin{equation}
    \label{defn: test statistic T^l}
    T^l(Y) \coloneqq \sqrt{n-l_n}\underset{1\leq i<j \leq l}{\max}N_{i,j}^l \frac{\Bar{Y}_{J_i^l}-\Bar{Y}_{J_j^l}}{\norm{Y-\Pi_{V_n}Y}},
\end{equation}
where $\Pi_{V_n}$ denotes the orthogonal projection from $\R^n$ onto $V_n$ and
\begin{equation*}
    N_{i,j}^l\coloneqq \left(\frac{1}{|J_i^l|}+\frac{1}{|J_j^l|}\right)^{-\frac{1}{2}}.
\end{equation*}
Further, let $q(l,u)$ denote the $(1-u)$-quantile of the random variable $T^l(\varepsilon)$ with $\varepsilon = (\varepsilon_1,\dots,\varepsilon_n)$.
Finally, the test statistic $T_{BHL}$ of \citet{baraud2005testing} is given by
\begin{equation*}
    T_{BHL}(Y,u_{\alpha})= \underset{l=2,\dots,l_n}{\max}(T^l(Y)-q(l,u_\alpha))
\end{equation*}
where $u_\alpha:=\sup\left\{u\in (0,1)\bigg| \; \P\left(T_{BHL}(\varepsilon,u)>0\right)\leq \alpha\right\}$ is calculated by Monte--Carlo simulation as follows:
\begin{enumerate}
    \item Find a suitable set of grid values $U = \{u_1,\dots,u_m\}$.
    \item Generate i.i.d samples $\varepsilon_i\sim \N(0,\sigma^2)$ for $i \in [n]$.
    \item Compute $T^l({\boldsymbol\varepsilon})$ in \eqref{defn: test statistic T^l} for all $l = 2,\dots,l_n$ with ${\boldsymbol\varepsilon} = (\varepsilon_i)_{i = 1}^n$ generated in the second step.
    \item Repeat the second and third steps $R$ times to compute $q(l,u_j)$ for all $l = 2,\dots,l_n$ and $j= 1,\dots,m$.
    \item Compute $p(u_j) \coloneqq \P(T_{BHL}(\varepsilon,u_j)>0)$ for all $j = 1,\dots,m$.
    \item Take $u_{\alpha} \coloneqq \max\{u_j|\; p(u_j)\leq \alpha\}$.
\end{enumerate}
The null hypothesis is rejected if $T(Y,u_{\alpha})>0$. 

\begin{theorem}
\label{Theorem: complexity of BHL}
Suppose that Assumptions \hyperref[A]{(A)} and \hyperref[B]{(B)} hold in the nonparametric regression model~\eqref{model}. Then, the testing procedure proposed by \citet{baraud2005testing} has a computational complexity of \( O(R(l_n^3 \vee n l_n)) \).
\end{theorem}

\begin{proof}
Let ${\boldsymbol\varepsilon} := (\varepsilon_1,\dots,\varepsilon_n)$ with i.i.d. $\varepsilon_i\sim \N(0,\sigma^2)$. Then,
\begin{equation*}
    T^l({\boldsymbol\varepsilon}) \coloneqq \sqrt{n-l_n} \underset{1\leq i<j \leq l}{\max} N_{i,j}^l \frac{\Bar{\varepsilon}_{J_i^l}-\Bar{\varepsilon}_{J_j^l}}{\norm{\varepsilon-\Pi_{V_n}\varepsilon}}.
\end{equation*}
Note that \( \{\mathbf{1}_{J_i} \mid i=1,\dots, l_n\} \) is an orthogonal basis of \( V_n \). The projection of \( {\boldsymbol\varepsilon} \) onto \( V_n \) is
\[
\Pi_{V_n}{\boldsymbol\varepsilon} = \sum_{i=1}^{l_n} \langle \varepsilon, {\boldsymbol u}_i \rangle {\boldsymbol u}_i,
\]
where \(\langle {\boldsymbol\varepsilon}, {\boldsymbol u}_i \rangle = \left \langle {\boldsymbol\varepsilon}, {\mathbf{1}_{J_i}}/{\sqrt{|J_i|}} \right\rangle =  \sum_{j \in J_i} \varepsilon_j/{\sqrt{|J_i|}}\) and \({\boldsymbol u}_i \coloneqq {\mathbf{1}_{J_i}}/{\sqrt{|J_i|}}\),  \(i=1,\dots,l_n\). Computing \( \langle {\boldsymbol\varepsilon}, {\boldsymbol u}_i \rangle \) requires \( O(|J_i|) \) operations, and summing over all \( i = 1, \dots, l_n \) gives a cost of  
\(
O\bigl(\sum_{i=1}^{l_n} |J_i|\bigr) = O(n).
\)
Thus, computing \( \Pi_{V_n} {\boldsymbol\varepsilon} \) and \( \norm{{\boldsymbol\varepsilon} - \Pi_{V_n} {\boldsymbol\varepsilon}} \) takes \( O(n) \) steps.

For any fixed \( l \in \{2,\dots,l_n\} \), computing \( \Bar{\varepsilon}_{J_i^l} \) for all \( i=1,\dots,l \) requires  
\(
O\bigl(\sum_{i=1}^{l} |J_i^l|\bigr) = O(n)
\)  
steps. Since the maximum in \( T^l({\boldsymbol\varepsilon}) \) is taken over \( O(l^2) \) pairs \( (i,j) \), this adds an additional \( O(l^2) \) computational cost. Hence, the total complexity of computing \( T^l(\varepsilon) \) is
\(O(l^2 \vee n).\)
Summing over all \( l \in \{2,\dots,l_n\} \), the total complexity for computing \( T^l(\varepsilon) \) for all \( l \) is
\(
\sum_{l=2}^{l_n} O(l^2 \vee n) = O(l_n^3 \vee n l_n).
\)
Since the test statistic relies on Monte--Carlo estimation of the critical values, we repeat the computation of \( T^l(\varepsilon) \) for \( R \) times. This results in a total complexity of
\(
O(R(l_n^3 \vee n l_n)).
\)
To estimate the quantiles \( q(l,u_j) \) for all \( u_j \in U \), we sort the \( R \) simulated values of \( T^l(\varepsilon) \), which requires \( O(mR\log R) \) operations, where \( m \) is the number of grid points in \( U \). Based on these quantiles, we compute the values \( T_{BHL}(\varepsilon, u_j) \) and the probabilities \( p(u_j) = \mathbb{P}(T_{BHL}(\varepsilon, u_j) > 0) \) for all \( u_j \in U \). This requires \( O(m R l_n) \) additional steps. Thus, computing \( T_{BHL}(Y, u_{\alpha}) \) takes \( O(l_n^3 \vee n l_n) \) steps.

Combining all the computational costs, the total complexity is
\begin{equation*}
    O(mR(l_n+\log(R)))+O(R(l_n^3 \vee n l_n))=O(R(l_n^3 \vee n l_n))
\end{equation*}
with a user-specified constant $m\in \N$.
\end{proof}

\begin{remark}
Note that $l_n \in[n]$, the computational complexity of $\Psi_{BHL}$ is lower bounded by $O(Rn)$. However, the statistical guarantee in \citet{baraud2005testing} becomes invalid for $l_n=O(1)$.   
\end{remark}  

\subsection{Proofs for \texorpdfstring{\cref{S: Adaptivity}}{S: Adaptivity}}
\subsubsection*{Proofs for \texorpdfstring{\cref{SS: A new Lepskii principle}}{SS: A new Lepskii principle}}
Let
\begin{equation}
    \label{e:C rho}
    C_{\rho} = \sqrt{\frac{32\sigma^2 \mu_2}{\lambda_0^2}},
\end{equation}
where $\mu_2 \coloneqq \int_{-\infty}^{\infty} K^2(u)du$ and $\lambda_0$ in \eqref{e: lambda_0}. Then for any $\beta \in (0,2]$, applying the same technique as in the proof of \citet[Theorem~1.8]{tsybakov2008introduction} we have 
\begin{equation}
    \label{ieq: rho_n(h)2}
    \mathbb{E}\left(\rho^2_{n}(h)\right)\leq q_2 \frac{\log n}{nh} \leq C_{\rho}^2 \frac{\log n}{nh}=: G_2^2(h), \quad \text{ for all } \, h\in (0,0.5],  
\end{equation}
where $q_2$ is given in \eqref{e:bw}.

Let $\{x_1,\dots,x_M\}$ be a finite set of random elements in a metric space $(\mathcal{X}, d)$, given on a probability space $(\Omega, \mathcal{F}, \mathbb{P})$. Let $G_1:[M] \to \mathbb{R}_+$ be a decreasing function. Suppose that $x \in \mathcal{M}$ is any fixed element. 
\begin{definition}[Admissibility, \citealp{mathe2006regularization}]
A nondecreasing function $G_1:[M]\to \mathbb{R}_+$ is called \emph{admissible} for $x$ if there exists a family of nonnegative random variables $\rho(m)$, for $m \in [M]$, for which
\begin{align*}
    d(x,x_m) \leq & \,G_1(m) + \rho (m), \quad m \in [M],\\
    \mathbb{E}(\rho^2(m)) \leq& \,G_2(m),\quad \quad \quad \quad \, \,m \in[M],\\
   \text{and} \qquad G_1(1) \leq& \,G_2(1).
\end{align*}
\end{definition}
Towards \cref{theorem: Lepskii}, we start with the deterministic oracle principle. Let $\{x_1,\dots,x_M\}$ be a finite set of \emph{deterministic} elements in a metric space $(\mathcal{X}, d)$ and we set $\rho(m)\equiv G_2(m)$. Then a nondecreasing function $G_1$ is admissible for $x$ if and only if 
\begin{equation*}
    d(x,x_m) \leq G_1(m)+G_2(m),\quad \text{ for all } m \in[M],
\end{equation*}
and $G_1(1)\leq G_2(1)$.
\begin{lemma}
\label{lemma: Mathe}
Let $\{x_1,\dots,x_M\}$ be a finite set of \emph{deterministic} elements in a metric space $(\mathcal{X}, d)$ and let $G_2:[M]\to \mathbb{R}_+$ be a decreasing function and $\kappa>1$. Define
\begin{equation*}
    m^* = \max \{m\,|\, \text{there is an admissible } G_1 \text{ with } G_1(m)\leq G_2(m)\}.
\end{equation*}
For 
\begin{equation*}
    \bar{m} = \min\{m\,|\,\text{ there exists }k \in [m-1] \text{ such that } d(x_k,x_m)> 4 \kappa G_2(k)\}-1,
\end{equation*}
we have
\begin{equation*}
    d(x,x_{\bar{m}})\leq 6G_2(m^*).
\end{equation*}
If additionally there exists $D\in (0,+\infty)$ such that $G_2(m)\leq D \cdot G_2(m+1)$, then 
\begin{equation*}
    d(x,x_{\bar{m}})\leq 6D \min\{G_1(m)+ \kappa G_2(m)\,|\,m\in[M],\, G_1 \text{ is admissible }\}
\end{equation*}
\end{lemma}

\begin{proof}
The proof follows in the same way as \citet[Lemma~4 and Corollary~1]{mathe2006regularization}. 
\end{proof}

\begin{lemma}\label{l:tail of rho}
Under the nonparametric regression model \eqref{model}, suppose that Assumptions \ref{K1}--\ref{K2} hold. Let $\beta \in(0,2] $ and $\rho_n(h)$ be the random term of LPE($\ceil{\beta}-1$) with bandwidth $h>0$, i.e.,
\begin{equation*}
    \rho_n(h) \coloneqq \underset{i= 1,\dots,n}{\max}\left| \mathbb{E}(\hat{f}_{n}(x_i;h))-\hat{f}_{n}(x_i;h) \right|.
\end{equation*} 
Then, for $G_2$ defined in \eqref{ieq: rho_n(h)2}, it holds
\begin{equation*}
    \mathbb{P}\bigl(\rho_n(h)\geq t G_2(h)\bigr) \leq  n^{-(t-1)^2{\mu_2 }/({4K^2_{\max}})},\quad \, \text{ for all } t\geq1.
\end{equation*}
\end{lemma}

\begin{proof}
We have
\begin{equation*}
    \rho_n(h) = \underset{i=1,\dots,n}{\max}\left|\sum_{k=1}^n W_{nk}(x_i)\varepsilon_k\right|=\sigma \underset{i=1,\dots,n}{\max}\left|\sum_{k=1}^n W_{nk}(x_i)\xi_k\right|\eqqcolon g(\xi),
\end{equation*}
where $\xi_i \coloneqq \varepsilon_i/\sigma$ are independent standard normal distributed. Let $\xi,\xi'\in \mathbb{R}^n$ and by $i^*$ denote the maximizer of $g(\xi)$. Then
\begin{align*}
    g(\xi)-g(\xi') =& \sigma \left|\sum_{k=1}^n W_{nk}(x_{i^*})\xi_k\right|-\sigma \underset{i=1,\dots,n}{\max}\left|\sum_{k=1}^n W_{nk}(x_{i})\xi_k'\right|\\
    \leq& \sigma \left|\sum_{k=1}^n W_{nk}(x_{i^*})\xi_k\right|-\sigma \left|\sum_{k=1}^n W_{nk}(x_{i^*})\xi_k'\right|\\
    \leq& \sigma \left|\sum_{k=1}^n W_{nk}(x_{i^*})(\xi_k-\xi_k')\right|\\
    \leq& \sigma \sqrt{\sum_{k=1}^n W_{nk}^2(x_{i^*})}\norm{\xi-\xi'}
    \leq \sigma \frac{C_*}{\sqrt{nh}}\norm{\xi-\xi'},
\end{align*}
where $C_*$ is given in Lemma~\ref{lemma of weights}. Exchanging $\xi$ and $\xi'$, we obtain that $g$ is Lipschitz continuous with $L_g = \sigma C_* /\sqrt{nh}$.
Furthermore, it follows from \eqref{ieq: rho_n(h)2} that
\begin{equation*}
    \mathbb{E}\bigl(g(\xi)\bigr) = \mathbb{E}\bigl(\rho_n(h)\bigr) \leq \sqrt{\mathbb{E}\bigl(\rho_n^2(h)\bigr)} \leq \sqrt{\frac{q_2\log n}{nh}}\leq G_2(h).
\end{equation*}
By \citet[Theorem~2.26]{wainwright2019high} we have
\begin{equation*}
    \label{ieq: tail of A}
    \mathbb{P}\bigl(g(\xi) - \mathbb{E}(g(\xi))\geq t\bigr) \leq \exp\left(-\frac{t^2}{2L_g^2}\right)\quad \text{ for all } t>0.
\end{equation*}
Therefore, for all $t\geq 1$, 
\begin{align*}
    \mathbb{P}\bigl(\rho_n(h)\geq tG_2(h)\bigr) 
    =& \,\mathbb{P}\bigl(g(\xi) - \mathbb{E}(g(\xi)) \geq tG_2(h) - \mathbb{E}(g(\xi))\bigr)\nonumber\\
    \leq& \,\mathbb{P}\bigl(g(\xi) - \mathbb{E}(g(\xi))\geq (t-1)G_2(h)\bigr) \nonumber\\
    \leq& \exp{\left(-\frac{\mu_2 (t-1)^2}{4K^2_{\max}}\log n \right)} = n^{-{(t-1)^2\mu_2 }/({4K^2_{\max}})},
\end{align*}
which concludes the proof.
\end{proof}

\begin{proof}[Proof of \cref{theorem: Lepskii}]
Let $\kappa>1$ and define
\begin{equation}\label{e:defn E kappa}
    \Pi(\omega) \coloneqq \underset{m\in [M]}{\max}\frac{\rho_n(m)(\omega)}{G_2(m)}, \;\text{ and }\;
    E_{\kappa} \coloneqq \{\omega\,|\, \Pi(\omega) \geq \kappa\}.
\end{equation}
Note that for any fixed $\omega \in E_{\kappa}^c$, $\rho_n(m)(\omega)$ is fixed and satisfies $\rho_n(m)(\omega) \leq \kappa G_2(m)$ for all $m \in [M]$. Thus, by Lemma \ref{lemma: Mathe}, we have
\begin{equation*}
    d_{\mathcal{A}}(f,\hat{f}_{n,\Bar{m}})\leq 12 \underset{m \in [M]}{\min}\{G_1(m)+\kappa G_2(m)\}\leq C_{\max} \left(\frac{\log n}{n}\right)^{\frac{\beta}{2\beta+1}},
\end{equation*}
for some constant $C_{\max}$ depends only on $G_1$, $G_2$ and $\kappa$.
Further, by Lemma~\ref{l:tail of rho},
\begin{equation}\label{ieq: E_kappa}
    \mathbb{P}(E_{\kappa})\;\leq\; \sum_{m = 1}^M \mathbb{P}\bigl(\rho_n(m)\geq \kappa G_2(m)\bigr)\;\leq\; M n^{-{(\kappa-1)^2\mu_2 }/({4K^2_{\max}})}.
\end{equation}
Therefore, 
\begin{equation*}
    \mathbb{P}\left(d_{\mathcal{A}}(f,\hat{f}_{n,\Bar{m}})\leq 12 \underset{m \in [M]}{\min}\{G_1(m)+\kappa G_2(m)\}\right)\geq 1- Mn^{-{(\kappa-1)^2\mu_2 }/({4K^2_{\max}})},
\end{equation*}
which implies, for $n \ge 2$, 
\begin{equation*}
    \mathbb{P}\left(d_{\mathcal{A}}(f,\hat{f}_{n,\Bar{m}})\leq C_{\max} \left(\frac{\log n}{n}\right)^{\frac{\beta}{2\beta+1}}\right)\geq 1- n^{-{(\kappa-1)^2\mu_2 }/({4K^2_{\max}})} \log n,
\end{equation*}
where the last inequality follows from the fact that $M \equiv 0.8 \log_4 n +0.2 \log_4 \log n\leq \log n$ for $n\ge 2$. Note that this inequality holds for all $\mathcal{A}\subseteq [n]$ and all $f\in \Sigma(\beta,L)$. Taking infimum, we have proven the statement. 
\end{proof}

\begin{proof}[Proof of Lemma \ref{theorem: Lepskii complexity}]
{\sc Part}~\ref{i:Lepskii compleixity:a}. 
For any fixed \( m \in [M] \), computing \( (\hat{f}_{n,m}(x_i))_{i\in \mathcal{A}} \) requires at most \( O(|\mathcal{A}| \cdot n h_m) \) operations.  Since \( (\hat{f}_{n,k}(x_i))_{i\in \mathcal{A}} \) has already been computed and stored in memory for all \( k \in [m-1] \), the additional computational cost incurred by executing lines 3–6 in the CALM  (\cref{Lepskii}) is \( O(|\mathcal{A}| \cdot m) \). Thus, the total computational complexity is upper bounded by
\begin{equation*}
    \sum_{m\in [M]} \bigl(|\mathcal{A}| \cdot (n h_m + m)\bigr) = O\left(|\mathcal{A}| \cdot n h_M\right) = O\left(|\mathcal{A}| \cdot n^{\frac{4}{5}} (\log n)^{\frac{1}{5}}\right).
\end{equation*}

{\sc Part}~\ref{i:Lepskii compleixity:b}. Recall from \eqref{e:defn E kappa} that 
\begin{align*}
    \Pi(\omega)\coloneqq \underset{m=1,\dots,M}{\max}\frac{\rho_n(m)}{G_2(m)} \quad\text{and}\quad
    E_{\kappa} \coloneqq \{\omega\,|\, \Pi(\omega) \geq \kappa\}.
\end{align*}
By \eqref{ieq: E_kappa}, we see that $\mathbb{P}(E_{\kappa})\leq M n^{-\mu_2(\kappa-1)^2 /(4K^2_{\max})}$. 

We consider the case of $h_{\Bar{m}}\gg (\log (n) /n)^{1 /(2\beta+1)}$. Then, on $E_{\kappa}^c$,
\begin{equation*}
    \rho_{\mathcal{A}}(m)\leq\rho_n(m) \leq \kappa G_2(m) = \kappa C_{\rho} \sqrt{\frac{\log n}{n h_m}},\quad m\in [M].
\end{equation*}
Since $h_{\Bar{m}}\gg (\log( n) /n)^{1/(2\beta+1)}$, we obtain 
\begin{equation}
    \label{ieq: random term}
    \rho_{\mathcal{A}}(\Bar{m})\leq \rho_n(\Bar{m})\leq\kappa C_{\rho} \sqrt{\frac{\log n}{n h_{\Bar{m}}}} \ll \left(\frac{\log n}{n}\right)^{\frac{\beta}{2\beta+1}}.
\end{equation}
Furthermore, since $f\in \mathcal{C}_{\mathcal{A},\beta}$, there exists $c_{\mathcal{A}}>0$ such that
\begin{equation}
    \label{ieq: deterministic term}
    B_{f,\mathcal{A}}(\Bar{m}) \geq c_{\mathcal{A}} h_{\Bar{m}}^{\beta}\gg \left(\frac{\log n}{n}\right)^{\frac{\beta}{2\beta+1}}. 
\end{equation}
Combining \eqref{ieq: random term} and \eqref{ieq: deterministic term}, we see that $d_{\mathcal{A}}(f,\hat{f}_{n,\Bar{m}})$ is dominated by the bias term and 
\begin{equation*}
    d_{\mathcal{A}}(f,\hat{f}_{n,\Bar{m}}) \gtrsim B_{f,\mathcal{A}}(\Bar{m}) \gg \left(\frac{\log n}{n}\right)^{\frac{\beta}{2\beta+1}}.
\end{equation*}
Therefore,
\begin{multline*}
\mathbb{P}\left(h_{\Bar{m}}\gg \left(\frac{\log n}{n}\right)^{\frac{1}{2\beta+1}}\right) 
    = \mathbb{P}\left(h_{\Bar{m}}\gg \left(\frac{\log n}{n}\right)^{\frac{1}{2\beta+1}},\, E_{\kappa} \right)+ \mathbb{P}\left( h_{\Bar{m}}\gg \left(\frac{\log n}{n}\right)^{\frac{1}{2\beta+1}},\, E_{\kappa}^c \right)\\
    \leq \mathbb{P}(E_{\kappa}) +\mathbb{P}\left(d_{\mathcal{A}}(f,\hat{f}_{n,\Bar{m}})\gg \left(\frac{\log n}{n}\right)^{\frac{\beta}{2\beta+1}}  \right)\leq  2Mn^{-\frac{\mu_2 }{4K^2_{\max}}(\kappa-1)^2},
\end{multline*}
where the last inequality follows from \eqref{ieq: E_kappa} and \cref{theorem: Lepskii}. Equivalently,
\begin{equation*}
    \mathbb{P}\left( \abs{\mathcal{A}} \cdot 4^{\Bar{m}} \lesssim \abs{\mathcal{A}}\cdot n^{\frac{2\beta}{2\beta+1}} \left( \log n\right)^{\frac{1}{2\beta+1}}\right)\geq 1- 2Mn^{-\frac{\mu_2 }{4K^2_{\max}}(\kappa-1)^2}\geq 1-2n^{-\frac{\mu_2 }{4K^2_{\max}}(\kappa-1)^2} \log n,
\end{equation*}
for sufficiently large $n$.
\end{proof}

\begin{proposition}
\label{p:late stopping}
Under the nonparametric regression model \eqref{model}, suppose that Assumptions \hyperref[M2]{(M2)} and \ref{K1}--\ref{K3} hold. Then with probability at least $1-2n^{-\mu_2(\kappa-1)^2/(4K^2_{\max})}\cdot (\log n)^2$, CALM (\cref{Lepskii}) terminates in at least 
\begin{equation*}
    O\left(\abs{\mathcal{A}}\cdot n^{\frac{2\beta}{(2\beta+1)}}\left(\log n\right)^{\frac{1}{2\beta+1}}\right)
\end{equation*}
steps.

\end{proposition}
\begin{proof}
Let 
\begin{equation}
\label{e:mstar}
    m^* \coloneqq \max\{m\,|\, G_1(m)\leq G_2(m),\,m = 1,\dots,M\}.
\end{equation}
Then, by straightforward computation, we have
$$4^{m^*} \asymp n^{\frac{2\beta}{2\beta+1}} (\log n)^{\frac{1}{2\beta +1}}.$$
Note that $G_1$ is increasing while $G_2$ is decreasing. We have for all $k,m\in [m^*]$ and $\kappa>1$,
\begin{equation*}
    B_{f,\mathcal{A}}(m) \leq G_1(m) \leq G_1(m^*) \leq G_2(m^*)\leq G_2(k) < \kappa G_2(k).
\end{equation*}
Thus, for any $m\in [m^*-1]$,
\begin{align*}
    \mathbb{P}(\Bar{m} = m) 
   \leq \, &\mathbb{P}\left( d_{\mathcal{A}}(\hat{f}_{n,k},\hat{f}_{n,m+1})>4\kappa G_2(k) \text{ for some } k\in[m]\right) \\
    \leq \, &\sum_{k\in [m]}  \mathbb{P}\left(d_{\mathcal{A}}(\hat{f}_{n,k},\hat{f}_{n,m+1})>4\kappa G_2(k)\right) \\
    \leq \, &\sum_{k\in [m]}  \mathbb{P}\left(\max\bigl\{d_{\mathcal{A}}(\hat{f}_{n,k},f), d_{\mathcal{A}}(\hat{f}_{n,m+1},f)\bigr\}>2\kappa G_2(k) \right) \\
    \leq \, &\sum_{k\in [m]}  \mathbb{P}\Big(\max\bigl\{B_{f,\mathcal{A}}(k)+\rho_{\mathcal{A}}(k),  B_{f,\mathcal{A}}(m+1)+\rho_{\mathcal{A}}(m+1)\bigr\}>2\kappa G_2(k)\Big) \\
    \leq \, &\sum_{k\in [m]}  \left(\mathbb{P}\Bigl(B_{f,\mathcal{A}}(k)+\rho_{\mathcal{A}}(k) >2\kappa G_2(k)\Bigr) + \mathbb{P}\Bigl(B_{f,\mathcal{A}}(m+1)+\rho_{\mathcal{A}}(m+1) >2\kappa G_2(k)\Bigr)\right)\\
    \leq \, &\sum_{k\in [m]}  \left(\mathbb{P}\Bigl(\rho_{\mathcal{A}}(k) >\kappa G_2(k)\Bigr) + \mathbb{P}\Bigl(\rho_{\mathcal{A}}(m+1) > \kappa G_2(k)\Bigr)\right)\\
    \leq \,& \sum_{k\in [m]} \Big(n^{-\frac{\mu_2}{4K^2_{\max}}(\kappa-1)^2} + n^{-\frac{\mu_2}{4K^2_{\max}} (2^{m+1-k}\kappa-1)^2}\Big)\leq \,2m n^{-\frac{\mu_2}{4K^2_{\max}}(\kappa-1)^2}.
\end{align*}
Consequently,
\begin{equation*}
    \mathbb{P}(\Bar{m} < m^* )\leq \sum_{m\in [m^*-1]}2m n^{-\frac{\mu_2(\kappa-1)^2}{4K^2_{\max}}} \leq M^2 n^{-\frac{\mu_2(\kappa-1)^2}{4K^2_{\max}}(\kappa-1)^2} \leq n^{-\frac{\mu_2(\kappa-1)^2}{4K^2_{\max}}} (\log n)^{2}.
\end{equation*}
Equivalently, we have
\begin{equation*}
    \mathbb{P}\left(4^{\Bar{m}} \geq 4^{m^*}\right) = \mathbb{P}\left( 4^{\Bar{m}} \gtrsim n^{\frac{2\beta}{2\beta+1}} (\log n)^{\frac{1}{2\beta+1}} \right) \geq 1- n^{-{(\kappa-1)^2\mu_2}/{(4K^2_{\max}})} (\log n)^{2}.
\end{equation*} 
Note that CALM terminates in $O(\abs{\mathcal{A}}\cdot 4^{\Bar{m}})$ steps, so we have proven the statement.
\end{proof}

\begin{remark}
\label{r:late stopping}
By definition, the parameter $m^*$ introduced in \eqref{e:mstar} depends solely on $G_1$, $G_2$, and the constant $\kappa$, and is independent of both $f$ and $\mathcal{A}$. Consequently, there exists a constant $C_{m^*} > 0$, independent of $f$ and $\mathcal{A}$, such that
\[
h_{m^*} =  \frac{4^{m^*-1}}{n} = C_{m^*}\left(\frac{\log n}{n}\right)^{\frac{1}{2\beta+1}}.
\]
By applying Proposition \ref{p:late stopping}, we obtain that the CALM procedure yields a bandwidth $h_{\Bar{m}} = h_{\Bar{m}}(f,\mathcal{A})$ satisfying
\[
\inf_{\mathcal{A} \subseteq [n]} \inf_{f \in \Sigma(\beta,L)} \mathbb{P}\Biggl( h_{\Bar{m}}(f,\mathcal{A}) \geq C_{m^*} \left(\frac{\log n}{n}\right)^{\frac{1}{2\beta+1}} \Biggr)
\ge 1 - n^{-(\kappa-1)^2{\mu_2}/{(4K^2_{\max})}} (\log n)^{2},
\]
for all $\beta \in (0,2]$ and all sufficiently large $n$.
\end{remark}

\subsubsection*{Proofs for \texorpdfstring{\cref{SS: Adaptive MC}}{SS: Adaptive MC}}

\begin{proof}[Proof of \cref{theorem: Lepskii consistency}]
Note that $\kappa> 1+2K_{\max}/\sqrt{\mu_2}$, and
\begin{align*}
    \Delta_{\beta,n} =& \left(\frac{\log n}{n}\right)^{\frac{\beta -\ceil{\beta}+1}{2\beta+1}},\\
    \gamma_n =&
    \begin{cases}
        \left(2C_{\max}+\sigma \sqrt{\frac{4W}{C_{m^*}}}+1\right) \Delta_{\beta,n}, &\quad \text{ for } \beta \in (0,1],\\
        \left(4LC_{h}^{\beta-1} + \sigma \sqrt{\frac{4W}{C^3_{m^*}}}+\sigma \sqrt{\frac{6W}{C^3_{m^*}}}\right)\Delta_{\beta,n}, &\quad \text{ for } \beta \in (1,2],
    \end{cases}\\
    \text{and} \quad C =& 
    \begin{cases}
        4C_{\max} + \sigma \sqrt{\frac{4W}{C_{m^*}}} +2LC_{m^*}^{\beta} + 2,&\quad \text{ for } \beta \in (0,1],\\
        2\sigma \sqrt{\frac{4W}{C^3_{m^*}}}+2\sigma \sqrt{\frac{6W}{C^3_{m^*}}}+10LC_{m^*}^{\beta-1},&\quad \text{ for } \beta \in (1,2]
    \end{cases}
\end{align*}
with $C_{\max}$ in \cref{theorem: Lepskii}, $C_{m^*}$ in Remark \ref{r:late stopping}, $W$ in \cref{Theorem: Bounds of variances} and $C_h$ in \eqref{defn Ch}. Let $(f_n)_{n\in \mathbb{N}}$ be a sequence of functions such that  $f_n\in \mathcal{F}_{\beta}(C\Delta_{\beta,n})$ for all $n \in \mathbb{N}$. Then by Corollary~\ref{p:>=hn}, we have $\varepsilon_{\ceil{\beta}-1,\gamma_n}(f_n)\geq h_{m^*}$. Consequently, for sufficiently large $n$, it holds that
\begin{equation*}
    \varepsilon_n^{-1}\coloneqq  \Ceil{-2\log\left(\frac{\alpha}{2}\right)\cdot \bigl(\varepsilon_{\ceil{\beta}-1,\gamma_n}(f_n)\bigr)^{-1}} \leq \Ceil{-2\log\left(\frac{\alpha}{2}\right) h_{m^*}^{-1}} \leq C_{n}(\alpha).
\end{equation*}
Let $I_1,\dots, I_{C_{n}(\alpha)}$ denote the repeatedly generated indices via uniform sampling in A-FOMT (\cref{alg:MC Lepskii}). For each $l \in [C_{n}(\alpha)]$, let $\mathcal{J}_l^+$ and $\mathcal{J}_l^-$ be the sets of generated random indices $J_k$ and $J_k'$ for left and right searches starting from $I_l$, respectively, and define the set of the pairs generated from $I_l$ as 
\begin{align*}
    \mathcal{P}_l := \left\{(I_l,I_l+J)\,\big|\, J\in \mathcal{J}_l^{+}\right\}\cup \left\{(I_l-J',I_l)\,\big|\, J'\in \mathcal{J}_l^{-}\right\}.
\end{align*}
Note that $\mathcal{P}_l$'s are i.i.d.\ distributed for all $l\in [C_{n}(\alpha)]$. Furthermore, define 
\begin{equation*}
    \mathcal{P} := \bigcup_{l \in [C_{n}(\alpha)]}\mathcal{P}_l,\quad \text{ and } \quad \mathcal{P}_{\varepsilon_n^{-1}} := \bigcup_{l \in [\varepsilon_n^{-1}]}\mathcal{P}_l,
\end{equation*}
i.e., $\mathcal{P}$ contains of all pairs in A-FOMT, while $\mathcal{P}_{\varepsilon_n^{-1}}$ consists of only pairs generated in the first $\varepsilon_n^{-1}$ rounds. Clearly, $\mathcal{P}_{\varepsilon_n^{-1}} \subseteq \mathcal{P}$. 
Moreover, for each $l\in [C_{n}(\alpha)]$, let $\mathcal{A}_l$ be the set of generated indices contained in $\mathcal{P}_l$. For each $n\in \mathbb{N}$ and $l\in [C_{n}(\alpha)]$, we define 
\begin{equation*}
    F\equiv  F(n) :=\bigcup_{l\in [C_{n}(\alpha)]} F_l,\,\, \text{ with } \,\, F_l \equiv F_{l,n}:= \left\{h_{\Bar{m}}(f_n, \mathcal{A}_l) < h_{m^*} \equiv C_{m^*}\left(\frac{\log n}{n}\right)^{\frac{1}{2\beta+1}}\right\}.
\end{equation*}
It follows from Remark \ref{r:late stopping} that 
\begin{equation}
    \label{ieq:late stopping}
    \mathbb{P}_{f_n}(F)\leq \sum_{l\in [C_{n}(\alpha)]} \mathbb{P}_{f_n}(F_l)\leq C_{n}(\alpha) n^{-\frac{(\kappa-1)^2\mu_2}{4K^2_{\max}}} (\log n)^{2} \le n^{1-\frac{(\kappa-1)^2\mu_2}{4K^2_{\max}}} (\log n)^{2},
\end{equation} 
which tends to zero as $n \to \infty$, provided $\kappa> 1+2K_{\max}/\sqrt{\mu_2}$. 
We now proceed to prove the statements for $\beta \in(0,1]$ and $\beta \in (1,2]$, separately.

\medskip

\textbf{Case $\beta \in (0,1]$.}
We define the events $E$ and $E_l$ as follows:
    \begin{equation*}
        E = E(n)= \bigcup_{l\in [C_{n}(\alpha)]} E_l,\,\,\, \text{ with } \,\,\, E_l = E_{l,n}= \left\{d_{\mathcal{A}_l}(\hat{f}_{n,\Bar{m}}, f_n) > C_{\max} \Delta_{\beta,n}\right\}.
    \end{equation*}
    Applying \cref{theorem: Lepskii}, we obtain the upper bound 
    \begin{equation} 
    \label{ieq:dA upper bound} 
    \mathbb{P}_{f_n}(E)\leq \sum_{l\in [C_{n}(\alpha)]} \mathbb{P}_{f_n}(E_l) \leq C_{n}(\alpha) \cdot n^{-\frac{\mu_2}{4K^2_{\max}}(\kappa-1)^2} \log n\leq n^{1-\frac{\mu_2}{4K^2_{\max}}(\kappa-1)^2} \log n. 
    \end{equation}
    Following a similar approach as in Corollary \ref{c:unknown sigma}, we decompose 
    \begin{align} 
    \label{e: p decom3} \mathbb{P}_{f_n}(\Phi_A = 0) &= \mathbb{P}_{f_n}(\Phi_{i,j} = 0, \text{ for all }(i,j)\in \mathcal{P}) \nonumber \\
    &\leq \mathbb{P}_{f_n}(\Phi_{i,j} = 0, \text{ for all } (i,j)\in \mathcal{P}_{\varepsilon_n^{-1}}) \nonumber \\ 
    &\leq \mathbb{P}_{f_n}(E) + \mathbb{P}_{f_n}(F) + \mathbb{P}_{f_n}(\Phi_{i,j} = 0, \text{ for all } (i,j)\in \mathcal{P}_{\varepsilon_n^{-1}}, E^c\cap F^c). 
    \end{align}
    From \eqref{ieq:late stopping} and \eqref{ieq:dA upper bound}, we conclude that both $\P_{f_n}(E)$ and $\P_{f_n}(F)$ vanish as $n\to \infty$, since $\kappa> 1+2K_{\max}/\sqrt{\mu_2}$. 
   Next, we define 
    \begin{equation*}
        \mathcal{I} = \mathcal{I}_n = \left\{(i,j)\,\bigg|\, f_n(x_i)-f_n(x_j)\geq \left(2C_{\max} +\sigma \sqrt{\frac{4W}{C_{m^*}}}\right)\Delta_{\beta,n}\right\},
    \end{equation*}
    and further decompose $\mathbb{P}_{f_n}(\Phi_{i,j} = 0,\,\text{ for all } (i,j)\in \mathcal{P}_{\varepsilon_n^{-1}},E^c\cap F^c)$ as follows:
    \begin{align}
        \label{eq:p Phi=0 EcFc3}
        &\,\mathbb{P}_{f_n}(\Phi_{i,j} = 0,\,\text{ for all } (i,j)\in \mathcal{P}_{\varepsilon_n^{-1}},\,E^c \cap F^c) \nonumber \\ 
        =&\,\mathbb{P}_{f_n}(\Phi_{i,j} = 0,\,\text{ for all } (i,j)\in \mathcal{P}_{\varepsilon_n^{-1}},\, \mathcal{P}_{\varepsilon_n^{-1}}\cap \mathcal{I} \neq \emptyset, \,E^c \cap F^c) \nonumber\\
        &\,+ \mathbb{P}_{f_n}(\Phi_{i,j} = 0,\,\text{ for all } (i,j)\in \mathcal{P}_{\varepsilon_n^{-1}},\, \mathcal{P}_{\varepsilon_n^{-1}}\cap \mathcal{I} = \emptyset,\,E^c \cap F^c)
    \end{align}
    For $\mathcal{P}_{\varepsilon_n^{-1}}\cap \mathcal{I}\neq\emptyset$, let $(i,j)$ be an arbitrary pair in $\mathcal{P}_{\varepsilon_n^{-1}}\cap \mathcal{I}$. Together with events $E^c$ and $F^c$ we have
    \begin{align*}
        \,\,&\underset{h\geq h_{m^*}}{\min} \Bigl(T_{i,j}(h) - C_{n,\alpha,i,j}(h)\Bigr)\\
        =& \underset{h\geq h_{m^*}}{\min}\Bigl(\hat{f}_{n}(x_i;h) - f_n(x_i) +f_n(x_i)-f_n(x_j) + f_n(x_j)- \hat{f}_{n}(x_j;h) - C_{n,\alpha,i,j}(h) \Bigr)\\
        \geq& -C_{\max}\cdot \Delta_{\beta,n} + \left(2C_{\max} +\sigma \sqrt{\frac{4W}{C_{m^*}}}\right)\cdot \Delta_{\beta,n} -C_{\max}\cdot \Delta_{\beta,n} - \underset{h\geq h_{m^*}}{\max} C_{n,\alpha,i,j}(h) \\
        \geq&\,\sigma\sqrt{\frac{4W}{C_{m^*}}}\left(\frac{\log n}{n}\right)^{\frac{\beta}{2\beta+1}}- C_{n,\alpha,i,j}(h_{m^*})  \geq 0,
    \end{align*}
    where $C_{n,\alpha,i,j}(h)$ is given in \eqref{C,n,alpha,I,J} with $h_n$ replaced by $h$. The first inequality follows from $E^c$ and $F^c$, and the last inequality can be proved by using the same technique as in Lemma~\ref{l: upper bound of C n,alpha}. Thus, the first probability in \eqref{eq:p Phi=0 EcFc3} is
    \begin{equation*}
        \mathbb{P}_{f_n}(\Phi_{i,j} = 0,\, \text{ for all }(i,j)\in \mathcal{P}_{\varepsilon_n^{-1}}\cap \mathcal{I},\,E^c\cap F^c) = 0.
    \end{equation*}
    For the second probability in \eqref{eq:p Phi=0 EcFc3}, noting that $\mathcal{P}_l$ are i.i.d.\ random sets, we derive 
    \begin{equation}
        \label{ieq:ub PIeEEc3}
        \mathbb{P}_{f_n}(\Phi_{i,j} = 0,\text{ for all }(i,j)\in\mathcal{P}_{\varepsilon_n^{-1}},\, \mathcal{P}_{\varepsilon_n^{-1}}\cap \mathcal{I} = \emptyset,\,E^c\cap F^c)
        \leq \,\mathbb{P}_{f_n}(\mathcal{P}_l\cap \mathcal{I} = \emptyset)^{\varepsilon_n^{-1}}.    
    \end{equation}
    For each $l\in [C_{n}(\alpha)]$, it follows from Lemma~\ref{l:epsilon_f}~\ref{l:eps_g:b} that
    \begin{align}
        \label{ieq:ub PlIeE3}
        \mathbb{P}_{f_n}(\mathcal{P}_l\cap \mathcal{I}= \emptyset) \leq \frac{1}{n}\sum_{i:\, x_i \in H_{f_n}^{1/n}(\gamma_n)}\mathbb{P}_{f_n}(\mathcal{P}_l\cap \mathcal{I}= \emptyset\,\big| \, I_l = i) + 1-\varepsilon_{0,\gamma_n}(f_n),
    \end{align}
    where $H_{f_n}(\gamma_n)$ is the set of $\gamma_n$-heavy points of $f_n$. Note that, for $x_i \in H_{f_n}^{1/n}(\gamma_n)$, there exists $a\in H_{f_n}(\gamma_n) =H_{f_n,R}(\gamma_n)\cup H_{f_n,L}(\gamma_n)$ satisfying $\abs{x_i-a}\leq 1/n$. Consider first $a\in H_{f_n,R} $, i.e., there exists $b\in (a,1]$ such that $\lambda(A)\geq (b-a)/2$, where 
    \begin{align*}
        A = \left\{x\in [a,b]\,\Big|\, f_n(a)-f_n(x)\geq \gamma_n \right\}.
    \end{align*}
    Using Lemma~\ref{l: 1/8}, we obtain
    \begin{equation*}
        \mathbb{P}_{f_n}\left(\text{ there exists }j\in \mathcal{J}_l^{+} \text{ satisfying }x_{i+j} \in A^{1/n}\right) \geq 1- \frac{1}{n}.
    \end{equation*}
    For any $j\in \mathcal{J}_l^{+}$ satisfying $x_{i+j} \in A^{1/n}$, there exists $c\in A$ fulfilling $\abs{c-x_{i+j}}\leq 1/n$ and $f_n(a)-f_n(c)\geq \gamma_n$. Then, we have 
    \begin{align*}
        f_n(x_i)-f_n(x_j)=& f_n(x_i)-f_n(a) + f_n(a)- f_n(c)+f_n(c)- f_n(x_j)\\
        \geq& \,\gamma_n -2L \left(\frac{1}{n}\right)^{\beta}\\
        \geq& \, \left(2C_{\max}+\sigma\sqrt{\frac{4W}{C_{\min}}}\right) \cdot \left(\frac{\log n}{n}\right)^{\frac{\beta}{2\beta+1}},
    \end{align*}
    which demonstrates that $(i,i+j)\in \mathcal{I}$. Thus, given $I_l = i$ and the $\gamma_n$-right-heaviness of $x_i$,
    \begin{equation*}
        \mathbb{P}_{f_n}(\mathcal{P}_l\cap \mathcal{I}\neq \emptyset\,\big| \, I_l = i) \geq \mathbb{P}_{f_n}\left(\text{ there exists }j\in \mathcal{J}_l^{+} \text{ satisfying }x_{i+j} \in A^{1/n}\right) \geq 1-\frac{1}{n}, 
    \end{equation*}
    and this implies that for all $\gamma_n$-right-heavy point $x_i$,
    \begin{equation}
        \label{ieq:ub PlIeE Il3}
        \mathbb{P}_{f_n}(\mathcal{P}_l\cap \mathcal{I}=  \emptyset\,\big| \, I_l = i)\leq \frac{1}{n}.
    \end{equation}
    By symmetry, we can achieve the same upper bound for any $\gamma_n$-left-heavy point $x_i$. 
    
    Combining \eqref{ieq:ub PIeEEc3}, \eqref{ieq:ub PlIeE3} and \eqref{ieq:ub PlIeE Il3}, we obtain 
    \begin{align}
        \label{ieq:ub p PIeE Ec3}
         & \limsup_{n\to \infty}\mathbb{P}_{f_n}(\Phi_{i,j} = 0,\text{ for all }(i,j)\in\mathcal{P}_{\varepsilon_n^{-1}},\, \mathcal{P}_{\varepsilon_n^{-1}}\cap \mathcal{I} = \emptyset,\,E^c\cap F^c) \nonumber \\
         \leq&\,\limsup_{n\to \infty}\mathbb{P}_{f_n}(\mathcal{P}_l\cap \mathcal{I} = \emptyset)^{C_{n}(\alpha)} \nonumber\\
         \leq& \,\limsup_{n\to \infty}\left( \frac{1}{n}\sum_{i:\, x_i \in H_{f_n}^{1/n}(\gamma_n)}\mathbb{P}_{f_n}(\mathcal{P}_l\cap \mathcal{I}= \emptyset\,\big| \, I_l = i) + 1-\varepsilon_{0,\gamma_n}(f_n)\right)^{\varepsilon_n^{-1}}\nonumber\\
         \leq& \,\limsup_{n\to \infty}\left( \frac{1}{n}\cdot n \cdot \max_{i:\, x_i \in H_{f_n}^{1/n}(\gamma_n)}\mathbb{P}_{f_n}(\mathcal{P}_l\cap \mathcal{I}= \emptyset\,\big| \, I_l = i) + 1-\varepsilon_{0,\gamma_n}(f_n)\right)^{\ceil{-2\log\left(\frac{\alpha}{2}\right)\cdot (\varepsilon_{0,\gamma_n}(f_n))^{-1}}} \nonumber\\ 
        \leq& \,\limsup_{n\to \infty}\left(\frac{1}{n} + 1-\varepsilon_{0,\gamma_n}(f_n)\right)^{\ceil{-2\log\left(\frac{\alpha}{2}\right)\cdot (\varepsilon_{0,\gamma_n}(f_n))^{-1}}} \nonumber \\
        \leq& \,\limsup_{n\to \infty}\left(  1-\frac{1}{2}\varepsilon_{0,\gamma_n}(f_n)\right)^{\ceil{-2\log\left(\frac{\alpha}{2}\right)\cdot (\varepsilon_{0,\gamma_n}(f_n))^{-1}}} \leq \frac{\alpha}{2}.
    \end{align}
    Summarizing \eqref{ieq:late stopping}--\eqref{ieq:ub PIeEEc3} and \eqref{ieq:ub p PIeE Ec3}, we see that for all $\beta \in (0,1]$
    \begin{equation*}
        \limsup_{n\to \infty}\mathbb{P}_{f_n}(\Phi_{A}=0)\leq \limsup_{n\to \infty}\mathbb{P}_{f_n}(\Phi_{i,j} = 0,\,\text{ for all } (i,j)\in \mathcal{P}_{\varepsilon_n^{-1}}) \leq \alpha,
    \end{equation*}
    i.e., with probability at least $1-\alpha$, A-FOMT can detect at least a violation of $f_n$ with only local tests depending on $O\bigl((\varepsilon_{0,\gamma_n}(f_n))^{-1}\bigr)$ uniformly distributed random indices $I_l$.
    
\medskip

\textbf{Case $\beta\in(1,2]$.} We define events $G$ and $G_{i,l}$ as follows:
    \begin{align*}
        G = \bigcup_{i\in[n-1],\, l\in [C_{n}(\alpha)]}G_{i,l}\qquad
        \text{with}\quad G_{i,j} = G_{i,j,n}= \left\{R_{i,i+1}(h_{\Bar{m}}(f_n,\mathcal{A}_l))\leq -\sigma \sqrt{\frac{6W}{C^3_{m^*}}} \frac{1}{n}\Delta_{\beta,n}\right\}.
    \end{align*}
    For any fixed $G_{i,l}$, we have
    \begin{align*}
        \mathbb{P}_{f_n}(G_{i,l}\cap F^c) =& \,\mathbb{P}_{f_n}\left( h_{\Bar{m}}(f_n,\mathcal{A}_l)\geq h_{m^*},\, R_{i,i+1}(h_{\Bar{m}}(f_n,\mathcal{A}_l))\leq -\sigma \sqrt{\frac{6W}{C^3_{m^*}}} \frac{1}{n}\Delta_{\beta,n} \right) \\
        \leq& \, \underset{h\geq h_{m^*}}{\max} \mathbb{P}_{f_n}\left(R_{i,i+1}(h)\leq -\sigma \sqrt{\frac{6W}{C^3_{m^*}}} \frac{1}{n}\Delta_{\beta,n} \right)\\
        \leq& \exp \left( -\frac{1}{2} \frac{\sigma^2\frac{6W}{C^3_{m^*}}\cdot \frac{1}{n^2} \Delta^2_{\beta,n}}{\underset{h\geq h_{m^*}}{\max}\V(R_{i,i+1}(h))} \right)
        \leq \exp \left( -\frac{1}{2} \frac{\sigma^2\frac{6W}{C^3_{m^*}}\cdot \frac{1}{n^2} \Delta^2_{\beta,n}}{\sigma^2 W n^{-3}h^{-3}_{m^*}}\right)
        \leq  \exp \left(-3\log n\right) = n^{-3}.
    \end{align*}
    Therefore, for sufficiently large $n$,
    \begin{equation*}
        \mathbb{P}_{f_n}(G\cap F^c) \leq \sum_{i\in[n-1],\, l\in [C_{n}(\alpha)]} \mathbb{P}_{f_n}(G_{i,l}\cap F^c) <n^{-1}.
    \end{equation*}
    Similar to \eqref{e: p decom3}, we obtain 
    \begin{align}
        \label{e: p decom4}
        \mathbb{P}_{f_n}(\Phi_A = 0) =& \, \mathbb{P}_{f_n}(\Phi_{i,j} = 0,\,\text{ for all } (i,j)\in \mathcal{P}) \nonumber \\
        \leq& \, \mathbb{P}_{f_n}(\Phi_{i,j} = 0,\,\text{ for all } (i,j)\in \mathcal{P}_{\varepsilon_n^{-1}}) \nonumber \\
        \leq&\, \mathbb{P}_{f_n}(F) + \mathbb{P}_{f_n}(G\cap F^c) + \mathbb{P}_{f_n}(\Phi_{i,j} = 0,\,\text{ for all } (i,j)\in \mathcal{P}_{\varepsilon_n^{-1}},G^c\cap F^c).
    \end{align}
    Clearly, both $\mathbb{P}_{f_n}(G\cap F^c)$ and $\mathbb{P}_{f_n}(F)$ converge to zero, as $n \to \infty$. We further decompose $\mathbb{P}_{f_n}(\Phi_{i,j} = 0,\,\text{ for all } (i,j)\in \mathcal{P} _{\varepsilon_n^{-1}},\,G^c \cap F^c)$ into
    \begin{align}
        \label{eq:p Phi=0 EcFc4}
        &\,\mathbb{P}_{f_n}(\Phi_{i,j} = 0,\,\text{ for all } (i,j)\in \mathcal{P}_{\varepsilon_n^{-1}},\,G^c \cap F^c) \nonumber \\ 
        =&\,\mathbb{P}_{f_n}(\Phi_{i,j} = 0,\,\text{ for all } (i,j)\in \mathcal{P}_{\varepsilon_n^{-1}},\, \mathcal{P}_{\varepsilon_n^{-1}}\cap \mathcal{I} \neq \emptyset, \,G^c \cap F^c) \nonumber\\
        &\,+ \mathbb{P}_{f_n}(\Phi_{i,j} = 0,\,\text{ for all } (i,j)\in \mathcal{P}_{\varepsilon_n^{-1}},\, \mathcal{P}_{\varepsilon_n^{-1}}\cap \mathcal{I} = \emptyset,\,G^c \cap F^c),
    \end{align}
    with 
    \begin{align*}
        \mathcal{I} = \mathcal{I}_n =  \bigg\{ (i,j)\;\bigg|\; h_{m^*} \leq x_i <x_{i+1}\leq 1-h_{m^*},\,
        \underset{h\geq h_{m^*}}{\min} D_{i,i+1}(h)-C_{n,\alpha,i,i+1}(h) - \sigma \sqrt{\frac{6W}{C_{m^*}}}\cdot \frac{1}{n}\Delta_{\beta,n} \geq 0 \bigg\}.
    \end{align*}
    Then, for any pair $(i,j)\in \mathcal{P}_{\varepsilon_n^{-1}}\cap \mathcal{I}$, together with $F^c \cap G^c$ we have
    \begin{align*}
        \underset{h\geq h_{m^*}}{\min}D_{i,i+1}(h)+R_{i,i+1}(h)-C_{n,\alpha,i,i+1}(h) 
        \geq \underset{h\geq h_{m^*}}{\min} D_{i,i+1}(h)-C_{n,\alpha,i,i+1}(h) - \sigma \sqrt{\frac{6W}{C_{m^*}}}\cdot \frac{1}{n}\Delta_{\beta,n}\geq 0,
    \end{align*}
    which demonstrates that $$\mathbb{P}_{f_n}(\Phi_{i,j} = 0,\,\text{ for all } (i,j)\in  \mathcal{P}_{\varepsilon_n^{-1}}\cap \mathcal{I}, \,G^c \cap F^c)=0.$$
    For the second term in \eqref{eq:p Phi=0 EcFc4}, applying the same technique in the proof of \cref{t:consistency1}, we can show, 
    for any $l\in [\varepsilon_n^{-1}]$,
    \begin{equation*}
        \mathbb{P}_{f_n}(\mathcal{P}_l\cap \mathcal{I}\neq \emptyset) \leq 1+\frac{1}{n}-\varepsilon_{1,\gamma_n}(f_n)\leq 1-\frac{1}{2}\varepsilon_{1,\gamma_n}(f_n).
    \end{equation*}
    Therefore, 
    \begin{align}
        \label{ieq:ub PIeEEc4}
        &\,\mathbb{P}_{f_n}(\Phi_{i,j} = 0,\text{ for all }(i,j)\in\mathcal{P}_{\varepsilon_n^{-1}},\, \mathcal{P}_{\varepsilon_n^{-1}}\cap \mathcal{I} = \emptyset,\,E^c\cap F^c) \nonumber\\
        \leq& \,\mathbb{P}_{f_n}(\mathcal{P}_l\cap \mathcal{I} = \emptyset)^{\varepsilon_n^{-1}} \nonumber\\
        \leq& \left(1-\frac{1}{2}\varepsilon_{1,\gamma_n}(f_n)\right)^{\ceil{-2\log \left(\frac{\alpha}{2}\right)\cdot(\varepsilon_{1,\gamma_n}(f_n))^{-1}}}\leq \frac{\alpha}{2}. 
    \end{align}
    It follows from \eqref{e: p decom4}--\eqref{ieq:ub PIeEEc4} that for $\beta \in (1,2]$,
    \begin{equation*}
        \limsup_{n\to \infty}\mathbb{P}_{f_n}(\Phi_{A}=0)\leq \limsup_{n\to \infty}\mathbb{P}_{f_n}(\Phi_{i,j} = 0,\,\text{ for all } (i,j)\in \mathcal{P}_{\varepsilon_n^{-1}}) \leq \alpha.
    \end{equation*}
    Equivalently, with probability at least $1-\alpha$, A-FOMT can detect at least a violation of $f_n$ with the first $O\bigl((\varepsilon_{1,\gamma_n}(f_n))^{-1} \cdot (\log n)^2 \bigr)$ local tests.
\end{proof}

\begin{proof}[Proof of \cref{t:MC Lepskii compleixity}]
{\sc Part}~\ref{i:tMC:a}. 
It follows from the proof of \cref{theorem: Lepskii consistency} that, for any $f\in \mathcal{F}_{\beta}(C\Delta_{\beta,n})$,  A-FOMT can detect a violation of $f$ by conducting local tests coming from $O\bigl((\varepsilon_{\ceil{\beta}-1,\gamma_n}(f))^{-1}\bigr)$ uniform indices $I_l$, with probability at least $1-\alpha$. For any fixed $l$, this leads to at most $O\bigl(\abs{\mathcal{A}_l}\cdot n^{4/5} (\log n)^{1/5}\bigr) = O\bigl(n^{4/5} (\log n)^{11/5}\bigr)$ steps (Lemma~\ref{theorem: Lepskii complexity}~\ref{i:Lepskii compleixity:a}). Thus, with probability at least $1-\alpha$, A-FOMT has computational complexity
\begin{equation*}
    O\left(\varepsilon^{-1}_{\ceil{\beta}-1,\gamma_n}(f) \cdot n^{\frac{4}{5}}(\log n)^{\frac{11}{5}}\right).
\end{equation*}

{\sc Part}~\ref{i:tMC:b}. 
If additionally $f\in \mathcal{C}_{\beta}$, then $f$ satisfies the self-similarity condition in \eqref{ieq:LB} on all $\mathcal{A}_l$ for $l\in [C_{n}(\alpha)]$. Thus, by Lemma~\ref{theorem: Lepskii complexity}, it holds, with probability at least $1-2n^{-\mu_2(\kappa-1)^2/(4K^2_{\max})} \log n\to 1$, as $n\to \infty$, that CALM computes estimates on $\mathcal{A}_l$ in $O\bigl(\abs{\mathcal{A}_l}\cdot n^{2\beta/(2\beta+1)} (\log n)^{1/(2\beta+1)}\bigr)$ steps for all $l\in [C_{n}(\alpha)]$. Consequently, A-FOMT requires computational complexity
\begin{equation*}
        O\left((\varepsilon_{\ceil{\beta}-1,\gamma_n}(f))^{-1} \cdot n^{\frac{2\beta}{2\beta+1}} (\log n)^{\frac{4\beta+3}{2\beta+1}}\right),
\end{equation*}
with probability at least $1-\alpha$. 

{\sc Part}~\ref{i:tMC:c}. 
In the worst case, the computational complexity of A-FOMT can be 
\begin{equation*}
    O\left(C_{n}(\alpha) \cdot n^{\frac{4}{5}}(\log n)^{\frac{11}{5}}\right) = O\left(n^{\frac{9}{5}}(\log n)^{\frac{6}{5}}\right).
\end{equation*}
\end{proof}

\algdef{SE}[REPEATN]{RepeatN}{End}[1]{\algorithmicrepeat\ #1 \textbf{times}}{\algorithmicend}
\alglanguage{pseudocode}
\begin{algorithm}[H]
\begin{algorithmic}[1]
\Statex {\textbf{Input:} sample size $n$ and index $i \in [n]$}
\State $\mathcal{P} = \emptyset,\,\mathcal{A} = \{i\}$ 
\If{$i \leq n-1$}
    \MRepeat{$\ceil{20 \log(n)}$} 
        \For{$1 \leq k \leq \ceil{\log_2(n-i)}$}
                \State Generate $J_k\sim \mathrm{Unif}([2^k\wedge(n-i)])$
                \If{$i+J_k\notin \mathcal{A}$}
                    \State $\mathcal{A} = \mathcal{A}\cup \{i+J_k\}$
                \EndIf
                \If{$(i,i+J_k)\notin \mathcal{P}$}
                    \State $\mathcal{P} = \mathcal{A}\cup \{(i,i+J_k)\}$
                \EndIf
        \EndFor
    \EndRepeat   
\EndIf

\If{$i \geq 2$}
    \MRepeat{$\ceil{20 \log(n)}$}        
    \For{$0 \leq k \leq \ceil{\log_2(i-1)}$}
                \State Generate $J_k'\sim \mathrm{Unif}([2^k\wedge(i-1)])$
                \If{$i-J_k'\notin \mathcal{A}$}
                    \State $\mathcal{A} = \mathcal{A}\cup \{i-J_k'\}$
                \EndIf
                \If{$(i-J_k',i)\notin \mathcal{P}$}
                    \State $\mathcal{P} = \mathcal{A}\cup \{(i-J_k',i)\}$
                \EndIf
        \EndFor
    \EndRepeat
\EndIf
\State \Return $(\mathcal{P},\mathcal{A})$\;
\end{algorithmic}
\caption{Indices Generator}
\label{alg:IG}
\end{algorithm}

\section{Additional materials of simulation study}
\begin{table}[ht]
\caption{Frequencies of rejections over $100$ repetitions of FOMT, A-FOMT, DS \citep{dumbgen2001multiscale}, ABD \citep{akakpo2014testing} and C \citep{chetverikov2019testing} for signals $f_i$, $i = 1, \dots, 4$. }
\label{T:signals}
\centering{\scriptsize
\begin{tabular}{| l| l|  l| l|l|l|}
\toprule
\hline \xrowht[()]{10pt}
\textbf{Sample size} & \textbf{Methods} & \textbf{$f_1$} & \textbf{$f_2$} & \textbf{$f_3$} & \textbf{$f_4$}\\[5pt]
    \hline \xrowht[()]{10pt}
    \multirow{5}{*}{$n = 400$} & FOMT   & $50\%$ & $94\%$ & $99\%$ & $70\%$ \\[5pt]
    & A-FOMT  & $41\%$ & $97\%$ & $92\%$ & $89\%$ \\[5pt]
    & DS   & $44\%$ & $87\%$ & $100\%$ & $58\%$ \\[5pt]
    & ABD   & $32\%$ & $58\%$ & $91\%$ & $44\%$ \\[5pt]
    & C   & $41\%$ & $80\%$ & $97\%$ & $80\%$ \\[5pt]
    \hline \xrowht[()]{10pt}
    \multirow{5}{*}{$n = 800$} & FOMT & $92\%$ & $100\%$ & $100\%$ & $98\%$ \\[5pt]
    & A-FOMT  & $78\%$ & $100\%$ & $97\%$ & $99\%$ \\[5pt]
    & DS  & $76\%$ & $99\%$ & $100\%$ & $97\%$ \\[5pt]
    & ABD  & $81\%$ & $96\%$ & $100\%$ & $96\%$ \\[5pt]
    & C  & $72\%$ & $95\%$ & $100\%$ & $99\%$ \\[5pt]
    \hline \xrowht[()]{10pt}
    \multirow{5}{*}{$n = 1200$} & FOMT   & $100\%$ & $100\%$ & $100\%$ & $99\%$ \\[5pt]
    & A-FOMT   & $96\%$ & $100\%$ & $100\%$ & $100\%$ \\[5pt]
    & DS  & $98\%$ & $100\%$ & $100\%$ & $100\%$ \\[5pt]
    & ABD & $83\%$ & $99\%$ & $100\%$ & $99\%$ \\[5pt]
    & C  & $94\%$ & $100\%$ & $100\%$ & $99\%$ \\[5pt]
    \hline \xrowht[()]{10pt}
    \multirow{5}{*}{$n = 1600$} & FOMT   & $100\%$ & $100\%$ & $100\%$ & $100\%$ \\[5pt]
    & A-FOMT  & $100\%$ & $100\%$ & $100\%$ & $100\%$ \\[5pt]
    & DS   & $99\%$ & $100\%$ & $100\%$ & $100\%$ \\[5pt]
    & ABD   & $100\%$ & $99\%$ & $100\%$ & $99\%$ \\[5pt]
    & C  & $95\%$ & $95\%$ & $100\%$ & $99\%$ \\[5pt]
    \hline \xrowht[()]{10pt}
    \multirow{4}{*}{$n = 2000$} &FOMT  & $100\%$ & $100\%$ & $100\%$ & $100\%$ \\[5pt]
    & A-FOMT  & $100\%$ & $100\%$ & $100\%$ & $100\%$ \\[5pt]
    & DS  & $100\%$ & $100\%$ & $100\%$ & $100\%$ \\[5pt]
    & ABD  & $100\%$ & $100\%$ & $100\%$ & $100\%$ \\[5pt]
    & C  & $99\%$ & $100\%$ & $100\%$ & $100\%$ \\[5pt]
    \hline \xrowht[()]{10pt}
    \multirow{4}{*}{$n = 2400$} &FOMT  & $100\%$ & $100\%$ & $100\%$ & $100\%$ \\[5pt]
    &A-FOMT & $100\%$ & $100\%$ & $100\%$ & $100\%$ \\[5pt]
    &DS   & $100\%$ & $100\%$ & $100\%$ & $100\%$ \\[5pt]
    & ABD   & $100\%$ & $100\%$ & $100\%$ & $100\%$ \\[5pt]
    & C  & $100\%$ & $100\%$ & $100\%$ & $100\%$ \\[5pt]
    \hline \xrowht[()]{10pt}
    \multirow{3}{*}{$n = 2800$} &FOMT   & $100\%$ & $100\%$ & $100\%$ & $100\%$ \\[5pt]
    &A-FOMT   & $100\%$ & $100\%$ & $100\%$ & $100\%$ \\[5pt]
    &DS   & $100\%$ & $100\%$ & $100\%$ & $100\%$ \\[5pt]
    & C  & $100\%$ & $100\%$ & $100\%$ & $100\%$ \\[5pt]
    \hline \xrowht[()]{10pt}
    \multirow{3}{*}{$n = 3200$} &FOMT   & $100\%$ & $100\%$ & $100\%$ & $100\%$ \\[5pt]
    &A-FOMT   & $100\%$ & $100\%$ & $100\%$ & $100\%$ \\[5pt]
    &DS  & $100\%$ & $100\%$ & $100\%$ & $100\%$ \\[5pt]
    & C   & $100\%$ & $100\%$ & $100\%$ & $100\%$ \\[5pt]
    \hline
\bottomrule
\end{tabular}
}
\end{table}

\begin{table}[ht]
\begin{minipage}{\textwidth}
\caption{The median of computational times (in seconds) of FOMT, A-FOMT, DS \citep{dumbgen2001multiscale}, ABD \citep{akakpo2014testing} and C \citep{chetverikov2019testing} over $100$ repetitions\tablefootnote{The computation times of $C$ (\citealp{chetverikov2019testing}) under $f_0$ with $n\in \{2000,2400,2800,3200\}$ are measured in one simulation.} for signals $f_i$, $i = 0, \dots, 4$.}
\label{T:time}
\centering{\scriptsize
\begin{tabular}{| l| l| S[table-format=5.3]|S[table-format=2.3]| S[table-format=2.3]|S[table-format=2.3]|S[table-format=3.3]|}
\toprule
\hline \xrowht[()]{10pt}
\textbf{Sample size} & \textbf{Methods} & \textbf{$f_0$} & \textbf{$f_1$} & \textbf{$f_2$} & \textbf{$f_3$} & \textbf{$f_4$}\\[5pt]
    \hline \xrowht[()]{10pt}
    \multirow{5}{*}{$n = 400$}  & FOMT  & 0.165        & 0.141        & 0.006        & 0.002        & 0.023 \\[5pt]
                                & A-FOMT &  0.244     & 0.238        & 0.015        & 0.008        & 0.029 \\[5pt]
                                & DS    & 0.209        & 0.198        & 0.187        & 0.152        & 0.203  \\[5pt]
                                & ABD   & 1.447        & 1.351        & 0.908        & 0.916        & 1.274 \\[5pt]
                                & C     & 13.619       & 13.443       & 0.171       & 0.087        & 0.286 \\[5pt]
    \hline \xrowht[()]{10pt}
    \multirow{5}{*}{$n = 800$}  & FOMT  & 0.259        & 0.011        & 0.003        & 0.002        & 0.012 \\[5pt]
                                & A-FOMT & 0.504      & 0.065        & 0.011        & 0.007        & 0.017 \\[5pt]
                                & DS    & 0.733        & 0.701        & 0.704        & 0.574        & 0.698 \\[5pt]
                                & ABD   & 9.057        & 2.197        & 2.230        & 2.908        & 2.798 \\[5pt]
                                & C     & 119.274      & 2.006        & 0.914        & 0.413        & 0.572 \\[5pt]
    \hline \xrowht[()]{10pt}
    \multirow{5}{*}{$n = 1200$} & FOMT  &  0.305       & 0.007        & 0.002        & 0.002        & 0.006 \\[5pt]
                                & A-FOMT & 0.859      & 0.030        & 0.010        & 0.005        & 0.015 \\[5pt]
                                & DS    & 1.585        & 4.199        & 1.653        & 1.114        & 1.469 \\[5pt]
                                & ABD   & 29.781       & 7.221        & 6.527        & 8.216        & 4.838 \\[5pt]
                                & C     & 384.963     & 2.433         & 0.892        & 0.769       & 0.615 \\[5pt]
    \hline \xrowht[()]{10pt}
    \multirow{5}{*}{$n = 1600$} & FOMT  & 0.367        & 0.021        & 0.002        & 0.002        & 0.005 \\[5pt]
                                & A-FOMT & 1.364      & 0.037        & 0.018        & 0.007        & 0.027 \\[5pt]
                                & DS    & 2.732        & 2.396        & 2.748        & 2.166        & 2.637 \\[5pt]
                                & ABD   & 75.097       & 6.772        & 10.141       & 14.358       & 7.547 \\[5pt]
                                & C     & 925.693     & 1.658         & 4.086        & 1.572       & 2.614 \\[5pt]
    \hline \xrowht[()]{10pt}
    \multirow{4}{*}{$n = 2000$} &FOMT  & 0.481         & 0.008        & 0.002        & 0.002        & 0.006 \\[5pt]
                                &A-FOMT & 1.604       & 0.021        & 0.010        & 0.007        & 0.013 \\[5pt]
                                &DS    & 4.185         & 3.684        & 3.637        & 3.112        & 3.773 \\[5pt]
                                & ABD   & 145.862       & 13.588        & 17.547       & 57.069       & 23.693 \\[5pt]
                                & C     & 2521.998     & 11.761         & 9.290        & 3.382       & 11.038 \\[5pt]
    \hline \xrowht[()]{10pt}
    \multirow{4}{*}{$n = 2400$} &FOMT  & 0.415         & 0.004        & 0.002        & 0.002        & 0.005 \\[5pt]
                                &A-FOMT & 1.854       & 0.025        & 0.012        & 0.006        & 0.016 \\[5pt]
                                &DS    & 6.744         & 5.705        & 5.887        & 4.803        & 5.583 \\[5pt]
                                & ABD   & 263.512       & 13.142        & 20.597       & 78.051       & 32.523 \\[5pt]
                                & C     & 5035.569     & 15.364         & 3.348        & 3.116       & 4.832 \\[5pt]

    \hline \xrowht[()]{10pt}
    \multirow{4}{*}{$n = 2800$} &FOMT  & 0.469         & 0.004        & 0.002        & 0.001        & 0.005 \\[5pt]
                                &A-FOMT & 2.771       & 0.030        & 0.017        & 0.008        & 0.019 \\[5pt]
                                &DS    & 9.108         & 8.271        & 7.104        & 6.869        & 7.784 \\[5pt]
                                & ABD   & 453.786       & 19.237        & 25.831       & 77.203       & 150.339 \\[5pt]
                                & C     & 7677.770     & 10.187         & 4.614        & 5.880       & 4.687 \\[5pt]

    \hline \xrowht[()]{10pt}
    \multirow{4}{*}{$n = 3200$} &FOMT  & 0.458         & 0.008        & 0.002        & 0.002        & 0.004 \\[5pt]
                                &A-FOMT& 3.265        & 0.025        & 0.017        & 0.005        & 0.025 \\[5pt]
                                &DS    & 11.713        & 10.751       & 9.173       & 9.259        & 10.254 \\[5pt]
                                & ABD   & 673.622       & 24.897        & 35.973       & 38.660       & 59.106 \\[5pt]                         & C     & 10022.263     & 52.434         & 11.785        & 7.691       & 7.461 \\[5pt]

    \hline
\bottomrule
\end{tabular}
}
\end{minipage}
\end{table}

\begin{table}[h]
\caption{The computation time (in seconds) for the critical values of the test statistics in DS \citep{dumbgen2001multiscale}, ABD \citep{akakpo2014testing} and C \citep{chetverikov2019testing} via Monte--Carlo simulations over $100$ repetitions. }
\label{T:critical}
\centering{\footnotesize
\begin{tabular}{| l| S[table-format=4.3]| S[table-format=5.3]| S[table-format=6.3]|}
\toprule
\hline \xrowht[()]{10pt}
\textbf{Sample size} & \text{DS} & \text{ABD} & \text{C}  \\[5pt]
    \hline \xrowht[()]{10pt}
    $n = 400$  & 21.309 & 145.645 & 1471.285   \\[5pt]
    \hline \xrowht[()]{10pt}
    $n = 800$  & 67.825 & 926.499 & 11579.844 \\[5pt]
    \hline \xrowht[()]{10pt}
    $n=1200$  & 150.645  & 3363.234 & 39986.723  \\[5pt]
    \hline \xrowht[()]{10pt}
    $n=1600$  & 235.340  & 7952.541 & 92650.752  \\[5pt]
    \hline \xrowht[()]{10pt}
    $n=2000$  & 361.777  & 14830.578 &  242213.872 \\[5pt]
    \hline \xrowht[()]{10pt}
    $n=2400$  & 611.162   & 24424.565 & 381613.095\\[5pt]
    \hline \xrowht[()]{10pt}
    $n=2800$  & 783.221  & 44426.914 & 655026.810   \\[5pt]
    \hline \xrowht[()]{10pt}
    $n=3200$  & 1012.866  & 68112.276 & 767809.500   \\[5pt]
    \hline 
\bottomrule
\end{tabular}
}
\end{table}

\begin{table}[ht]
\caption{The median of computational times (in seconds) of FOMT and A-FOMT for large scale datasets over $100$ and $10$  repetitions, respectively,  with signals $f_i$ for $i = 0, \dots, 4$.}
\label{T:time large scale}
\centering{\scriptsize
\begin{tabular}{|l|l|S[table-format=5.3] | S[table-format=1.3] | S[table-format=1.3]|S[table-format=1.3]|S[table-format=1.3]|}
\toprule
\hline \xrowht[()]{10pt}
\textbf{Sample size} & \textbf{Methods} & \textbf{$f_0$} & \textbf{$f_1$} & \textbf{$f_2$} & \textbf{$f_3$} & \textbf{$f_4$}\\[5pt]
    \hline \xrowht[()]{10pt}
    \multirow{2}{*}{$n = 10^5$}  & FOMT  & 5.541        & 0.047        & 0.034        & 0.015        & 0.035 \\[5pt]
                                & A-FOMT & 383.723     & 0.185        & 0.125        & 0.044        & 0.100 \\[5pt]
    \hline \xrowht[()]{10pt}
    \multirow{2}{*}{$n = 2\times10^5$}  & FOMT  & 10.620        & 0.070        & 0.052        & 0.027        & 0.056 \\[5pt]
                                & A-FOMT & 2071.370      & 0.230        & 0.137        & 0.150        & 0.272 \\[5pt]
    \hline \xrowht[()]{10pt}
    \multirow{2}{*}{$n = 3\times10^5$} & FOMT  & 14.954       & 0.106        & 0.083        & 0.045        & 0.090 \\[5pt]
                                & A-FOMT & 5105.517      & 0.225        & 0.197        & 0.368        & 0.443 \\[5pt]
    \hline \xrowht[()]{10pt}
    \multirow{2}{*}{$n = 4\times10^5$} & FOMT  & 18.541        & 0.121        & 0.090        & 0.054        & 0.092 \\[5pt]
                                & A-FOMT & 6054.229      & 0.220        & 0.285        & 0.205        & 0.402 \\[5pt]
    \hline \xrowht[()]{10pt}
    \multirow{2}{*}{$n = 5\times10^5$} &FOMT  & 23.250         & 0.156        & 0.108        & 0.063        & 0.125 \\[5pt]
                                &A-FOMT & 13547.665       & 0.244        & 0.266        & 0.341        & 0.675 \\[5pt]
    \hline \xrowht[()]{10pt}
    \multirow{2}{*}{$n = 6\times10^5$} &FOMT  & 26.307         & 0.149        & 0.136        & 0.077        & 0.164 \\[5pt]
                                &A-FOMT & 15136.516       & 0.334        & 0.268        & 0.351        & 0.634 \\[5pt]
    \hline \xrowht[()]{10pt}
    \multirow{2}{*}{$n = 7\times10^5$} &FOMT  & 30.724         & 0.178        & 0.143        & 0.085        & 0.148 \\[5pt]
                                &A-FOMT & 33094.646       & 0.322        & 0.413        & 1.009        & 1.153 \\[5pt]
    \hline \xrowht[()]{10pt}
    \multirow{2}{*}{$n = 8\times10^5$} &FOMT  & 35.518         & 0.293        & 0.162        & 0.097        & 0.201 \\[5pt]
                                &A-FOMT& 35503.989        & 0.419        & 0.378        & 0.974        & 1.244 \\[5pt]
    \hline \xrowht[()]{10pt}
    \multirow{2}{*}{$n = 9\times10^5$} &FOMT  & 39.899         & 0.241        & 0.184        & 0.107        & 0.187 \\[5pt]
                                &A-FOMT& 37707.014        & 0.348        & 0.490        & 0.763        & 1.222 \\[5pt]
    \hline \xrowht[()]{10pt}
    \multirow{2}{*}{$n = 10^6$} &FOMT  & 44.858         & 0.363        & 0.198        & 0.118        & 0.205 \\[5pt]
                                &A-FOMT& 39699.549        & 0.578        & 0.373        & 0.584        & 1.068 \\[5pt]
    \hline
\bottomrule
\end{tabular}
}
\end{table}

\end{document}